\documentclass[twoside]{article}
\usepackage{amssymb, floatflt}
\input knot.mac
\renewcommand{\copyrightheading}[1]{}

\runninghead{Enumerating the Prime Alternating Links}
{Enumerating the Prime Alternating Links}

\begin{document}
\copyrightheading{}		    %{Vol.~0, No.~0 (2001) 00--00}

%\publisher{     }

\fpage{1}
\centerline{\bf Enumerating the Prime Alternating Links}
\centerline{\footnotesize Stuart Rankin and Ortho Smith\fnm{$\dagger$}}
\centerline{\footnotesize Department of Mathematics, University of Western Ontario}
\centerline{\footnotesize srankin@uwo.ca, dsmith6@uwo.ca}
\fnt{$\dagger$}{Partially supported by Wolfgang Struss, and by NSERC Grant
R3046A02.}
\vskip20pt
%% reference labels
\def\calvo{1}
\def\conway{2}
\def\FrayMendez{3}
\def\lickorish{4}
\def\smithrankinI{5}  % Part I
\def\smithrankinII{6} % Part II

\section{Abstract}
 In [\smithrankinI], four knot operators were introduced and used
 to construct all prime alternating knots of a given crossing size.
 An efficient implementation of this construction was made possible
 by the notion of the master array of an alternating knot. The master
 array and an implementation of the construction appeared in
 [\smithrankinII]. The basic scheme (as described in [\smithrankinI])
 is to apply two of the operators, $D$ and $ROTS$, to the prime alternating
 knots of minimal crossing size $n-1$, which results in a
 large set of prime alternating knots of minimal crossing size $n$, and
 then the remaining two operators, $T$ and $OTS$, are applied to these
 $n$ crossing knots to complete the production of the set of prime
 alternating knots of minimal crossing size $n$.

 In this paper, we show how to obtain all prime alternating links of
 a given minimal crossing size. More precisely, we shall establish that 
 given any two prime alternating links of minimal crossing size $n$,
 there is a finite sequence of $T$ and $OTS$ operations that 
 transforms one of the links into the other. Consequently, one may
 select any prime alternating link of minimal crossing size $n$ (which
 is then called the seed link), and repeatedly apply only the operators $T$
 and $OTS$ to obtain all prime alternating links of minimal crossing size $n$
 from the chosen seed link. The process may be standardized by specifying
 the seed link to be (in the parlance of [\smithrankinI]) the unique
 link of $n$ crossings with group number 1, the $(n,2)$ torus link.

% We observe that if the set of all prime alternating knots of minimal
% crossing size $n$ is available, then the process can start with this set,
% The $OTS$ operator will not change the number of components in a link,
% while $T$ may increase or decrease the number of components by 1. Thus
% we may start with the prime alternating knots of minimal crossing size $n$,
% that is, the prime alternating links of one component, and apply $T$ and
% $OTS$--next paper, I think.
 
\section{Introduction}
 In [\smithrankinI], four knot operators were introduced. Of the four,
 the two called $D$ and $ROTS$ were simply specific instances of the
 general splicing operation (see Calvo [\calvo] for an extensive discussion
 of the splicing operation). A form of $D$ was also used by H. de Fraysseix
 and P. Ossona de Mendez (see [\FrayMendez]) in their work to characterize
 Gauss codes, and both $T$ and $OTS$ appeared in Conway's seminal paper
 [\conway]. These operators were used in [\smithrankinI] to present
 a method for the construction of all prime alternating knots of
 a given minimal crossing size. An efficient implementation of this method
 was presented in [\smithrankinII].

 When the operators $D$ and $ROTS$ are applied to the prime
 alternating knots of minimal crossing size $n-1$, the result is
 a large set of prime alternating knots of minimal crossing size $n$ (in
 the computational work we have done, about $98\%$ of the total number
 of knots have been constructed by $D$ and $ROTS$).
 If the remaining two operators, $T$ and $OTS$, are then applied to these
 $n$ crossing knots, the remaining prime alternating knots of minimal
 crossing size $n$ are obtained.

 In this paper, we extend this work to show how all prime alternating links
 of a given minimal crossing size may be obtained. More precisely, we shall
 establish that given any two prime alternating links of minimal crossing
 size $n$, there is a finite sequence of $T$ and $OTS$ operations that will
 transform one of the links into the other. Consequently, one may select any prime
 alternating link of minimal crossing size $n$ (which is then called the
 seed link), and repeatedly apply the operators $T$ and $OTS$ to obtain all
 prime alternating links of minimal crossing size $n$ from the chosen seed
 link. The process may be standardized by specifying the seed link to be
 the $(n,2)$ torus link.

\begin{definition}\label{link def}
 A {\it link} of $n$ components is a smooth embedding of the disjoint union
 of $n$ copies of $S^1$ into $R^3$. The image of each copy of $S^1$ is
 called a {\it component} of the link.
  A link $L$ is said to be {\it split} if
 there exist disjoint open 3-balls $U$ and $V$ such that $L\subset U\cup V$,
 $L\cap U\ne\emptyset$ and $L\cap V\ne\emptyset$.
 A link diagram is a projection of the link into a plane such
 that the preimage of any point is of size at most two, and a point has
 preimage of size two only if the point is the image of a crossing, which
 also displays the over/under behaviour of both strands of each crossing.
 A crossing is said to be a {\it link crossing} if the two strands of
 the crossing belong to different components, otherwise it is said to
 be a {\it component crossing}.
\end{definition}

 We shall regard a link diagram as a 4-regular plane graph by 
 considering each crossing as a vertex of the graph and the portion of
 the curve between two consecutive crossings as an edge between the two
 vertices. If a link is split, then there exists a diagram of
 the link which is not connected (as a graph). From now on, by link we
 shall mean non-split link.

 By a link traversal, we mean the process of assigning an orientation to
 each component, then traversing each component in the direction of its
 orientation. A link is said to be {\it alternating} if there exists a
 diagram of the link such that as the link is traversed, the sequence of
 overpasses and underpasses alternates. Such a diagram is said to be an
 {\it alternating diagram} for the link. An alternating diagram of a link
 is said to be {\it reduced} if it is loop-free.

 \begin{definition}\label{tangle def}
   For any positive integer $m$, an {\it $m$-tangle} is a connected
   plane graph in which there exists exactly one face, called the
   {\it edge face} of the $m$-tangle, for which some number $k$ of its
   boundary vertices have degree less than four, with the sum of the degrees
   of the boundary vertices equal to $4k-m$, while all other vertices
   have degree 4.
 \end{definition}

  We shall consider each vertex $v$ of degree less than four in an
  $m$-tangle to have $4-\text{deg}(v)$ arcs lying in the edge face and
  incident to $v$ (which we then refer to as arcs incident to the tangle),
  such that except
  for the endpoint $v$, the arcs lie in the interior of the edge face, and
  no two incident arcs meet each other other than at a vertex in the boundary
  of the edge face in the case of two such curves incident to the same
  boundary vertex.  With this convention, we have recovered the plane
  projection equivalent of the usual notion of a tangle.

\begin{definition}\label{tangle in link}
  Let $G$ be a 4-regular plane graph. An induced subgraph $T$ of $G$ which
  is an $m$-tangle for some positive integer $m$ is called an $m$-tangle of
  $G$.
\end{definition}

  A reduced diagram of an alternating link is a 4-regular plane graph, and
  our interest will be the $m$-tangles found in such a graph. For each
  vertex $v$ of degree less than four in a given $m$-tangle in a reduced
  diagram of an alternating link, we shall regard each arc that is incident
  to $v$ in the diagram but is not an edge of the $m$-tangle
  as an initial segment of an edge curve incident to $v$.

  We remark that since a tangle is a graph, it must contain at least one
  vertex. As well, the requirement of connectedness removes the unwanted
  situation of what were called pass-through arcs in [\smithrankinI].
 
  With this terminology, an alternating link $L$ is prime if and only if
  a reduced alternating diagram $D$ of $L$ has no 2-tangles (see
  Theorem 4.4 of [\lickorish]); equivalently, if and only if $D$ is
  3-edge-connected (that is, the removal of at most two edges does not
  disconnect the graph).

\begin{lemma}\label{incident even}
 If $T$ is an $m$-tangle in a reduced alternating diagram of an alternating
 link, then $m$ is even.
\end{lemma}

\begin{proof}%
 Suppose that $T$ contains $k$ crossings. Then the sum of the vertex
 degrees in $T$ is $4k-m$, and by the handshake lemma, this sum is twice
 the number of edges of $T$. Thus $m$ is even.
\end{proof} 

 Note that in a reduced alternating diagram of a prime alternating
 link, there are no $2$-tangles, so $m=4$ is the smallest positive integer
 for which there exists an $m$-tangle in such a diagram. We shall be
 primarily interested in $4$-tangles and $6$-tangles.

 We are now in a position to introduce the general tangle turn operation,
 of which $T$ and $OTS$ are just particular instances.
 
\begin{definition}\label{tangle turning}
 Let $D$ be a reduced alternating diagram of an alternating prime
 link, and let $T$ be an $m$-tangle of $L$. Choose an edge incident
 to $T$ and, starting with the selected edge, proceed in a clockwise
 direction around the tangle, labelling the incident edges as $1,2,\ldots,m$
 in order as they are encountered. Then for each edge incident to $T$,
 move further away from $T$ on the edge, and if the edge was labelled $i$,
 place the label $i^\prime$ on the edge, and cut the edge between the labels
 $i$ and $i^\prime$. For each $i$ from $1$ to $m-1$, attach the cut edge
 labelled $i$ to the cut edge labelled $(i+1)^\prime$, and attach the edge
 labelled $m$ to the cut edge labelled $1$. Finally, perform the unknotting
 surgery on each crossing of $T$. The result is a diagram of a link
 that is said to have been obtained by {\it turning} $T$.
\end{definition}

 It is always the case that the result of turning an $m$-tangle in a
 alternating diagram of an alternating link is again an alternating diagram
 of some alternating link. This may be seen by considering any two
 consecutive edges incident to the tangle
 being turned. They must be edges in the boundary walk of a face of the
 plane graph that is the reduced alternating diagram of the alternating link
 in question. Start at one of the edges and follow the boundary walk into
 the tangle being turned until we reach the other incident edge, labelling
 the ends of each edge traversed with either $u$ if the end of the edge is
 an underpass, or $o$ if the end of the edge is an overpass. It should now
 be apparent that the two incident edges will have opposite labelling at the
 vertices of the tangle to which they are attached. The result of turning
 the tangle will then cause every crossing of the turned tangle to have the
 wrong over/under behaviour. The unknotting surgery remedies this situation.

\begin{figure}[ht]
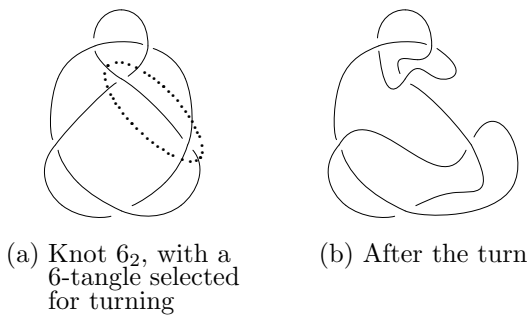

\centering
\begin{tabular}{c@{\hskip30pt}c}
$\vcenter{\xy /r10pt/:,
 (0,0)="b",
 (2.5,2.5)="b2",
 (-2.5,2.5)="b1",
 (0,5)="m",
 (1,6.25)="t2",
 (-1,6.25)="t1",
 "b";"b2"+(.2,-.2)**\crv{"b"+(.75,-.3) & "b"+(2.5,-.25) & "b2"+(.75,-.85)},
 "b2"+(-.2,.25);"m"**\crv{"b2"+(-.8,1.13) & "m"+(.6,-.35)}, 
 "m";"t1"+(0,-.2)**\crv{"m"+(-.3,.14) & "t1"+(.05,-.5)},
 "t1"+(-.05,.3);"t2"**\crv{"t1"+(0,1) & (0,8) & "t2"+(0.01,.95)},
 "t2";"m"+(.25,.1)**\crv{"t2"+(0,-.25) & "t2"+(-.3,-.8) & "m"+(.5,.3)},
 "m"+(-.2,-.15);"b1"**\crv{"m"+(-.5,-.3) & "b1"+(.5,.7)},
 "b1";"b"+(-.4,-.2)**\crv{"b1"+(-.75,-.85) & "b"+(-2.5,-.25) & "b"+(-.75,-.3)},
 "b"+(.4,.2);"b2"**\crv{"b"+(.65,.3) & "b"+(1.8,1) & "b2"+(-.15,-.5)},
 "b2";"t2"+(.2,-.1)**\crv{"b2"+(.25,.5) & "b2"+(.2,3) & "t2"+(.5,-.2)},
 "t2"+(-.2,.05);"t1"**\crv{"t2"+(-.5,.2) & "t1"+(.5,.2)},
 "t1";"b1"+(-.1,.2)**\crv{"t1"+(-.5,-.2) & "b1"+(-.2,3) & "b1"+(-.25,.5)},
 "b1"+(.1,-.2);"b"**\crv{"b1"+(.15,-.5) & "b"+(-1.8,1) & "b"+(-.65,.3)},
 (1,4)="x",
 "x"+(-.6,-.7);"x"+(.6,.7)**\crv{~*{.}"x"+(-.6,-.7)+(1,-1.1) & "x"+(3,-3) 
        & "x"+(.6,.7)+(1,-.9)},
 "x"+(-.6,-.7);"x"+(.6,.7)**\crv{~*{.}"x"+(-.6,-.7)+(-1,1.1) & "x"+(-2,2) 
        & "x"+(.6,.7)+(-1,.9)},
 \endxy}$
&
$\vcenter{\xy /r10pt/:,
 (0,0)="b",
 (2.5,2.5)="b2",
 (-2.5,2.5)="b1",
 (0,5)="m",
 (1,6.25)="t2",
 (-1,6.25)="t1",
 "b";"b2"+(.1,.3)**\crv{"b"+(.75,-.3) & "b"+(6,-.25) & "b2"+(1,2)},
 "b2";"m"+(.2,-.2)**\crv{"b2"+(-.8,1.13) & "m"+(.6,-.35)}, 
 "m";"t1"+(0,-.2)**\crv{"m"+(-.8,-1.2) & "t1"+(.05,-.5)},
 "t1"+(-.05,.3);"t2"**\crv{"t1"+(0,1) & (0,8) & "t2"+(0.01,.95)},
 "t2";"m"+(-.2,.2)**\crv{"t2"+(0,-.25) & "t2"+(-.3,-.8) & "m"+(-.5,1)},
 "b1";"b2"+(-.2,-.2)**\crv{"b1"+(1.2,2.1) & "b2"+(-1.2,-2.4)},
 "b1";"b"+(-.4,-.2)**\crv{"b1"+(-.75,-.85) & "b"+(-2.5,-.25) & "b"+(-.75,-.3)},
 "b"+(.4,.2);"b2"**\crv{"b"+(.65,.3) & "b"+(1.8,1) & "b2"+(1,-2)},
 "m";"t2"+(.2,-.1)**\crv{"m"+(.5,.8) & "m"+(1,-.5)  & "t2"+(1.5,-1)},
 "t2"+(-.2,.05);"t1"**\crv{"t2"+(-.5,.2) & "t1"+(.5,.2)},
 "t1";"b1"+(-.1,.2)**\crv{"t1"+(-.5,-.2) & "b1"+(-.2,3) & "b1"+(-.25,.5)},
 "b1"+(.1,-.2);"b"**\crv{"b1"+(.15,-.5) & "b"+(-1.8,1) & "b"+(-.65,.3)},
\endxy}$\\
\noalign{\vskip8pt}
(a) \vtop{\leftskip=0pt\hsize=.2\hsize\noindent Knot $6_2$, with a 
        6-tangle selected for turning} & (b) After the turn
\end{tabular}
\caption{Turning an $m$-tangle need not result in a
               prime link.}\label{turning in a prime link}
\end{figure}

 While it is natural to enquire as to whether turning a tangle in a
 prime alternating link will always result in a prime link, it is not
 difficult to see that this is not necessarily so. For example,
 in Figure \ref{turning in a prime link}, we turn a 6-tangle in a prime
 alternating knot and obtain an alternating knot which is the sum of
 two trefoils.

\begin{proposition}\label{prime tangle turning}
 Let $D$ be a reduced alternating diagram of an alternating prime
 link $L$, let $T$ be a 2-edge-connected $m$-tangle of $D$
 with $m\le 6$, and let $D^\prime$ be the alternating diagram of an
 alternating link $L^\prime$ that is obtained by turning $T$. Then
 $L^\prime$ is a prime alternating link and $D^\prime$ is reduced.
\end{proposition}

\begin{proof}
 Observe that $T$ is a tangle in $D^\prime$, and that $D^\prime-T=D-T$.
 Suppose that $D^\prime$ is not prime, and that $T_1$ is a 2-tangle in
 $D^\prime$. By definition of tangle, there are vertices that
 don't belong to $T$, whence
 $T_1$ must meet both $T$ and $D^\prime-T=D-T$. Furthermore, since
 $T$ is connected, there must be at least one edge of $T$ joining a
 vertex of $T_1\cap T$ to a vertex of $T-T_1$. Since any such edge is
 incident to $T_1$, there can be at most two such edges. If there were
 exactly one such edge, then that edge would be a cut-edge in $T$, which
 by hypothesis is not possible. Thus both edges incident to $T_1$ join
 vertices in $T_1\cap T$ to vertices in $T-T_1$. Now $T_1\cap
 (D^\prime-T)=T_1\cap (D-T)$ must have at least 4 incident edges
 since $D$ is prime, and since no edge incident to $T_1\cap(D^\prime-T)$
 can have an endpoint outside of $T_1$, there must be at least
 four edges in $T_1$ joining vertices in $T_1\cap (D^\prime-T)$ to vertices
 in $T_1\cap T$. If all edges incident to $T$ are found among these (that
 is to say, if all edges (four or six, as the case may be) incident to $T$
 actually join vertices in $T_1\cap(D^\prime-T)$ to vertices in
 $T_1\cap T$), then $T-T_1$ is a 2-tangle contained in $T$, which is not
 possible since $L$ is prime. Thus $T$ must be a 6-tangle, and there must
 be two edges incident to $T$ that join vertices in $T-T_1$ to vertices in
 $D^\prime-(T\cup T_1)=D-(T\cup T_1)$. But then $D-(T\cup T_1)$
 is a 2-tangle in $D$, which is not possible. Thus each case results
 in a contradiction, whence $D^\prime$ must be prime.
\end{proof}
 
\par
 \section{The $T$ and $OTS$ Operators}
 Both of these operators act on tangles; 4-tangles in the case of $T$, and
 6-tangles for $OTS$. An application of either operator to an $n$-crossing
 reduced alternating diagram results in an $n$-crossing reduced alternating
 diagram. Furthermore, in every application of $OTS$ and most applications
 of $T$, the number of link components in the resulting link will be
 unchanged from that of the original link. Certain applications of $T$ will
 result in either an increase or a decrease of 1 in the number of link
 components. As we shall see, this is enough to allow us to obtain all
 prime alternating links of minimal crossing size $n$.
 
\subsection{The $OTS$ operator}
\begin{figure}[ht]
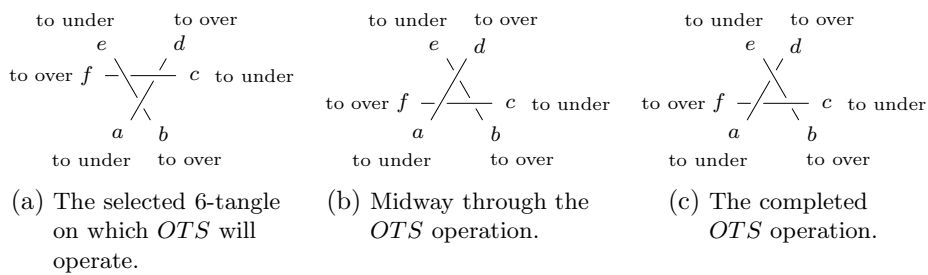

\centering
\begin{tabular}{ccc}
$\vcenter{\xy /l7pt/:,
{\xypolygon3"a"{~>{}~<{}}},
"a1";"a2"**\dir{}?(-.6)="x1"?(-.2)="x2"?(.2)="x3"
    ?(.8)="x4"?(1.2)="x5"?(1.6)="x6",
    "x1";"x4"**\dir{-},
    "x5";"x6"**\dir{-},
"a2";"a3"**\dir{}?(-.6)="y1"?(-.2)="y2"?(.2)="y3"
    ?(.8)="y4"?(1.2)="y5"?(1.6)="y6",
    "y1";"y4"**\dir{-},
    "y5";"y6"**\dir{-},
"a3";"a1"**\dir{}?(-.6)="z1"?(-.2)="z2"?(.2)="z3"
    ?(.8)="z4"?(1.2)="z5"?(1.6)="z6",
    "z1";"z4"**\dir{-},
    "z5";"z6"**\dir{-},
"x1"*!<5pt,5pt>{\hbox{\footnotesize $a$}}*!<15pt,15pt>{\hbox{\scriptsize to under}},
"z6"*!<-5pt,5pt>{\hbox{\footnotesize $b$}}*!<-15pt,15pt>{\hbox{\scriptsize to over}},
"y1"*!<-7pt,0pt>{\hbox{\footnotesize $c$}}*!CL(2){\hbox{\scriptsize to under}},
"x6"*!<-5pt,-6pt>{\hbox{\footnotesize $d$}}*!<-15pt,-15pt>{\hbox{\scriptsize to over}},
"z1"*!<5pt,-5pt>{\hbox{\footnotesize $e$}}*!<15pt,-15pt>{\hbox{\scriptsize to under}},
"y6"*!<7pt,0pt>{\hbox{\footnotesize $f$}}*!CR(2){\hbox{\scriptsize to over}},        
\endxy}$
&
$\vcenter{\xy /r7pt/:,
{\xypolygon3"a"{~>{}~<{}}},
"a1";"a2"**\dir{}?(-.6)="x1"?(-.2)="x2"?(.2)="x3"
    ?(.8)="x4"?(1.2)="x5"?(1.6)="x6",
    "x1";"x6"**\dir{-},
"a2";"a3"**\dir{}?(-.6)="y1"?(-.2)="y2"?(.2)="y3"
    ?(.8)="y4"?(1.2)="y5"?(1.6)="y6",
    "y1";"y2"**\dir{-},
    "y3";"y6"**\dir{-},
"a3";"a1"**\dir{}?(-.6)="z1"?(-.2)="z2"?(.2)="z3"
    ?(.8)="z4"?(1.2)="z5"?(1.6)="z6",
    "z1";"z2"**\dir{-},
    "z3";"z4"**\dir{-},    
    "z5";"z6"**\dir{-},
"x1"*!<-5pt,-5pt>{\hbox{\footnotesize $d$}}*!<-15pt,-15pt>{\hbox{\scriptsize to over}},
"z6"*!<5pt,-5pt>{\hbox{\footnotesize $e$}}*!<15pt,-15pt>{\hbox{\scriptsize to under}},
"y1"*!<7pt,0pt>{\hbox{\footnotesize $f$}}*!CR(2){\hbox{\scriptsize to over}},
"x6"*!<5pt,5pt>{\hbox{\footnotesize $a$}}*!<15pt,15pt>{\hbox{\scriptsize to under}},
"z1"*!<-5pt,5pt>{\hbox{\footnotesize $b$}}*!<-15pt,15pt>{\hbox{\scriptsize to over}},
"y6"*!<-7pt,0pt>{\hbox{\footnotesize $c$}}*!CL(2){\hbox{\scriptsize to under}},            
\endxy}$
&
$\vcenter{\xy /r7pt/:,
{\xypolygon3"a"{~>{}~<{}}},
"a1";"a2"**\dir{}?(-.6)="x1"?(-.2)="x2"?(.2)="x3"
    ?(.8)="x4"?(1.2)="x5"?(1.6)="x6",
    "x1";"x2"**\dir{-},
    "x3";"x6"**\dir{-},    
"a2";"a3"**\dir{}?(-.6)="y1"?(-.2)="y2"?(.2)="y3"
    ?(.8)="y4"?(1.2)="y5"?(1.6)="y6",
    "y1";"y2"**\dir{-},
    "y3";"y6"**\dir{-},
"a3";"a1"**\dir{}?(-.6)="z1"?(-.2)="z2"?(.2)="z3"
    ?(.8)="z4"?(1.2)="z5"?(1.6)="z6",
    "z1";"z2"**\dir{-},
    "z3";"z6"**\dir{-},    
"x1"*!<-5pt,-5pt>{\hbox{\footnotesize $d$}}*!<-15pt,-15pt>{\hbox{\scriptsize to over}},
"z6"*!<5pt,-5pt>{\hbox{\footnotesize $e$}}*!<15pt,-15pt>{\hbox{\scriptsize to under}},
"y1"*!<7pt,0pt>{\hbox{\footnotesize $f$}}*!CR(2){\hbox{\scriptsize to over}},
"x6"*!<5pt,5pt>{\hbox{\footnotesize $a$}}*!<15pt,15pt>{\hbox{\scriptsize to under}},
"z1"*!<-5pt,5pt>{\hbox{\footnotesize $b$}}*!<-15pt,15pt>{\hbox{\scriptsize to over}},
"y6"*!<-7pt,0pt>{\hbox{\footnotesize $c$}}*!CL(2){\hbox{\scriptsize to under}},            
\endxy}$\\
\noalign{\vskip8pt}
(a) \vtop{\hsize=1.25in\leftskip=0pt\noindent\small The selected 6-tangle on which $OTS$ will
operate.}
 & (b) \vtop{\hsize=1.25in\leftskip=0pt\noindent\small Midway through the $OTS$ operation.}
 & (c) \vtop{\hsize=1in\leftskip=0pt\noindent\small The completed $OTS$ operation.}
\end{tabular}
\caption{The $OTS$ operator}\label{ots operator}
\end{figure}

 A full, proper, alternating 6-tangle which is a cycle graph on 3 vertices
 shall be called an $OTS$ 6-tangle, and the $OTS$ operator is simply the
 general tangle turn operation applied to an $OTS$ 6-tangle. However, we
 prefer to visualize the $OTS$ operation not as a tangle turn, but rather as
 the act of moving an arc joining two crossings of the $OTS$ 6-tangle
 across the third crossing of the $OTS$ 6-tangle, much like the Reidemeister III move. This point of view is
 illustrated in Figure \ref{ots operator}. In (a), we have shown an
 $OTS$ 6-tangle, and (b) and (c) illustrate how an $OTS$ operation is
 performed on such a 6-tangle. The operator can be considered to consist of
 two stages: in the first stage, one of the three strands $ad$, $be$, or
 $cf$ is chosen. The chosen strand is cut, moved to the other side of the
 crossing formed by the other two strands, and rejoined so as to preserve
 the over/under pattern on the strand that has been moved. This stage has
 been completed in Figure \ref{ots operator} (b), where we have illustrated
 the situation if strand $cf$ was chosen.  At this point, the
 crossing formed by the other two strands is now an over-pass when it should
 be an under-pass or vice-versa, so to complete the $OTS$ operation, we must
 apply the unknotting surgery to this crossing (the crossing formed by strands
 $d$ and $e$). The completed $OTS$ operation is shown in Figure
 \ref{ots operator} (c).

The tangle turn interpretation of the $OTS$ operation makes it clear that for a given $OTS$
6-tangle, the $OTS$ operation yields the same diagram (not just flype equivalent)
independently of which of the three strands is selected to move across the remaining crossing.
A moment's reflection also reveals that $OTS$ is self-inverse.

Since an $OTS$ operation can be viewed as the cutting and rejoining of one strand,
followed by the cutting and rejoining of a second strand, it follows that the
end result is a diagram of a link with the same components as
the original link, while the tangle turn interpretation establishes that the
result of applying $OTS$ to an alternating diagram is again an alternating
diagram.

Finally, since an $OTS$ 6-tangle has no cut-edges, it follows from
Proposition \ref{prime tangle turning} that the result of applying $OTS$
to an $n$-crossing prime (hence reduced) alternating diagram is again a
prime alternating diagram.

\subsection{The $T$ operator}
 The $T$ operator works on tangles that we call (sub)groups of a link. The
 notion of group and subgroup was introduced in [\smithrankinI] for knots,
 but the same definition applies to links as well, although we do encounter
 a new scenario when we have a link of two or more components, which leads
 us to the notion of a link group. 

 \begin{definition}\label{group}
  A {\it group} in a link is a 4-tangle which is a maximal 2-braid in the
  link, and a {\it subgroup} of a group is simply a 2-braid contained in
  the group. At each end of a 2-braid, the two strands of the braid (both
  arcs incident to the tangle) are called {\it end arcs} of the subgroup. A
  group is said to be a {\it link group} if the two strands of the group
  belong to different link components, and a group that is not a link group
  is called a {\it component group}, or simply a group. A group is referred
  to as a $k$-group if it contains $k$ crossings, and a group consisting of
  a single crossing is simply referred to as a loner. A component group of
  at least two crossings is said to be {\it positive} if during a traversal
  of the component, the two strands of the group are traversed in the same
  direction, otherwise it is said to be {\it negative}.
\end{definition}  
  
\begin{figure}[ht]
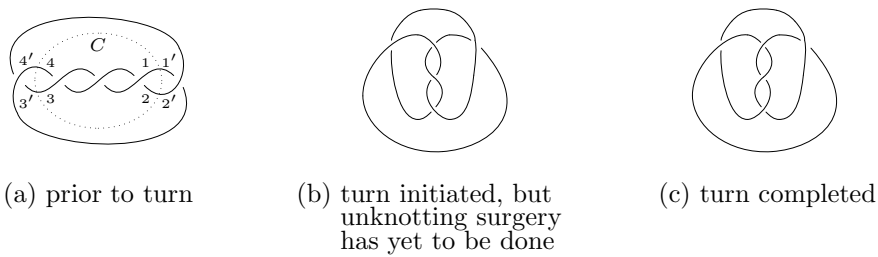

\centering
\begin{tabular}{c@{\hskip30pt}c@{\hskip30pt}c}
 $\vcenter{\xy /r30pt/:,
(1.5,.45)*{\hbox{$\ssize C$}},
(2.1,.23)*{\hbox{\tiny 1}},
(2.4,.26)*{\hbox{\tiny $1^\prime$}},
(2.1,-.23)*{\hbox{\tiny 2}},
(2.4,-.26)*{\hbox{\tiny $2^\prime$}},
(.9,.23)*{\hbox{\tiny 4}},
(.6,.28)*{\hbox{\tiny $4^\prime$}},
(.9,-.23)*{\hbox{\tiny 3}},
(.6,-.26)*{\hbox{\tiny $3^\prime$}},
(1.5,0),{\ellipse(.8,.6){.}},
(0,0)="x1",
(1,0)="x2",
"x1";"x2"**\dir{}?(.3)="x3",
%"x3"+(0,-.4)="y1";"x2"+(-.05,.065)**\crv{"x2"+(-.4,.4)},
"x3"+(.18,-.41)="y1";"x2"+(-.08,.06)**\crv{"x3"+(0,-.1) & "x2"+(-.4,.4)},
(.5,0)="x1",
(1.5,0)="x2",
"x1"+(.08,-.06);"x2"+(-.05,.065)**\crv{"x1"+(.4,-.4) & "x2"+(-.4,.4)},
(1,0)="x1",
(2,0)="x2",
"x1"+(.08,-.06);"x2"+(-.05,.065)**\crv{"x1"+(.4,-.4) & "x2"+(-.4,.4)},
(1.5,0)="x1",
(2.5,0)="x2",
"x1"+(.08,-.06);"x2"+(-.05,.065)**\crv{"x1"+(.4,-.4) & "x2"+(-.4,.4)},
(2,0)="x1",
(3,0)="x2",
"x1";"x2"**\dir{}?(.7)="x3",
"x1"+(.08,-.06);"x3"+(-.18,.41)="x4"**\crv{"x1"+(.4,-.4) & "x3"+(0,.1)},
(.5,0)="x1",
"x1"+(-.07,.085);"x4"**\crv{"x1"+(-.45,1.04) & "x4"+(-.16,.51)},
"x2"+(-.5,0)="y2",
"y1";"y2"+(.07,-.085)**\crv{"y1"+(.16,-.51) & "y2"+(.45,-1.04)},
 \endxy}$
 &
  $\vcenter{\xy /r20pt/:,
    (0,0)="o",
%% starting lower 3-group
%% lower right of 3-group
    "o"+(0,-.6)="a",
    "o";"a"+(.04,.04)**\crv{"o"+(.3,-.25)},
    "a"+(-.05,-.06)="x1";"o"+(-.3,-.8)="x3"**\crv{"x1"+(-.08,-.08) & "x3"+(.1,0)},
%% lower left of 3-group 
    "o";"a"+(-.1,.08)**\crv{~*{\null}~**\dir{}"o"+(-.3,-.3)}?(.1)="x",
    "x";"o"+(.3,-.8)**\crv{"o"+(-.3,-.4) & "o"+(.3,-.8)+(-.2,0)},
%% top left of 3-group 
    "o";"a"+(.1,-.08)**\crv{~*{\null}~**\dir{}"o"+(.3,.3)}?(.1)="x",
    "x";"o"+(-.3,.8)**\crv{"o"+(.3,.4) & "o"+(-.3,.8)+(.2,0)},    
%% top right of 3-group
    "o"+(0,.6)="a",
    "o";"a"+(-.04,-.04)**\crv{"o"+(-.3,.25)},
    "a"+(.05,.06)="x1";"o"+(.3,.8)="x3"**\crv{"x1"+(.1,.1) & "x3"+(-.05,0)},
%% left crossing
    "o"+(-.3,-.8)="x3";(-.8,.5)="b"**\crv{"x3"+(-.4,0) &  "b"+(-.03,-.3)},
    (.8,.5)+(0,.2)="d",
    "b"+(0,.2)="s";"d"+(-.02,0)**\crv{"s"+(.01,.1) & "s"+(.2,.4) & "o"+(-.3,1.3) 
         & "o"+(.3,1.3)  & "d"+(-.2,.4) & "d"+(-.01,.1)},
    "o"+(.3,-.8)="x4";"d"+(-.02,0)**\crv{"x4"+(.35,0) & "d"+(-.05,-.2)+(.03,-.3) 
        & "d"+(.015,-.15)},
    "o"+(-.3,.8);"d"+(.1,-.15)**\crv{"o"+(-.7,.85) & "o"+(-2,-.4) & "o"+(0,-2) 
         & "o"+(2,-.7)},
    "d"+(-.1,0);"o"+(.3,.8)**\crv{"d"+(-.1,0)+(-.15,.1) & "o"+(.4,.8)},
 \endxy}$
&
  $\vcenter{\xy (0,0)="o",
    {/r20pt/:,
    %% left crossing
    "o"+(-.3,-.8)="x3";(-.8,.5)="b"**\crv{"x3"+(-.4,0) &  "b"+(-.03,-.3)},
    (.8,.5)+(0,.2)="d",
    "b"+(0,.2)="s";"d"+(-.02,0)**\crv{"s"+(.01,.1) & "s"+(.2,.4) & "o"+(-.3,1.3) 
         & "o"+(.3,1.3)  & "d"+(-.2,.4) & "d"+(-.01,.1)},
    "o"+(.3,-.8)="x4";"d"+(-.02,0)**\crv{"x4"+(.35,0) & "d"+(-.05,-.2)+(.03,-.3) 
        & "d"+(.015,-.15)},
    "o"+(-.3,.8);"d"+(.1,-.15)**\crv{"o"+(-.7,.85) & "o"+(-2,-.4) & "o"+(0,-2) 
         & "o"+(2,-.7)},
    "d"+(-.1,0);"o"+(.3,.8)**\crv{"d"+(-.1,0)+(-.15,.1) & "o"+(.4,.8)}},
       {0;<-20pt,0pt>:<0pt,20pt>::,
%% starting lower 3-group
%% lower right of 3-group
    "o"+(0,-.6)="a",
    "o";"a"+(.04,.04)**\crv{"o"+(.3,-.25)},
    "a"+(-.05,-.06)="x1";"o"+(-.3,-.8)="x3"**\crv{"x1"+(-.08,-.08) & "x3"+(.1,0)},
%% lower left of 3-group 
    "o";"a"+(-.1,.08)**\crv{~*{\null}~**\dir{}"o"+(-.3,-.3)}?(.1)="x",
    "x";"o"+(.3,-.8)**\crv{"o"+(-.3,-.4) & "o"+(.3,-.8)+(-.2,0)},
%% top left of 3-group 
    "o";"a"+(.1,-.08)**\crv{~*{\null}~**\dir{}"o"+(.3,.3)}?(.1)="x",
    "x";"o"+(-.3,.8)**\crv{"o"+(.3,.4) & "o"+(-.3,.8)+(.2,0)},    
%% top right of 3-group
    "o"+(0,.6)="a",
    "o";"a"+(-.04,-.04)**\crv{"o"+(-.3,.25)},
    "a"+(.05,.06)="x1";"o"+(.3,.8)="x3"**\crv{"x1"+(.1,.1) & "x3"+(-.05,0)}},
 \endxy}$
\\
\noalign{\vskip8pt}
(a) prior to turn & (b) \vtop{\leftskip=0pt\hsize=.25\hsize\noindent turn initiated,
  but unknotting surgery has yet to be done}
   & (c) turn completed
\end{tabular}
\caption{The $T$ operator: turning a subgroup of size 3.}\label{T operator}
\end{figure}

 The $T$ operator is simply the general tangle turn operation, but performed
 only on full, proper 4-tangles that are (sub)groups of the link. Since a
 (sub)group does not contain any cut-edges, it follows from Proposition
 \ref{prime tangle turning} that the result of applying $T$ to an $n$-crossing
 reduced alternating diagram again an $n$-crossing reduced alternating
 diagram. Furthermore, it is evident that $T$ is its own inverse. We
 illustrate the $T$ operator in Figure \ref{T operator}, where in
 Figure \ref{T operator} (a), the torus knot of five crossings is shown,
 with a subgroup of size three singled out for turning. Note that the group
 of five crossings (and thus the subgroup of three crossings) is a positive
 group. Figure \ref{T operator} (b), the turn has been initiated, but the
 unknotting has yet to be done. Finally, in Figure \ref{T operator} (c), the
 completed turn is shown. The result is a prime alternating knot with a
 negative 2-group and the turned group, now a negative 3-group.

\begin{figure}[ht]
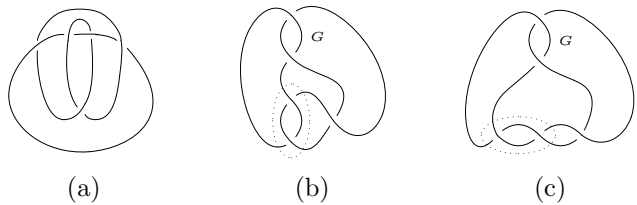

\centering
\begin{tabular}{c@{\hskip30pt}c@{\hskip30pt}c}
  $\vcenter{\xy (0,0)="o",
    /r20pt/:,
    %% left crossing
    "o"+(-.3,-.8)="x3";(-.8,.5)="b"**\crv{"x3"+(-.4,0) &  "b"+(-.03,-.3)},
    (.8,.5)+(0,.2)="d",
    "b"+(0,.2)="s";"d"+(-.02,0)**\crv{"s"+(.01,.1) & "s"+(.2,.4) & "o"+(-.3,1.3)
         & "o"+(.3,1.3)  & "d"+(-.2,.4) & "d"+(-.01,.1)},
    "o"+(.3,-.8)="x4";"d"+(-.02,0)**\crv{"x4"+(.35,0) & "d"+(-.05,-.2)+(.03,-.3) 
        & "d"+(.015,-.15)},
    "o"+(-.3,.81);"d"+(.1,-.15)**\crv{"o"+(-.7,.85) & "o"+(-2,-.4) & "o"+(0,-2) 
         & "o"+(2,-.7)},
    "d"+(-.1,0);"o"+(.3,.81)**\crv{"d"+(-.1,0)+(-.15,.1) & "o"+(.4,.8)},
    "x3";"x3"+(.4,.3)**\crv{"x3"+(.25,0)},
    "x3"+(.4,.3);"o"+(.22,.72)**\crv{"x3"+(.6,.6)},
    "o"+(.2,.88);"o"+(-.22,.7)**\crv{"o"+(.2,1) &"o"+(0,1.2) & "o"+(-.2,1)},
    "x4";"o"+(.1,-.7)**\crv{"x4"+(-.15,0)},
    "o"+(-.045,-.56);"o"+(-.22,.8)+(0,-.1)**\crv{"o"+(-.2,-.38) & "o"+(-.3,.1)},
    "o"+(.35,.81);"o"+(-.1,.81)**\crv{"o"+(.20,.83) },
 \endxy}$
 &
 $\vcenter{
 \xy /r8pt/:,
 (0,0)="c1",
 (0,2)="c2",
 (2,1)="x",
 (0,4)="d1",
 (0,6)="d2",
 "c1"*=<4pt>{\hbox to 4pt{\hfil}}="c1b",
 "c2"*=<4pt>{\hbox to 4pt{\hfil}}="c2b",
 "d1"*=<4pt>{\hbox to 4pt{\hfil}}="d1b",
 "d2"*=<4pt>{\hbox to 4pt{\hfil}}="d2b",
 "x"*=<4pt>{\hbox to 4pt{\hfil}}="xb",  
 "c1";"c2b"**\crv{(-1,1)},
 "c2b";"x"**\crv{(1,3) & (1.5,1.65)},
 "x";"d2b"**\crv{(3,-.25) & (5,1) & (4,5) & (1,7)},
 "d2b";"d1"**\crv{(-1,5)},
 "d1";"xb"**\crv{(.5,3.5) & (1.5,3.25) & (3,2.4)},
 "xb";"c1"**\crv{(1.5,.35) & (1,-1)},
 "c1b";"c2"**\crv{(1,1)},
 "c2";"d1b"**\crv{(-1,3)},
 "d1b";"d2"**\crv{(1,5)},
 "d2";"c1b"**\crv{(-.75,6.75) & (-2.5,5.5) & (-2.5,1.5) & (-1,-.75)},
  "c1";"c2"**\dir{}?(.5)="c0",
 "c1";"c0",{\ellipse<14pt,7pt>{.}},
  "d1";"d2"**\dir{}?(.5)="d0",
 "d0"*!<-10pt,0pt>{\hbox{\tiny$G$}},
\endxy}$
&
 $\vcenter{
 \xy /r8pt/:,
 (-2,.5)="c1",
 (0,.5)="c2",
 (2,.5)="x",
 (0,4)="d1",
 (0,6)="d2",
 "c1"*=<4pt>{\hbox to 4pt{\hfil}}="c1b",
 "c2"*=<4pt>{\hbox to 4pt{\hfil}}="c2b",
 "d1"*=<4pt>{\hbox to 4pt{\hfil}}="d1b",
 "d2"*=<4pt>{\hbox to 4pt{\hfil}}="d2b",
 "x"*=<4pt>{\hbox to 4pt{\hfil}}="xb",  
 "c1";"c2b"**\crv{(-1,-.5)},
 "c2b";"x"**\crv{ (1,1.5)},
 "x";"d2b"**\crv{(3,-.25) & (5,1) & (4,5) & (1,7)},
 "d2b";"d1"**\crv{(-1,5)},
 "d1";"xb"**\crv{(.5,3.5) & (1.5,3.25) & (3,2.4)},
 "xb";"c2"**\crv{(1,-.5)},
 "c1b";"c2"**\crv{(-1,1.5)},
 "c1";"d1b"**\crv{(-2.4,1.2) & (-2,2.5) & (-1,3)},
 "d1b";"d2"**\crv{(1,5)},
 "d2";"c1b"**\crv{(-.75,6.75) & (-2.5,5.5) & (-4,1.5) & (-3,-.75)},
 "c1";"c2"**\dir{}?(.5)="c0",
 "c1";"c0",{\ellipse<14pt,7pt>{.}},
 "d1";"d2"**\dir{}?(.5)="d0",
 "d0"*!<-10pt,0pt>{\hbox{\tiny$G$}},
\endxy}$ \\
\noalign{\vskip6pt}
(a) & (b) & (c)
\end{tabular}
\caption{$T$ can change the link component count.}\label{T makes link}
\end{figure}

 \noindent The $T$ operator can change the component count for a link. For
 example, when we apply $T$ to the subgroup consisting of the top two
 crossings of the negative 3-group in the knot shown in Figure \ref{T operator}
 (c), the result is the 2-component prime link shown in Figure \ref{T makes link} (a).
 
In  Figure \ref{T makes link} (b), a 2-component prime alternating link is
shown, with a link 2-group selected for turning. The effect of performing
$T$ on this link 2-group is seen in  Figure \ref{T makes link} (c). Note
that in this example, the turned 2-group has become a subgroup of a negative
3-group of the 1-component link. Note also that the link 2-group labelled by
$G$ in (b) has necessarily become a component 2-group in (c), positive as it
turns out. It is a simple matter to modify the example in (b) by adding a crossing
so that the negative group that results when the link 2-group is turned is a
4-group, but now $G$ is a negative group.

In the preceding examples, we have seen several of the situations that may result
when a group is turned. The next proposition, the proof of which is straightforward
and is omitted, describes all of the situations that may occur when a group is turned.

\begin{proposition}\label{all turn scenarios}
 Let $D$ be a reduced alternating diagram of an alternating link $L$,
 and let $G$ be a (sub)group of $D$. Let $D^\prime$ denote the reduced
 alternating diagram that results when $T$ is applied to $G$, and let
 $L^\prime$ denote the prime alternating link that is represented by
 $D^\prime$. Then
 \begin{alphlist}
  \item if $G$ is an even link (sub)group, then $D^\prime$ is a prime
  alternating link with one
  more component than $L$, and $G$ is now a negative group in $D^\prime$
  (possibly a subgroup of a larger negative group);
  \item if $G$ is an odd link (sub)group, then $L^\prime$ is a prime
  alternating link with the same number of components as $L$, and $G$ is
  a link group in $D^\prime$ (possibly a subgroup of a larger link group);
  \item if $G$ is an even component (sub)group, then $L^\prime$ is a prime
  alternating link with one fewer components than $L$, and $G$ is now a
  link group in $D^\prime$ (possibly now a subgroup of a larger link group);
  \item if $G$ is an odd component (sub)group, then $L^\prime$ is a prime
  alternating link with the same number of components as $L$, and $G$ is
  still a component (sub)group, though now with the opposite sign to that
  which it had in $D$ (and again, $G$ might now be a subgroup of a larger
  group in $D^\prime$.
  \end{alphlist}
\end{proposition}

 We also remark that if the (sub)group that is to have $T$ applied to it
 is a loner, then the result is the identical diagram. Thus in
 practice, we shall never apply $T$ to a loner. Furthermore, if an
 $(n-1)$-subgroup of the $n$-group of the $(n,2)$ torus link is turned, the
 result is the mirror image of the $(n,2)$ torus link.
 
\par
 \section{$OTS$ and $T$ make all prime alternating links}
 
  We are now in a position to establish our main objective in this paper; namely,
  to establish that given any two prime alternating links of the same minimal
  crossing size, there is a finite sequence of $T$ and $OTS$ operations that will
  transform one into the other. Since $T$ and $OTS$ are each self-inverse, it
  will suffice to prove that for any prime alternating link $L$ of a given minimal
  crossing size $n$, there is a finite sequence of $T$ and $OTS$ operations that
  will tranform the $(n,2)$ torus link into $L$. This work relies on Theorem 5 of
  [\smithrankinI], a graph-theoretic result which is applicable to reduced
  alternating diagrams of prime alternating links as well as to those of prime
  alternating knots, as was the case in [\smithrankinI].

  For the reader's convenience, we present below the main graph-theoretic notions
  and results from [\smithrankinI] that will be required for the subsequent work.

\begin{definition}
 A plane graph $G$ whose edges are piecewise smooth curves is called a 2-region,
 respectively a minimal loop, if the following conditions are satisfied:
 \begin{alphlist}
  \item $G$ is connected;
  \item $G$ has a face $F$ whose boundary is a cycle $C$;
  \item there are exactly two vertices on $C$ that have degree 2 in $G$
  (called the base vertices of the 2-region), respectively there is exactly
  one vertex on $C$ that has degree 2 in $G$ (called the base vertex of the
  minimal loop);
  \item every non-base vertex on $C$ has degree 3 in $G$;
  \item each vertex that lies in the interior of the region $\mathbb{R}^2-F$
   has degree 4 in $G$.
 \end{alphlist}

 \noindent
  The cycle $C$ is called the {\em boundary} of the 2-region or minimal loop,
  respectively, while the interior of the region $\mathbb{R}^2-F$ is
  called the {\em interior} of the 2-region, respectively minimal loop (note
  that in the event that all vertices of $G$ lie on $C$, then $F$ is not
  uniquely determined--in such a case, let $F$ be the bounded region).
  In the case of a 2-region with boundary cycle $C$ and base vertices $p$ and
  $q$, the two paths between $p$ to $q$ that $C$ determines are called the
  boundary paths of $G$.
  
  More generally, if $G$ is a plane graph with vertex set $V$, then
  a subgraph $G'$ of $G$ with vertex set $V'\subseteq V$ is said to be
  a 2-region (respectively, minimal loop) of $G$ if $G'$ is
  a 2-region (respectively, minimal loop) and every vertex of $G$ that lies in
  the interior of $G'$ belongs to $V'$, and every edge of $G$ that meets the
  interior of $G'$ belongs to $V'$. If $G'$ is a 2-region or minimal loop
  of $G$, then $G$ is said to contain the 2-region or minimal loop $G'$.

  A 2-region $G$ is said to be {\it minimal} if $H$ is a 2-region of $G$ implies
  that $H=G$, and a 2-region that has no vertices in its interior is called a
  2-{\it group}. A minimal loop with a single vertex is said to be trivial.
\end{definition}

 Note that the interior of a 2-region or minimal loop could be the unbounded
 region determined by the boundary of the 2-region, respectively minimal loop.
 Further note that a 2-group has exactly two vertices, which are multiply
 connected by two edges.

 \begin{proposition}[Proposition 8, {[\smithrankinI}]\label{existence of 2-region}
  Every non-trivial minimal loop contains a 2-region, and every 2-region contains
  a minimal 2-region.
 \end{proposition}

\begin{definition}
 Let $D$ be a reduced alternating diagram of a prime alternating link $L$, and
 let $v$ be a component crossing of $D$. Choose any of the four edges incident
 to $v$ and construct a closed walk based at $v$ by following a link
 traversal in the direction of the chosen arc until $v$ is reached for
 the first time. Such a closed walk is called a {\it component circuit based
 at $v$}.
\end{definition}

The proof of the next result as presented in [\smithrankinI] for prime alternating
knots is actually valid for prime alternating links, with the modifications to the
statement as shown below.

\begin{proposition}[Proposition 9, {[\smithrankinI}]\label{existence of min loops in cl}
 Let $D$ be a reduced alternating diagram of a prime alternating link $L$. 
 Then
 \begin{alphlist}
  \item For every component crossing $v$ of $D$, and each component circuit 
   $C$ based at $v$, there is a minimal loop of $D$ whose boundary
   cycle is a subwalk of $C$;
  \item every minimal loop of $D$ contains a 2-region of $D$, and
  \item every 2-region of $D$ contains a minimal 2-region of $D$.
 \end{alphlist} 
\end{proposition}

\begin{corollary}
 Every reduced alternating diagram of a prime alternating link different from the
 unknot has a minimal 2-region.
\end{corollary}

\begin{proof}
  To begin with, recall that by link we mean proper link, so our link is
  not the unknot, nor does it consist of two or more unlinked unknots.
  If it contains a component crossing, then by Proposition \ref{existence of
  min loops in cl}, it contains a minimal 2-region. Suppose then that our
  diagram contains only link crossings. Choose any crossing, and select
  two adjacent edges incident to the crossing. Follow these edges out from the
  crossing. The two strands belong to different components, and since the
  configuration has no component crossings, each component forms a simple
  closed curve in the plane. Thus the two strands must meet again, and we
  proceed until a point of intersection of the two strands is encountered.
  The two paths that we have followed together form a simple closed curve in
  the plane, and the two paths, together with all vertices and edges in
  one of the regions determined by the closed curve forms a 2-region. By
  Proposition \ref{existence of min loops in cl} (c), the configuration
  contains a minimal 2-region.
\end{proof}

Note that an empty 2-region is simply a 2-subgroup.

Our next observation is that the $OTS$ operation on an alternating link diagram
has an analog for 4-regular plane graphs.

\begin{definition}\label{ots in 4-reg graph}
 Let $G$ be a 4-regular plane graph. If $C$ is a 3-cycle in $G$ such that
 no two vertices of $C$ are multiply-connected, then $C$ is called an
 $ots$-triangle in $G$. The 4-regular plane graph $G^\prime$ that results from $G$
 upon modifying an $ots$-triangle $C$ as shown in Figure \ref{ots operations
 on 3-cycle} is said to have been obtained from $G$ by applying $ots$ to $C$.
\end{definition}

\begin{figure}[ht]
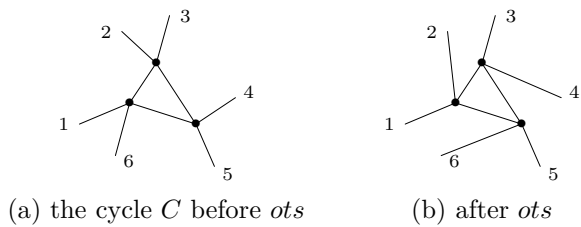

 \centering
   \begin{tabular}{c@{\hskip30pt}c}
    $\vcenter{\xy /r20pt/:,
       (0,0)="a"*=0{\hbox{$\ssize\bullet$}},
       (1.25,-.4)="b"*=0{\hbox{$\ssize\bullet$}},
       (.5,.75)="c"*=0{\hbox{$\ssize\bullet$}},
       "a";"b"**\dir{-};"c"**\dir{-};"a"**\dir{-},
       "a";"a"+(-.95,-.4)="ax"**\dir{-},
       "a";"a"+(-.27,-1)="ay"**\dir{-},
       "ax"*!<6pt,0pt>{\hbox{$\ssize 1$}},
       "ay"*!<-5pt,2pt>{\hbox{$\ssize 6$}},
       "c";"c"+(-.65,.6)="cx"**\dir{-},
       "c";"c"+(.25,.9)="cy"**\dir{-},
       "cx"*!<6pt,0pt>{\hbox{$\ssize 2$}},
       "cy"*!<-6pt,0pt>{\hbox{$\ssize 3$}},
       "b";"b"+(.75,.5)="bx"**\dir{-},
       "b";"b"+(.35,-.8)="by"**\dir{-},
       "bx"*!<-5pt,-2pt>{\hbox{$\ssize 4$}},
       "by"*!<-5pt,3pt>{\hbox{$\ssize 5$}},       
     \endxy}$
   &
    $\vcenter{\xy /r20pt/:,
       (0,0)="a"*=0{\hbox{$\ssize\bullet$}},
       (1.25,-.4)="b"*=0{\hbox{$\ssize\bullet$}},
       (.5,.75)="c"*=0{\hbox{$\ssize\bullet$}},
       "a";"b"**\dir{-};"c"**\dir{-};"a"**\dir{-},
       "a";"a"+(-.95,-.4)="ax"**\dir{-},
       "b";"a"+(-.27,-1)="ay"**\dir{-},
       "ax"*!<6pt,0pt>{\hbox{$\ssize 1$}},
       "ay"*!<-5pt,2pt>{\hbox{$\ssize 6$}},
       "a";"c"+(-.65,.6)="cx"**\dir{-},
       "c";"c"+(.25,.9)="cy"**\dir{-},
       "cx"*!<6pt,0pt>{\hbox{$\ssize 2$}},
       "cy"*!<-6pt,0pt>{\hbox{$\ssize 3$}},
       "c";"b"+(.75,.5)="bx"**\dir{-},
       "b";"b"+(.35,-.8)="by"**\dir{-},
       "bx"*!<-5pt,-2pt>{\hbox{$\ssize 4$}},
       "by"*!<-5pt,3pt>{\hbox{$\ssize 5$}},       
     \endxy}$ \\
     \noalign{\vskip 6pt}
     (a) the cycle $C$ before $ots$ & (b) after $ots$
   \end{tabular}
  \caption{The $ots$ operation on a 4-regular plane graph}
  \label{ots operations on 3-cycle}
\end{figure}

Furthermore, since a 4-regular, 3-edge-connected plane graph can be considered
as an reduced alternating diagram of a prime alternating link, it follows from
Proposition \ref{prime tangle turning} that the result of an $ots$ operation on a
4-regular, 3-edge-connected plane graph is again a 4-regular, 3-edge-connected
plane graph.

We shall consider $ots$ operations on 3-cycles in 2-regions and minimal loops,
and for this purpose, it is convenient to introduce specialized versions of $ots$
for 2-regions and minimal loops. These amount to the restrictions of $ots$
to the various situations involving a 3-cycle in a 2-region or a minimal loop.

\begin{definition}\label{ots in minimal region}
 Let $G$ be a 2-region or a minimal loop. An $ots$-triangle in $G$ is a
 face of degree 3 none of whose boundary edges belong to the boundary of a
 2-group of $G$.

 Let $O$ be an $ots$-triangle in $G$, say with boundary edges $B_1$, $B_2$
 and $B_3$. There are three possible situations: none, exactly one, or
 exactly two of the edges $B_1$, $B_2$, $B_3$ is a boundary edge of $G$.
 We define the $ots$ operation in each of these three cases. Let
 $V$ denote the vertex set of $G$ and $E$ denote the edge set of $G$. 
 
 \noindent Case 1: none of $B_1$, $B_2$, $B_3$ is a boundary edge of $C$.
 Now each pair of these edges have exactly one endpoint in common.
 Since a boundary vertex has at most one edge lying in the interior of $G$
 incident to it, we see that none of these three common endpoints is a
 boundary vertex of $G$, and so the compact set $B_1\cup B_2\cup B_3$ is
 contained in the interior of $G$. Furthermore, each of the three common
 endpoints has degree 4, so at each there are two additional incident edges.
 Since the edge curves are smooth, there is an open neighborhood $U$ of $O$
 whose intersection with each of these six additional edges is connected,
 and which does not meet any other edge curve of $G$. Arbitrarily choose
 one pair of boundary edges of $O$, and suppose that the edges were labelled
 so that $B_1$ and $B_2$ are the chosen edges. Let $a$ denote the common
 endpoint of $B_1$ and $B_2$, and let $e$ and $f$ denote the two edges
 that are incident to $a$ in addition to $B_1$ and $B_2$. Further suppose
 that all labelling has been done so that in a clockwise scan at $a$, the edges
 are encountered in the order $B_1$, $e$, $f$ and $B_2$. Choose a point $x\ne a$
 in $e \cap U$, and choose a point $y\ne a$ in $f\cap U$. Let $B_3'$ be a
 smooth curve from $x$ to $y$ within $U$ that does not meet any curve of $G$
 other than $e$ at $x$ and $f$ at $y$. Let $b$ denote the common endpoint of
 $B_2$ and $B_3$, and let $c$ denote the common endpoint of $B_3$ and $B_1$.
 Further, let $g$ and $h$ denote the two edges incident to $b$ other than $B_2$
 and $B_3$, labelled in order in the clockwise direction, and let $i$ and $j$
 denote the two edges incident to $c$ other than $B_3$ and $B_1$, labelled in
 order in the clockwise direction. Let $j'$ denote a smooth curve with endpoint
 $x$ that agrees with $j$ outside of $U$ and within $U$ meets no edge curve of
 $G$ except for $e$ at $x$, and let $g'$ denote a smooth curve with endpoint
 $y$ that agrees with $g$ outside of $U$ and within $U$ meets no edge curve of
 $G$ except for $f$ at $y$. Let the portion of $e$ from $a$ to $x$ be denoted
 by $B_1'$ and denote the remaining portion of $e$ by $e'$. Similarly, let
 the portion of $f$ from $a$ to $y$ be denoted by $B_2'$ and denote the remaining
 portion of $f$ by $f'$. Finally, let $i'=i\cup B_1$ and $h'=h\cup B_2$. Let
 $G'$ denote the plane graph with piecewise smooth edge curves whose vertex
 set is $\left(V-\{\,b,c\,\}\right)\cup\{\,x,y\,\}$ and edge set $\left(E-
 \{\,B_1,B_2,B_3,e,f,g,h,i,j\,\}\right)\cup \{B_1',B_2',B_3',e',f',g',h',i',j'\,\}$.
 Since $G'$ agrees with $G$ outside of $U$, it follows that $G'$ is a 2-region,
 respectively minimal loop, said to be obtained from $G$ by an $ots$ operation
 (on $O$).
 
 \noindent Case 2: exactly one of $B_1$, $B_2$, $B_3$ is a boundary edge of $G$.
 Suppose that the curves were labelled so that $B_3$ is the boundary curve of $G$.
 Label the vertices and edges as in Case 1, with the only differences stemming
 from the fact that the endpoints of $B_3$ are boundary vertices of $G$ (non-base,
 since neither $B_1$ nor $B_2$ is a boundary edge of $G$), so there is no
 edge $h$ or $i$. Carry out the construction as in Case 1 (so there is no
 corresponding $h'$ or $i'$), and let $G'$ denote the plane graph with piecewise
 smooth edge curves whose vertex set is $\left(V-\{\,a,b,c\,\}\right)\cup
 \{\,x,y\,\}$ and edge set $\left(E-\{\,B_1,B_2,B_3,e,f,g,i\,\}\right)\cup
 \{B_3',e',f',g',i'\,\}$. Since $G'$ agrees with $G$ outside of $U$, it follows
 that $G'$ is a 2-region, respectively minimal loop, said to be obtained from $G$
 by an $ots$ operation (on $O$). Note that the number of vertices of $G'$ is one
 less than the number of vertices of $G$. In this case, we say that vertex $a$
 has been $ots$-ed out of $G$.
 
\noindent Case 3: two of the edges $B_1$, $B_2$, $B_3$ are boundary edges of $G$.
 Suppose that the curves have been labelled so that $B_3$ is not a boundary edge
 of $G$. Let the endpoints of $B_3$ be $b$ and $c$. Then $b$ and $c$ are
 non-base boundary vertices, while the common endpoint of $B_1$ and $B_2$
 is a base boundary vertex. Each of $b$ and $c$ have one additional boundary
 edge incident to them. Let $h$ be the additional boundary edge incident to
 $b$ and let $i$ be the additional boundary edge incident to $c$, and set
 $h'=h\cup B_2$ and $i'=i\cup B_1$. Let $G'$ denote the plane graph with
 piecewise smooth edge curves whose vertex set is $V-\{\,b,c\,\}$ and edge
 set $\left(E-\{\,B_1,B_2,B_3,h,i\,\}\right)\cup \{h',i'\,\}$. Since $G'$
 agrees with $G$ outside of $U$, it follows that $G'$ is a 2-region, respectively
 minimal loop, said to be obtained from $G$ by an $ots$ operation (on $O$).

\noindent
 In Cases 2 and 3, we say that the $ots$-triangle is on the boundary of $G$.
\end{definition}

\begin{figure}[ht]
 \centering
   \begin{tabular}{c@{\hskip30pt}c}
%   \noalign{\vskip6pt}
$\vcenter{\xy /r40pt/:,
(.625,.65)="a11";"a11";{\ellipse(.9){.}},
(0,0)="a1";(1,1.5)="a2"**\dir{-},
(.25,1.5)="a3";(1.25,0)="a4"**\dir{-},
(-.25,.5)="a5";(1.5,.5)="a6"**\dir{-},
(.25,1.25)*{\hbox{$\ssize e$}},
(1,1.25)*{\hbox{$\ssize f$}},
(1.375,.625)*{\hbox{$\ssize g$}},
(1,.125)*{\hbox{$\ssize h$}},
(.25,.125)*{\hbox{$\ssize i$}},
(-.125,.625)*{\hbox{$\ssize j$}},
(.625,1.125)*{\hbox{$\sssize a$}},
(1,.625)*{\hbox{$\sssize b$}},
(.25,.625)*{\hbox{$\sssize c$}},
(.375,.86)*{\hbox{$\sssize {B_1}$}},
(.875,.86)*{\hbox{$\sssize {B_2}$}},
(.625,.375)*{\hbox{$\sssize {B_3}$}},
(.625,.95)*{\hbox{\small$\bullet$}},
(.325,.5)*{\hbox{\small$\bullet$}},
(.925,.5)*{\hbox{\small$\bullet$}},
\endxy}$
&
$\vcenter{\xy /r40pt/:,
(.625,.65)="a11";"a11";{\ellipse(.9){.}},
(0,0)="a1";(1.25,1.25)="a2"**\dir{-},
(0,1.25)="a3";(1.25,0)="a4"**\dir{-},
(-.25,.5)="a5";(.25,1)="a7"**\dir{-},
"a7";(1,1)="a8"**\dir{-},
"a8";(1.5,.5)="a6"**\dir{-},
(.25,1.25)*{\hbox{$\ssize e'$}},
(1.1,1.225)*{\hbox{$\ssize f'$}},
(1.375,.425)*{\hbox{$\ssize g'$}},
(.9,.125)*{\hbox{$\ssize h'$}},
(.25,.125)*{\hbox{$\ssize i'$}},
(-.125,.5)*{\hbox{$\ssize j'$}},
(.625,.45)*{\hbox{$\sssize a$}},
(.125,1)*{\hbox{$\sssize x$}},
(1.13,1)*{\hbox{$\sssize y$}},
(.3,.75)*{\hbox{$\sssize {B'_1}$}},
(.925,.75)*{\hbox{$\sssize {B'_2}$}},
(.625,1.15)*{\hbox{$\sssize {B'_3}$}},
(.625,.625)*{\hbox{\small$\bullet$}},
(.25,1)*{\hbox{\small$\bullet$}},
(1,1)*{\hbox{\small$\bullet$}},
\endxy}$ \\
\noalign{\vskip 6pt}
(a) An interior $ots$-triangle & (b) After the $ots$ operation \\
\noalign{\vskip 6pt}
$\vcenter{\xy /r40pt/:,
(.625,.65)="a11";"a11";{\ellipse(.9){.}},
(.325,.5)="a1";(1,1.5)="a2"**\dir{-},
(.25,1.5)="a3";(.925,.5)="a4"**\dir{-},
(-.25,.5)="a5";(1.5,.5)="a6"**\dir{-},
(.25,1.25)*{\hbox{$\ssize e$}},
(1,1.25)*{\hbox{$\ssize f$}},
(1.375,.625)*{\hbox{$\ssize g$}},
(-.125,.625)*{\hbox{$\ssize j$}},
(.625,1.125)*{\hbox{$\sssize a$}},
(1,.625)*{\hbox{$\sssize b$}},
(.25,.625)*{\hbox{$\sssize c$}},
(.375,.86)*{\hbox{$\sssize {B_1}$}},
(.875,.86)*{\hbox{$\sssize {B_2}$}},
(.625,.375)*{\hbox{$\sssize {B_3}$}},
(.625,.95)*{\hbox{\small$\bullet$}},
(.325,.5)*{\hbox{\small$\bullet$}},
(.925,.5)*{\hbox{\small$\bullet$}},
\endxy}$
&
$\vcenter{\xy /r40pt/:,
(.625,.65)="a11";"a11";{\ellipse(.9){.}},
(1,1)="a1";(1.25,1.25)="a2"**\dir{-},
(0,1.25)="a3";(.25,1)="a4"**\dir{-},
(-.25,.5)="a5";(.25,1)="a7"**\dir{-},
"a7";(1,1)="a8"**\dir{-},
"a8";(1.5,.5)="a6"**\dir{-},
(.25,1.25)*{\hbox{$\ssize e'$}},
(1.1,1.225)*{\hbox{$\ssize f'$}},
(1.375,.425)*{\hbox{$\ssize g'$}},
(-.125,.5)*{\hbox{$\ssize j'$}},
(.125,1)*{\hbox{$\sssize x$}},
(1.13,1)*{\hbox{$\sssize y$}},
(.625,.875)*{\hbox{$\sssize {B'_3}$}},
(.25,1)*{\hbox{\small$\bullet$}},
(1,1)*{\hbox{\small$\bullet$}},
\endxy}$\\
\noalign{\vskip 6pt}
(c) An $ots$-triangle with one boundary edge & (d) After the $ots$ operation \\
\noalign{\vskip 6pt}
$\vcenter{\xy /r40pt/:,
(.625,.75)="a11";"a11";{\ellipse(.9){.}},
(-.25,.75)="a1";(1.2,1.45)="a2"**\dir{-},
"a1";(1.2,.05)="a3"**\dir{-},
(.75,1.225)="a4";(.75,.275)="a5"**\dir{-},
(-.25,.75)*{\hbox{\small$\bullet$}},
(.75,1.225)*{\hbox{\small$\bullet$}},
(.75,.275)*{\hbox{\small$\bullet$}},
(-.375,.75)*{\hbox{$\ssize a$}},
(.85,1.15)*{\hbox{$\sssize b$}},
(.85,.35)*{\hbox{$\sssize c$}},
(.9,1.45)*{\hbox{$\ssize h$}},
(.9,.05)*{\hbox{$\ssize i$}},
(.25,.375)*{\hbox{$\sssize {B_1}$}},
(.25,1.13)*{\hbox{$\sssize {B_2}$}},
(.88,.75)*{\hbox{$\sssize {B_3}$}},
\endxy}$
&
$\vcenter{\xy /r40pt/:,
(.625,.75)="a11";"a11";{\ellipse(.9){.}},
(-.25,.75)="a1";(1.2,1.45)="a2"**\dir{-},
"a1";(1.2,.05)="a3"**\dir{-},
(-.25,.75)*{\hbox{\small$\bullet$}},
(.5,1.25)*{\hbox{$\ssize h'$}},
(.5,.25)*{\hbox{$\ssize i'$}},
(-.375,.75)*{\hbox{$\ssize a$}},
\endxy}$\\
\noalign{\vskip 6pt}
(e) An $ots$-triangle with two boundary edges & (f) After the $ots$ operation
\end{tabular}
\caption{$ots$ operations in a 2-region or a minimal loop}
\label{ots operations on 2-region}
\end{figure}

For example, in Figure \ref{ots operations on 2-region} (a) and (b), an $ots$-triangle
in the interior of a 2-region or minimal loop and the result of applying an
$ots$ operation to the $ots$-triangle are shown, while in
(c) and (d), an $ots$-triangle with exactly one edge on the boundary of the 2-region
or minimal loop, together with the outcome of an $ots$ operation to this $ots$-triangle
are shown. Finally, in (e) and (f), an $ots$-triangle with two edges on the
boundary of the 2-region or minimal loop and the effect of applying an $ots$
operation to the $ots$-triangle are shown.

We remark that the plane graphs that result from each of the three possible $ots$
operations that may be performed on an $ots$-triangle that is contained in the
interior of a 2-region or minimal loop are isomorphic via an isotopy of the plane
that fixes all points in an open neighborhood of the compact set that consists of
the $ots$-triangle together with the six additional edges incident to the vertices
of the $ots$-triangle.

\begin{theorem}[Theorem 5, {[\smithrankinI]}]\label{empty min 2-region}
 Given a minimal 2-region, there is a finite sequence of $ots$ operations which 
 will transform the minimal 2-region into an empty 2-region.
\end{theorem}

As we mentioned above, if $H$ is a minimal 2-region of a 4-regular plane graph
$G$, then each $ots$ operation on $H$, as defined in Definition \ref{ots in
minimal region}, is the restriction of a (unique) $ots$ operation on $G$. 

\begin{definition}\label{collapse 2-group}
 Let $G$ be a 4-regular plane graph, and let $v$ and $w$ be the two vertices of
 a 2-(sub)group $H$ of $G$, with edges $e$ and $f$ being the edges of the subgroup.
 Form a new plane graph $G^*$ from $G$ by the following process: continuously
 shrink the two edges $e$ and $f$ to cause $v$ and $w$ to become identified,
 forming a new vertex $\overline{vw}$ located in the interior of the empty region
 bounded by $e$ and $f$. At the same time, permit every edge of $G$
 that is incident to either $v$ or $w$ to continuously extend, ultimately to have
 $v$ and/or $w$ become $\overline{vw}$ when the contraction of $e$ and $f$ is
 complete. $G^*$ has one fewer vertices and two fewer edges than $G$. We say
 $G^*$ has been created by {\it collapsing the (sub)group $H$}.
\end{definition}

\begin{proposition}\label{collapse 2-group preserves structure}
 Let $H$ be a 2-(sub)group of a 4-regular, 3-edge-connected plane graph $G$, and
 let $G^*$ denote the graph obtained by collapsing $H$. Then $G^*$ is a
 4-regular, 3-edge-connected plane graph. 
\end{proposition}

\begin{proof}
 $G^*$ is a plane graph by construction, and it is immediate that $G^*$ is
 4-regular and connected. Suppose that $G^*$ is not 3-edge-connected. Then
 there must exist two edges, say $e_1$ and $e_2$, whose deletion from $G^*$
 results in a disconnected graph. But then there must be distinct vertices
 $x$ and $y$ of $G^*$ such that every path from $x$ to $y$ uses either $e_1$ or
 $e_2$. Suppose that neither $x$ nor $y$ is
 the vertex representing the collapsed group $H$. Then $x$ and $y$ are vertices
 of $G$, and since $G$ is 3-edge-connected, there exists a path $P$ from $x$ to $y$
 that uses neither $e_1$ nor $e_2$ (where we consider $e_1$ and $e_2$ as edges
 of $G$, identifying an edge incident to $v$ or $w$ with its extended image in
 $G^*$ if necessary). If $P$ contains one of the edges of $H$, then replace
 the segment of $P$ that consists of the edge and the two end points by the
 vertex that represents $H$ to obtain a path in $G^*$ from $x$ to $y$ that uses
 neither $e_1$ nor $e_2$. Since this is not possible, it must be that one of
 $x$ or $y$ is the vertex representing $H$. Without loss of generality, suppose
 that $x$ is this vertex. Let $v$ be one of the vertices of $H$. Since $y\ne x$,
 $y$ is a vertex of $G$, and $y$ is not a vertex of $H$. As before, since $G$
 is 3-edge-connected, there is a path $P$ from $v$ to $y$ that uses neither $e_1$
 nor $e_2$. Let $w$ be the other vertex of $H$. If $w$ does not appear in $P$,
 then $P$ is a path in $G^*$ from $x$ to $y$ that uses neither $e_1$ nor $e_2$.
 Since this is not possible, $w$ must appear in $P$. But then the segment of $P$
 from $w$ to $y$ provides a path from $x$ to $y$ in $G^*$ that uses neither $e_1$
 nor $e_2$. Thus the assumption that $G^*$ is not 3-edge-connected has led to a
 contradiction, and so it followes that $G^*$ is 3-edge-connected, as required.
\end{proof}

\begin{proposition}\label{have loop}
 Let $G$ be a 4-regular, 3-edge-connected plane graph without 2-(sub)groups.
 Then either $G$ is simple or else $G$ consists of a single vertex with two
 loops.
\end{proposition}

\begin{proof}
 Suppose that $G$ is not simple. By hypothesis, $G$ has no multiply-connected
 vertices, so it must have a loop $e$ at some vertex $v$. Since $G$ is 4-regular,
 either there is a second loop based at $v$, or else that are two additional
 edges incident to $v$. But $G$ is 3-edge-connected, so there can't be two
 additional edges incident to $v$, and so there is a second loop based at $v$.
 Since $G$ is connected and 4-regular, it follows that $v$ is the only vertex of
 $G$, and the two loops are the only edges of $G$.
\end{proof}

\subsection{Condensing an alternating link diagram}

\begin{definition}\label{condensation}
 The {\it condensation} of a reduced alternating diagram $D$ of an 
 alternating link $L$ is the 4-regular plane graph $G$ that is obtained
 from $D$ by repeatedly replacing each group by a single crossing.
\end{definition}

For example, the condensation of the reduced alternating diagram of the
$(n,2)$ torus link is a single vertex with two loops. This graph shall be
denoted by $G_0$.

It was established in Proposition \ref{collapse 2-group preserves structure}
that for any 4-regular, 3-edge-connected plane graph $D$, the result of
collapsing a 2-(sub)group is again a 4-regular, 3-edge-connected plane graph.
Since the condensation of $D$ can be obtained by repeatedly collapsing
2-(sub)groups until none of the original 2-(sub)groups of $D$ remain, it
follows that the condensation of $D$ is again a 4-regular, 3-edge-connected
plane graph. Moreover, we may then repeat the process to find the
condensation of that graph. By iterating the condensation operator, we
will eventually arrive at a 4-regular, 3-edge-connected plane graph without
2-groups. By Proposition \ref{have loop}, such a graph is either $G_0$
(a single vertex with two loops), or else it is simple. Suppose that the
graph is not $G_0$. Then it may be considered to be a reduced alternating
diagram of some prime alternating link and so by Proposition \ref{existence of
min loops in cl}, it contains a minimal 2-region (but no 2-group, since all
2-groups have been collapsed). By Theorem \ref{empty min 2-region}, there is a finite
sequence of $ots$ operations that will empty the minimal 2-region. The result
of this is a 4-regular, 3-edge-connected plane graph with at least one 2-group. The
entire process can now be repeated. Since each condensation step results in a
decrease in the number of vertices, the process must eventually terminate,
with $G_0$ as the result.

The process outlined above provides a blueprint for the transformation of an $n$-crossing
reduced alternating diagram of a prime alternating link to the reduced alternating
diagram of the $(n,2)$ torus link, as we shall show next that each of the steps
in the graph transformation described above is supported by a corresponding
sequence of $T$ and/or $OTS$ operations on an $n$-crossing reduced alternating
diagram. The sequence of condensations and/or $ots$ operations on the 
graph which is the reduced alternating diagram $D$ of a prime
alternating link $L$ thereby gives rise to a (possible longer) sequence of
$T$ and/or $OTS$ operations on $D$ which transforms $D$ into the reduced alternating
diagram of the $(n,2)$ torus link.

\begin{figure}[ht]
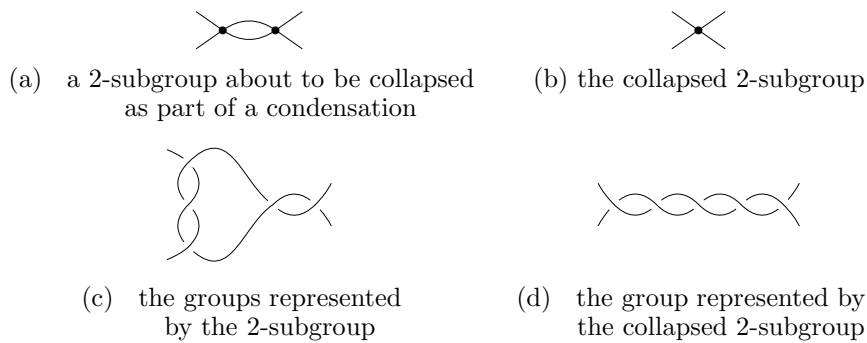

 \centering
   \begin{tabular}{c@{\hskip10pt}c}
    $\vcenter{\xy /r20pt/:,
       (0,0)="a"*=0{\hbox{$\ssize\bullet$}},
       (1,0)="b"*=0{\hbox{$\ssize\bullet$}},
       "a";"b"**\crv{(.5,.35)},
       "a";"b"**\crv{(.5,-.35)},
       "a";"a"+(-.5,.35)**\dir{-},
       "a";"a"+(-.5,-.35)**\dir{-},
       "b";"b"+(.5,.35)**\dir{-},
       "b";"b"+(.5,-.35)**\dir{-},       
       "a"*!<0pt,-6pt>{\hbox{}},
       "b"*!<0pt,-6pt>{\hbox{}},       
     \endxy}$
   &
    $\vcenter{\xy /r20pt/:,
       (0,0)="a"="b"*=0{\hbox{$\ssize\bullet$}},
       "a";"a"+(-.5,.35)**\dir{-},
       "a";"a"+(-.5,-.35)**\dir{-},
       "b";"b"+(.5,.35)**\dir{-},
       "b";"b"+(.5,-.35)**\dir{-},       
       "a"*!<0pt,-6pt>{\hbox{}},
     \endxy}$ \\
     \noalign{\vskip 6pt}
     (a) \begin{tabular}[t]{c} a 2-subgroup about to be collapsed\\
     \noalign{\vskip-1pt} as part of a condensation\end{tabular}
      & (b) the collapsed 2-subgroup\\
        \noalign{\vskip 6pt}
          $\vcenter{
           \xy /r8pt/:,
           (0,0)="c1",
           (0,2)="c2",
           (0,4)="c3",           
           (4,2)="d1",
           (6,2)="d2",
           "c1"*=<4pt>{\hbox to 4pt{\hfil}}="c1b",
           "c2"*=<4pt>{\hbox to 4pt{\hfil}}="c2b",
           "c3"*=<4pt>{\hbox to 4pt{\hfil}}="c3b",           
           "d1"*=<4pt>{\hbox to 4pt{\hfil}}="d1b",
           "d2"*=<4pt>{\hbox to 4pt{\hfil}}="d2b",
           "c1"+(-1,-.6);"c1"**\crv{(-.3,-.3)},
           "c1";"c2b"**\crv{(1,1)},
           "c2b";"c3"**\crv{(-1,3)},
           "c3";"d1b"**\crv{(2,6) & (2.5,3)},
           "d1b";"d2"**\crv{(5,1)},
           "d2";"d2"+(.75,1)**\crv{"d2"+(.5,.5)},
           "c1b";"d1"**\crv{(2,-2) & (2.5,1)},
           "d1";"d2b"**\crv{(5,3)},
           "d2b";"d2"+(.75,-1)**\crv{"d2"+(.5,-.5)},
           "c1b";"c2"**\crv{(-1,1)},
           "c2";"c3b"**\crv{(1,3)},
           "c3b";"c3"+(-1,.6)**\crv{(-.3,4.3)},
           \endxy}$
          &
    $\vcenter{\xy /r8pt/:,
           (0,0)="c1",
           (2,0)="c2",
           (4,0)="c3",           
           (6,0)="d1",
           (8,0)="d2",
           "c1"*=<4pt>{\hbox to 4pt{\hfil}}="c1b",
           "c2"*=<4pt>{\hbox to 4pt{\hfil}}="c2b",
           "c3"*=<4pt>{\hbox to 4pt{\hfil}}="c3b",           
           "d1"*=<4pt>{\hbox to 4pt{\hfil}}="d1b",
           "d2"*=<4pt>{\hbox to 4pt{\hfil}}="d2b",
           "c1"+(-.75,-1);"c1b"**\crv{"c1"+(-.5,-.5)},
           "c1b";"c2"**\crv{(1,1)},
           "c2";"c3b"**\crv{(3,-1)},
           "c3b";"d1"**\crv{(5,1)},
           "d1";"d2b"**\crv{(7,-1)},
           "d2b";"d2"+(.75,1)**\crv{"d2"+(.5,.5)},
           "c1"+(-.75,1);"c1"**\crv{"c1"+(-.5,.5)},
           "c1";"c2b"**\crv{(1,-1)},
           "c2b";"c3"**\crv{(3,1)},
           "c3";"d1b"**\crv{(5,-1)},
           "d1b";"d2"**\crv{(7,1)},
           "d2";"d2"+(.75,-1)**\crv{"d2"+(.5,-.5)},           
         \endxy}$ \\
     \noalign{\vskip 6pt}
     (c) \begin{tabular}[t]{c}the groups represented\\ \noalign{\vskip-1pt}
       by the 2-subgroup\end{tabular}
     & (d) \begin{tabular}[t]{c}the group represented by\\ \noalign{\vskip-1pt}
     the collapsed 2-subgroup\end{tabular}
   \end{tabular}
  \caption{Collapsing may require turning a subgroup}
  \label{turn and collapse}
\end{figure}

To begin with, suppose that $D$ is an $n$-crossing reduced alternating diagram
of a prime alternating link $L$. Set $G=D$, so that $G$ is a 4-regular, 3-edge-connected plane
graph with $n$ vertices, and carry out the
reduction of $G$ to the graph $G_0$ consisting of two loops on a single vertex.
Suppose that at a certain point in this process, we have arrived at an
$n$-crossing reduced alternating diagram
$D^\prime$, and the next step in the graph transformation is a condensation.
If this is the very first step in the process, then nothing need be done to $D$.
Otherwise, it will be the case that in each 2-subgroup  of the group
being condensed, each of the two edges of the 2-subgroup represents one edge from
each end of a group in $D^\prime$.
Such a situation is illustrated in Figure \ref{turn and collapse}. As shown there,
the act of collapsing the 2-subgroup in the graph corresponds to turning one
(or both) of the two subgroups of $D^\prime$ to form a larger subgroup in the 
resulting link. 

Thus we see how to manipulate a link diagram to mirror the collapse of
a 2-subgroup in the graph and hence the condensation of the link diagram.
Suppose now that we are at a stage where no further condensation is possible.
At this point, we have an $n$-crossing reduced alternating diagram $D''$ and its
condensation, $G''$, and either $G''=G_0$, and we are done, or else $G''$ is a simple graph on $k$ vertices for some
$k$ with $2\le k\le n$. Suppose that $G''$ is simple. Since $G''$ is a 4-regular, 3-edge-connected plane
graph, it is a reduced alternating diagram of some prime alternating link with $k$
crossings. By Proposition \ref{existence of min loops in cl}, $G''$ contains a
minimal 2-region, and by Theorem \ref{empty min 2-region}, there is a finite
sequence of $ots$ operations that, when applied to $G''$, empties the minimal
2-region. The result is a 4-regular, 3-edge-connected plane graph, $G_1$, on $k$
vertices that contains a 2-group (possibly as many as three 2-groups). It remains
for us to demonstrate that each of these $ots$ operations is mirrored by a sequence
of $T$ and/or $OTS$ operations on $D''$, resulting in a reduced alternating
diagram $D_1$ of a prime alternating $n$-crossing link whose condensation
is $G_1$. In fact, there are many different such sequences, and we shall
describe only one, which we shall refer to as $TOTS$.

\subsection{$TOTS$: representing $ots$ operations at the link diagram level}

\begin{figure}[ht]
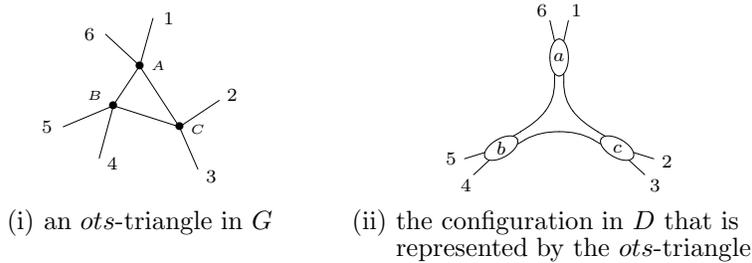

 \centering
   \begin{tabular}{c@{\hskip30pt}c}
    $\vcenter{\xy /r20pt/:,
       (0,0)="a"*=0{\hbox{$\ssize\bullet$}},
       (1.25,-.4)="b"*=0{\hbox{$\ssize\bullet$}},
       (.5,.75)="c"*=0{\hbox{$\ssize\bullet$}},
       "a";"b"**\dir{-};"c"**\dir{-};"a"**\dir{-},
       "a";"a"+(-.95,-.4)="ax"**\dir{-},
       "a";"a"+(-.27,-1)="ay"**\dir{-},
       "ax"*!<6pt,0pt>{\hbox{$\ssize 5$}},  % was 1
       "ay"*!<-5pt,2pt>{\hbox{$\ssize 4$}},  % was 6
       "c";"c"+(-.65,.6)="cx"**\dir{-},
       "c";"c"+(.25,.9)="cy"**\dir{-},
       "cx"*!<6pt,0pt>{\hbox{$\ssize 6$}},  % was 2
       "cy"*!<-6pt,0pt>{\hbox{$\ssize 1$}},  % was 3
       "b";"b"+(.75,.5)="bx"**\dir{-},
       "b";"b"+(.35,-.8)="by"**\dir{-},
       "bx"*!<-5pt,-2pt>{\hbox{$\ssize 2$}}, % was 4
       "by"*!<-5pt,3pt>{\hbox{$\ssize 3$}},  % was 5
       "a"*!<7pt,-4pt>{\hbox{\tiny$B$}},
       "b"*!<-7pt,1pt>{\hbox{\tiny$C$}},
       "c"*!<-7pt,0pt>{\hbox{\tiny$A$}},
     \endxy}$
   &
    $\vcenter{\xy /r7pt/:,
      (0,0)="b1",
      (.8660254,.5)="b2",
      "b1";"b2",{\ellipse(,.5){}},
      (8,0)="c1",
      "c1"-(.8660254,-.5)="c2",
      "c1";"c2",{\ellipse(,.5){}},
      (4,6.4)="a1",
      (4,5.4)="a2",
      "a1";"a2",{\ellipse(,.5){}},
      "b1"+(-.125,.21650635)="b1u",
      "b1"+(.125,-.21650635)="b1d",
      "b1u"+(-.8660254,-.5)+(-.125,.21650635)="b1ul"*!<5pt,0pt>{\hbox{$\ssize5$}},
      "b1d"+(-.8660254,-.5)+(.125,-.21650635)="b1dl"*!<3pt,4pt>{\hbox{$\ssize4$}},
      "b1u"+(.08660254,.05)="b1ur",
      "b1d"+(.08660254,.05)="b1dr",      
      "b1dl";"b1dr"**\dir{-},
      "b1ul";"b1ur"**\dir{-},
      "b2"+"b2"-"b1"+(-.125,.21650635)="b2u",
      "b2"+"b2"-"b1"+(.125,-.21650635)="b2d",
      "b2u"+(.8660254,.5)="b2ur",
      "b2d"+(.8660254,.5)="b2dr",
      "b2u"+(-.08660254,-.05)="b2ul",
      "b2d"+(-.08660254,-.05)="b2dl",
      "b2dl";"b2dr"**\dir{}?(.4)="b2dm",
      "b2ul";"b2ur"**\dir{}?(.4)="b2um",
      "b2dl";"b2dm"**\dir{-},
      "b2ul";"b2um"**\dir{-},      
%% done b      
      "c1"+(.125,.21650635)="c1u",
      "c1"+(-.125,-.21650635)="c1d",
      "c1u"+(.8660254,-.5)+(.125,.21650635)="c1ur"*!<-5pt,1pt>{\hbox{$\ssize2$}},
      "c1d"+(.8660254,-.5)+(-.125,-.21650635)="c1dr"*!<-4pt,4pt>{\hbox{$\ssize3$}},
      "c1u"+(-.08660254,.05)="c1ul",
      "c1d"+(-.08660254,.05)="c1dl",      
      "c1dr";"c1dl"**\dir{-},
      "c1ur";"c1ul"**\dir{-},
      "c2"+"c2"-"c1"+(.125,.21650635)="c2u",
      "c2"+"c2"-"c1"+(-.125,-.21650635)="c2d",
      "c2u"+(-.8660254,.5)="c2ul",
      "c2d"+(-.8660254,.5)="c2dl",
      "c2u"+(.08660254,-.05)="c2ur",
      "c2d"+(.08660254,-.05)="c2dr",
      "c2dr";"c2dl"**\dir{}?(.4)="c2dm",
      "c2ur";"c2ul"**\dir{}?(.4)="c2um",
      "c2dr";"c2dm"**\dir{-},
      "c2ur";"c2um"**\dir{-},      
% done c
      "a1"+(.25,0)="a1u",
      "a1"+(-.25,0)="a1d",
      "a1u"+(0,1)+(.25,0)="a1ul"*!<-3pt,-4pt>{\hbox{$\ssize1$}},
      "a1d"+(0,1)+(-.25,0)="a1dl"*!<3pt,-4pt>{\hbox{$\ssize6$}},
      "a1u"+(0,-.1)="a1ur",
      "a1d"+(0,-.1)="a1dr",      
      "a1dl";"a1dr"**\dir{-},
      "a1ul";"a1ur"**\dir{-},
      "a2"+"a2"-"a1"+(.25,0)="a2u",
      "a2"+"a2"-"a1"+(-.25,0)="a2d",
      "a2u"+(0,.1)="a2ul",
      "a2d"+(0,.1)="a2dl",
      "a2u"+(0,-1)="a2ur",
      "a2d"+(0,-1)="a2dr",
      "a2dl";"a2dr"**\dir{}?(.4)="a2dm",
      "a2ul";"a2ur"**\dir{}?(.4)="a2um",
      "a2dl";"a2dm"**\dir{-},
      "a2ul";"a2um"**\dir{-},
% done a      
      "b1";"a1"**\dir{}?(.5)="ba",
      "c1";"a1"**\dir{}?(.5)="ca",
      "b1";"c1"**\dir{}?(.5)="bc",      
      {"b1";"ca":"c1";"ba",x}="O",
      "ba";"O"**\dir{}?(.7)="bamo",
      "b2um";"a2dm"**\crv{"b2ur" & "bamo" & "a2dr"},
      "bc";"O"**\dir{}?(.7)="bcmo",
      "b2dm";"c2dm"**\crv{"b2dr" & "bcmo" & "c2dl"},
      "ca";"O"**\dir{}?(.7)="camo",
      "c2um";"a2um"**\crv{"c2ul" & "camo" & "a2ur"},             
      "b2"*{\hbox{\scriptsize$b$}},
      "c2"*{\hbox{\scriptsize$c$}},
      "a2"*{\hbox{\scriptsize$a$}},            
     \endxy}$ \\
     \noalign{\vskip 6pt}
     (i) an $ots$-triangle in $G$ & (ii) \begin{tabular}[t]{@{}l}
     the configuration in $D$ that is\\[-2pt]
     represented by the $ots$-triangle\end{tabular}
   \end{tabular}
  \caption{An $ots$ triangle in $G$ and the configuration it represents in $D$}
  \label{ots diagram}
\end{figure}

Suppose that $L$ is an $n$-crossing prime alternating link with reduced
alternating diagram $D$, and let $G$ denote the condensation of $D$.
Suppose further that $G$ is a simple graph on $k$ vertices for some
$k$ with $2\le k\le n$, and that in $G$ we have an $ots$-triangle with nodes
$A$, $B$, and $C$ representing groups in $D$ of sizes $a$, $b$, and $c$,
respectively, as shown in Figure \ref{ots diagram} (i) and (ii). Note that
we are using symbolism for groups that was introduced in [\smithrankinI], where
the major axis of the ellipse is intended to represent the two strands that are
wrapping around each other to form the group. We wish to find a finite 
sequence of $T$ and/or $OTS$ operations which, when applied to $D$, results
in a reduced alternating diagram $D'$ of a prime alternating link $L'$ such
that the condensation of $D'$ is the result of applying $ots$ to the 
$ots$-triangle with nodes $A$, $B$, and $C$ in $G$. Note that it may have 
been necessary to turn any or all of the three groups represented by $A$, 
$B$, and $C$ in order that they be aligned as shown. In such a case, these 
turn operations would be included at the beginning of the sought-after
sequence of $T$ and/or $OTS$ operations.

We begin with the case when one of the groups labelled by $A$, $B$, and $C$
is a loner.

\begin{proposition}\label{ssc}
 If the group labelled by $A$ is a loner, then there exists a finite
 sequence of $T$ and $OTS$ operations on $D$ that results in a reduced
 alternating diagram $D'$ of a prime alternating link $L'$ whose
 condensation is the result of applying $ots$ in $G$ to the $ots$-triangle 
 with nodes $A$, $B$, a group of size $b$, and $C$, a group of size $c$.
 Specifically, if $b$ is odd, then there
 exists a sequence of $T$ and $OTS$ operations that transforms Figure
 \ref{1 loner, 2 nonloners} (a) into Figure \ref{1 loner, 2 nonloners} (b),
 while if $b$ is even, then there exists a sequence of $T$ and $OTS$
 operations that transforms Figure \ref{1 loner, 2 nonloners} (a) into
 Figure \ref{1 loner, 2 nonloners} (c).
 
\begin{figure}[ht]
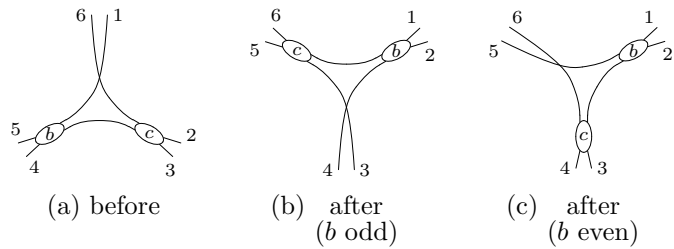

 \centering
   \begin{tabular}{c@{\hskip20pt}c@{\hskip20pt}c}
    $\vcenter{\xy /r6pt/:,
      (0,0)="b1",
      (.8660254,.5)="b2",
      "b1";"b2",{\ellipse(,.5){}},
      (8,0)="c1",
      "c1"-(.8660254,-.5)="c2",
      "c1";"c2",{\ellipse(,.5){}},
      (4,7)="a1",
      (4,6)="a2",
 %     "a1";"a2",{\ellipse(,.5){}},
      "b1"+(-.125,.21650635)="b1u",
      "b1"+(.125,-.21650635)="b1d",
      "b1u"+(-.8660254,-.5)+(-.125,.21650635)="b1ul",
      "b1d"+(-.8660254,-.5)+(.125,-.21650635)="b1dl",
      "b1u"+(.08660254,.05)="b1ur",
      "b1d"+(.08660254,.05)="b1dr",      
      "b1dl";"b1dr"**\dir{-},
      "b1ul";"b1ur"**\dir{-},
      "b2"+"b2"-"b1"+(-.125,.21650635)="b2u",
      "b2"+"b2"-"b1"+(.125,-.21650635)="b2d",
      "b2u"+(.8660254,.5)="b2ur",
      "b2d"+(.8660254,.5)="b2dr",
      "b2u"+(-.08660254,-.05)="b2ul",
      "b2d"+(-.08660254,-.05)="b2dl",
      "b2dl";"b2dr"**\dir{}?(.4)="b2dm",
      "b2ul";"b2ur"**\dir{}?(.4)="b2um",
      "b2dl";"b2dm"**\dir{-},
      "b2ul";"b2um"**\dir{-},      
%% done b      
      "c1"+(.125,.21650635)="c1u",
      "c1"+(-.125,-.21650635)="c1d",
      "c1u"+(.8660254,-.5)+(.125,.21650635)="c1ur",
      "c1d"+(.8660254,-.5)+(-.125,-.21650635)="c1dr",
      "c1u"+(-.08660254,.05)="c1ul",
      "c1d"+(-.08660254,.05)="c1dl",      
      "c1dr";"c1dl"**\dir{-},
      "c1ur";"c1ul"**\dir{-},
      "c2"+"c2"-"c1"+(.125,.21650635)="c2u",
      "c2"+"c2"-"c1"+(-.125,-.21650635)="c2d",
      "c2u"+(-.8660254,.5)="c2ul",
      "c2d"+(-.8660254,.5)="c2dl",
      "c2u"+(.08660254,-.05)="c2ur",
      "c2d"+(.08660254,-.05)="c2dr",
      "c2dr";"c2dl"**\dir{}?(.4)="c2dm",
      "c2ur";"c2ul"**\dir{}?(.4)="c2um",
      "c2dr";"c2dm"**\dir{-},
      "c2ur";"c2um"**\dir{-},      
% done c
      "a1"+(.25,0)="a1u",
      "a1"+(-.25,0)="a1d",
      "a1u"+(0,1)+(.25,0)="a1ul",
      "a1d"+(0,1)+(-.25,0)="a1dl",
      "a1u"+(0,-.1)="a1ur",
      "a1d"+(0,-.1)="a1dr",      
      "a1dl";"a1dr"**\dir{},
      "a1ul";"a1ur"**\dir{},
      "a2"+"a2"-"a1"+(.25,0)="a2u",
      "a2"+"a2"-"a1"+(-.25,0)="a2d",
      "a2u"+(0,.1)="a2ul",
      "a2d"+(0,.1)="a2dl",
      "a2u"+(0,-1)="a2ur",
      "a2d"+(0,-1)="a2dr",
      "a2dl";"a2dr"**\dir{}?(.4)="a2dm",
      "a2ul";"a2ur"**\dir{}?(.4)="a2um",
      "a2dl";"a2dm"**\dir{},
      "a2ul";"a2um"**\dir{},
% done a      
      "b1";"a1"**\dir{}?(.5)="ba",
      "c1";"a1"**\dir{}?(.5)="ca",
      "b1";"c1"**\dir{}?(.5)="bc",      
      {"b1";"ca":"c1";"ba",x}="O",
      "ba";"O"**\dir{}?(.7)="bamo",
      "b2um";"a1ul"**\crv{"b2ur" & "bamo" & "a2ur"},
      "bc";"O"**\dir{}?(.6)="bcmo",
      "b2dm";"c2dm"**\crv{"b2dr" & "bcmo" & "c2dl"},
      "ca";"O"**\dir{}?(.7)="camo",
      "c2um";"a1dl"**\crv{"c2ul" & "camo" & "a2dr"},             
      "b2"*{\hbox{\scriptsize$b$}},
      "c2"*{\hbox{\scriptsize$c$}},
%      "a2"*{\hbox{\tiny$a$}},            
      "a1ul"*!<-4pt,0pt>{\hbox{\scriptsize$1$}},
      "a1dl"*!<4pt,0pt>{\hbox{\scriptsize$6$}},
      "c1ur"*!<-4pt,-3pt>{\hbox{\scriptsize$2$}},
      "c1dr"*!<1pt,5pt>{\hbox{\scriptsize$3$}},
      "b1dl"*!<-3pt,4pt>{\hbox{\scriptsize$4$}},
      "b1ul"*!<1pt,-5pt>{\hbox{\scriptsize$5$}},            
     \endxy}$
    &
    $\vcenter{\xy /l6pt/:,
      (0,0)="b1",
      (.8660254,.5)="b2",
      "b1";"b2",{\ellipse(,.5){}},
      (8,0)="c1",
      "c1"-(.8660254,-.5)="c2",
      "c1";"c2",{\ellipse(,.5){}},
      (4,7)="a1",
      (4,6)="a2",
 %     "a1";"a2",{\ellipse(,.5){}},
      "b1"+(-.125,.21650635)="b1u",
      "b1"+(.125,-.21650635)="b1d",
      "b1u"+(-.8660254,-.5)+(-.125,.21650635)="b1ul",
      "b1d"+(-.8660254,-.5)+(.125,-.21650635)="b1dl",
      "b1u"+(.08660254,.05)="b1ur",
      "b1d"+(.08660254,.05)="b1dr",      
      "b1dl";"b1dr"**\dir{-},
      "b1ul";"b1ur"**\dir{-},
      "b2"+"b2"-"b1"+(-.125,.21650635)="b2u",
      "b2"+"b2"-"b1"+(.125,-.21650635)="b2d",
      "b2u"+(.8660254,.5)="b2ur",
      "b2d"+(.8660254,.5)="b2dr",
      "b2u"+(-.08660254,-.05)="b2ul",
      "b2d"+(-.08660254,-.05)="b2dl",
      "b2dl";"b2dr"**\dir{}?(.4)="b2dm",
      "b2ul";"b2ur"**\dir{}?(.4)="b2um",
      "b2dl";"b2dm"**\dir{-},
      "b2ul";"b2um"**\dir{-},      
%% done b      
      "c1"+(.125,.21650635)="c1u",
      "c1"+(-.125,-.21650635)="c1d",
      "c1u"+(.8660254,-.5)+(.125,.21650635)="c1ur",
      "c1d"+(.8660254,-.5)+(-.125,-.21650635)="c1dr",
      "c1u"+(-.08660254,.05)="c1ul",
      "c1d"+(-.08660254,.05)="c1dl",      
      "c1dr";"c1dl"**\dir{-},
      "c1ur";"c1ul"**\dir{-},
      "c2"+"c2"-"c1"+(.125,.21650635)="c2u",
      "c2"+"c2"-"c1"+(-.125,-.21650635)="c2d",
      "c2u"+(-.8660254,.5)="c2ul",
      "c2d"+(-.8660254,.5)="c2dl",
      "c2u"+(.08660254,-.05)="c2ur",
      "c2d"+(.08660254,-.05)="c2dr",
      "c2dr";"c2dl"**\dir{}?(.4)="c2dm",
      "c2ur";"c2ul"**\dir{}?(.4)="c2um",
      "c2dr";"c2dm"**\dir{-},
      "c2ur";"c2um"**\dir{-},      
% done c
      "a1"+(.25,0)="a1u",
      "a1"+(-.25,0)="a1d",
      "a1u"+(0,1)+(.25,0)="a1ul",
      "a1d"+(0,1)+(-.25,0)="a1dl",
      "a1u"+(0,-.1)="a1ur",
      "a1d"+(0,-.1)="a1dr",      
      "a1dl";"a1dr"**\dir{},
      "a1ul";"a1ur"**\dir{},
      "a2"+"a2"-"a1"+(.25,0)="a2u",
      "a2"+"a2"-"a1"+(-.25,0)="a2d",
      "a2u"+(0,.1)="a2ul",
      "a2d"+(0,.1)="a2dl",
      "a2u"+(0,-1)="a2ur",
      "a2d"+(0,-1)="a2dr",
      "a2dl";"a2dr"**\dir{}?(.4)="a2dm",
      "a2ul";"a2ur"**\dir{}?(.4)="a2um",
      "a2dl";"a2dm"**\dir{},
      "a2ul";"a2um"**\dir{},
% done a      
      "b1";"a1"**\dir{}?(.5)="ba",
      "c1";"a1"**\dir{}?(.5)="ca",
      "b1";"c1"**\dir{}?(.5)="bc",      
      {"b1";"ca":"c1";"ba",x}="O",
      "ba";"O"**\dir{}?(.7)="bamo",
      "b2um";"a1ul"**\crv{"b2ur" & "bamo" & "a2ur"},
      "bc";"O"**\dir{}?(.6)="bcmo",
      "b2dm";"c2dm"**\crv{"b2dr" & "bcmo" & "c2dl"},
      "ca";"O"**\dir{}?(.7)="camo",
      "c2um";"a1dl"**\crv{"c2ul" & "camo" & "a2dr"},             
      "b2"*{\hbox{\scriptsize$b$}},
      "c2"*{\hbox{\scriptsize$c$}},
%      "a2"*{\hbox{\scriptsize$a$}},            
      "a1ul"*!<4pt,0pt>{\hbox{\scriptsize$4$}},
      "a1dl"*!<-4pt,0pt>{\hbox{\scriptsize$3$}},
      "c1ur"*!<4pt,3pt>{\hbox{\scriptsize$5$}},
      "c1dr"*!<-1pt,-5pt>{\hbox{\scriptsize$6$}},
      "b1dl"*!<3pt,-4pt>{\hbox{\scriptsize$1$}},
      "b1ul"*!<-1pt,5pt>{\hbox{\scriptsize$2$}},            
     \endxy}$
    &
    $\vcenter{\xy /r6pt/:, a(60):,
      (0,0)="b1",
      (.8660254,.5)="b2",
      "b1";"b2",{\ellipse(,.5){}},
      (8,0)="c1",
      "c1"-(.8660254,-.5)="c2",
      "c1";"c2",{\ellipse(,.5){}},
      (4,7)="a1",
      (4,6)="a2",
 %     "a1";"a2",{\ellipse(,.5){}},
      "b1"+(-.125,.21650635)="b1u",
      "b1"+(.125,-.21650635)="b1d",
      "b1u"+(-.8660254,-.5)+(-.125,.21650635)="b1ul",
      "b1d"+(-.8660254,-.5)+(.125,-.21650635)="b1dl",
      "b1u"+(.08660254,.05)="b1ur",
      "b1d"+(.08660254,.05)="b1dr",      
      "b1dl";"b1dr"**\dir{-},
      "b1ul";"b1ur"**\dir{-},
      "b2"+"b2"-"b1"+(-.125,.21650635)="b2u",
      "b2"+"b2"-"b1"+(.125,-.21650635)="b2d",
      "b2u"+(.8660254,.5)="b2ur",
      "b2d"+(.8660254,.5)="b2dr",
      "b2u"+(-.08660254,-.05)="b2ul",
      "b2d"+(-.08660254,-.05)="b2dl",
      "b2dl";"b2dr"**\dir{}?(.4)="b2dm",
      "b2ul";"b2ur"**\dir{}?(.4)="b2um",
      "b2dl";"b2dm"**\dir{-},
      "b2ul";"b2um"**\dir{-},      
%% done b      
      "c1"+(.125,.21650635)="c1u",
      "c1"+(-.125,-.21650635)="c1d",
      "c1u"+(.8660254,-.5)+(.125,.21650635)="c1ur",
      "c1d"+(.8660254,-.5)+(-.125,-.21650635)="c1dr",
      "c1u"+(-.08660254,.05)="c1ul",
      "c1d"+(-.08660254,.05)="c1dl",      
      "c1dr";"c1dl"**\dir{-},
      "c1ur";"c1ul"**\dir{-},
      "c2"+"c2"-"c1"+(.125,.21650635)="c2u",
      "c2"+"c2"-"c1"+(-.125,-.21650635)="c2d",
      "c2u"+(-.8660254,.5)="c2ul",
      "c2d"+(-.8660254,.5)="c2dl",
      "c2u"+(.08660254,-.05)="c2ur",
      "c2d"+(.08660254,-.05)="c2dr",
      "c2dr";"c2dl"**\dir{}?(.4)="c2dm",
      "c2ur";"c2ul"**\dir{}?(.4)="c2um",
      "c2dr";"c2dm"**\dir{-},
      "c2ur";"c2um"**\dir{-},      
% done c
      "a1"+(.25,0)="a1u",
      "a1"+(-.25,0)="a1d",
      "a1u"+(0,1)+(.25,0)="a1ul",
      "a1d"+(0,1)+(-.25,0)="a1dl",
      "a1u"+(0,-.1)="a1ur",
      "a1d"+(0,-.1)="a1dr",      
      "a1dl";"a1dr"**\dir{},
      "a1ul";"a1ur"**\dir{},
      "a2"+"a2"-"a1"+(.25,0)="a2u",
      "a2"+"a2"-"a1"+(-.25,0)="a2d",
      "a2u"+(0,.1)="a2ul",
      "a2d"+(0,.1)="a2dl",
      "a2u"+(0,-1)="a2ur",
      "a2d"+(0,-1)="a2dr",
      "a2dl";"a2dr"**\dir{}?(.4)="a2dm",
      "a2ul";"a2ur"**\dir{}?(.4)="a2um",
      "a2dl";"a2dm"**\dir{},
      "a2ul";"a2um"**\dir{},
% done a      
      "b1";"a1"**\dir{}?(.5)="ba",
      "c1";"a1"**\dir{}?(.5)="ca",
      "b1";"c1"**\dir{}?(.5)="bc",      
      {"b1";"ca":"c1";"ba",x}="O",
      "ba";"O"**\dir{}?(.7)="bamo",
      "b2um";"a1ul"**\crv{"b2ur" & "bamo" & "a2ur"},
      "bc";"O"**\dir{}?(.6)="bcmo",
      "b2dm";"c2dm"**\crv{"b2dr" & "bcmo" & "c2dl"},
      "ca";"O"**\dir{}?(.7)="camo",
      "c2um";"a1dl"**\crv{"c2ul" & "camo" & "a2dr"},             
      "b2"*{\hbox{\scriptsize$c$}},
      "c2"*{\hbox{\scriptsize$b$}},
%      "a2"*{\hbox{\scriptsize$a$}},            
      "a1ul"*!<-3pt,-4pt>{\hbox{\scriptsize$6$}},
      "a1dl"*!<3pt,4pt>{\hbox{\scriptsize$5$}},
      "c1ur"*!<3pt,-4pt>{\hbox{\scriptsize$1$}},
      "c1dr"*!<-1pt,5pt>{\hbox{\scriptsize$2$}},
      "b1dl"*!<-4pt,0pt>{\hbox{\scriptsize$3$}},
      "b1ul"*!<4pt,0pt>{\hbox{\scriptsize$4$}},            
     \endxy}$    
     \\
     \noalign{\vskip 6pt}
     (a) before  & (b) \begin{tabular}[t]{@{}c}after\\[-2pt]
      ($b$ odd)\end{tabular} & (c) \begin{tabular}[t]{@{}c}after\\[-2pt]
      ($b$ even)\end{tabular}
   \end{tabular}
  \caption{The representation of $ots$ with at least one loner}
  \label{1 loner, 2 nonloners}
\end{figure}
\end{proposition}

\begin{proof}
 The proof will be by induction on $b$. If $b=1$, then the $ots$-triangle
 corresponds to an $OTS$ 6-tangle in the link, where two of the three
 crossings of the $OTS$ 6-tangle are loners, while the third is an end
 crossing of a group of size $c$ (it is possible that $c>1$). The arc
 with endpoints the groups labelled by $A$ and $B$ may be $OTS$'ed across
 each crossing in the group labelled $C$ in turn, requiring $c$ $OTS$
 operations in all. The result is as shown in Figure \ref{1 loner, 2 nonloners}
 (b) (with $b=1$), which is the desired result since $b=1$ is odd.

 Suppose now that $b\ge1$ is an integer for which the statement of the
 proposition is valid for all groups $B$ of size at most $b$, and consider 
 an $ots$-triangle in which vertex
 $A$ represents a loner, the crossing $x$ in Figure \ref{induction for
 loner ots} (a), $B$ represents a group of size $b+1$, and $C$
 represents a group of arbitrary size $c$. Let $B_1$ represent the subgroup of size
 $b$ that is obtained from the group $B$ by separating off the end crossing
 (denoted by $y$) whose arcs are connected to groups $A$ and $C$, as
 shown in Figure \ref{induction for loner ots} (a). Then $y$, $x$ and $C$
 form an $ots$-triangle with $y$ and $x$ single crossings. This is handled
 in the same way as in the base case, namely the arc with endpoints
 $y$ and $x$ is $OTS$'ed across the group $C$. The result is as shown in
 Figure \ref{induction for loner ots} (b).
 
\begin{figure}[ht]
 \centering
   \begin{tabular}{c@{\hskip20pt}c@{\hskip20pt}c}
    $\vcenter{\xy /r6pt/:,
      (-1,-1)="b1",
      (.8660254,.5)+(-1,-1)="b2",
      "b1";"b2",{\ellipse(,.5){}},
      (8,0)="c1",
      "c1"-(.8660254,-.5)="c2",
      "c1";"c2",{\ellipse(,.5){}},
      (4,8)="a1",
      (4,7)="a2",
 %     "a1";"a2",{\ellipse(,.5){}},
      "b1"+(-.125,.21650635)="b1u",
      "b1"+(.125,-.21650635)="b1d",
      "b1u"+(-.8660254,-.5)+(-.125,.21650635)="b1ul",
      "b1d"+(-.8660254,-.5)+(.125,-.21650635)="b1dl",
      "b1u"+(.08660254,.05)="b1ur",
      "b1d"+(.08660254,.05)="b1dr",      
      "b1dl";"b1dr"**\dir{-},
      "b1ul";"b1ur"**\dir{-},
      "b2"+"b2"-"b1"+(-.125,.21650635)="b2u",
      "b2"+"b2"-"b1"+(.125,-.21650635)="b2d",
      "b2u"+(.8660254,.5)="b2ur",
      "b2d"+(.8660254,.5)="b2dr",
      "b2ur"+(.8660254,.5)="b2ury",
      "b2dr"+(.8660254,.5)+(.3,.3)="b2dry",
      "b2u"+(-.08660254,-.05)="b2ul",
      "b2d"+(-.08660254,-.05)="b2dl",
      "b2dl";"b2dr"**\dir{}?(.4)="b2dm",
      "b2ul";"b2ur"**\dir{}?(.4)="b2um",
      "b2dl";"b2dm"**\dir{-},
      "b2ul";"b2um"**\dir{-},      
%% done b      
      "c1"+(.125,.21650635)="c1u",
      "c1"+(-.125,-.21650635)="c1d",
      "c1u"+(.8660254,-.5)+(.125,.21650635)="c1ur",
      "c1d"+(.8660254,-.5)+(-.125,-.21650635)="c1dr",
      "c1u"+(-.08660254,.05)="c1ul",
      "c1d"+(-.08660254,.05)="c1dl",      
      "c1dr";"c1dl"**\dir{-},
      "c1ur";"c1ul"**\dir{-},
      "c2"+"c2"-"c1"+(.125,.21650635)="c2u",
      "c2"+"c2"-"c1"+(-.125,-.21650635)="c2d",
      "c2u"+(-.8660254,.5)="c2ul",
      "c2d"+(-.8660254,.5)="c2dl",
      "c2u"+(.08660254,-.05)="c2ur",
      "c2d"+(.08660254,-.05)="c2dr",
      "c2dr";"c2dl"**\dir{}?(.4)="c2dm",
      "c2ur";"c2ul"**\dir{}?(.4)="c2um",
      "c2dr";"c2dm"**\dir{-},
      "c2ur";"c2um"**\dir{-},      
% done c
      "a1"+(.25,0)="a1u",
      "a1"+(-.25,0)="a1d",
      "a1u"+(0,1)+(.25,0)="a1ul",
      "a1d"+(0,1)+(-.25,0)="a1dl",
      "a1u"+(0,-.1)="a1ur",
      "a1d"+(0,-.1)="a1dr",      
      "a1dl";"a1dr"**\dir{},
      "a1ul";"a1ur"**\dir{},
      "a2"+"a2"-"a1"+(.25,0)="a2u",
      "a2"+"a2"-"a1"+(-.25,0)="a2d",
      "a2u"+(0,.1)="a2ul",
      "a2d"+(0,.1)="a2dl",
      "a2u"+(0,-1)="a2ur",
      "a2d"+(0,-1)="a2dr",
      "a2dl";"a2dr"**\dir{}?(.4)="a2dm",
      "a2ul";"a2ur"**\dir{}?(.4)="a2um",
      "a2dl";"a2dm"**\dir{},
      "a2ul";"a2um"**\dir{},
% done a      
      "b1";"a1"**\dir{}?(.5)="ba",
      "c1";"a1"**\dir{}?(.5)="ca",
      "b1";"c1"**\dir{}?(.5)="bc",      
      {"b1";"ca":"c1";"ba",x}="O",
      "ba";"O"**\dir{}?(.7)="bamo",
      "b2um";"b2dry"**\crv{"b2ur"},
      "b2dm";"b2ury"**\crv{"b2dr"},      
      "b2ury";"a1ul"**\crv{ "bamo" & "a2ur"},
      "bc";"O"**\dir{}?(.6)="bcmo",
      "b2dry";"c2dm"**\crv{ "bcmo" & "c2dl"},
      "ca";"O"**\dir{}?(.7)="camo",
      "c2um";"a1dl"**\crv{"c2ul" & "camo" & "a2dr"},             
      "b2"*{\hbox{\scriptsize$b$}},
      "c2"*{\hbox{\scriptsize$c$}},
%      "a2"*{\hbox{\tiny$a$}},
      "b2"*!<-12pt,-13pt>{\hbox{\scriptsize$y$}},
      "a2"*!<5pt,10pt>{\hbox{\scriptsize$x$}},      
      "a1ul"*!<-4pt,0pt>{\hbox{\scriptsize$1$}},
      "a1dl"*!<4pt,0pt>{\hbox{\scriptsize$6$}},
      "c1ur"*!<-4pt,-3pt>{\hbox{\scriptsize$2$}},
      "c1dr"*!<1pt,5pt>{\hbox{\scriptsize$3$}},
      "b1dl"*!<-3pt,4pt>{\hbox{\scriptsize$4$}},
      "b1ul"*!<1pt,-5pt>{\hbox{\scriptsize$5$}},            
     \endxy}$
    &
  %%%%%%%%%%%%%%%%%%%%%%%%%%%%%%%%%%%%%%%%%%%%%%%%%%%%%%%%%%
    $\vcenter{\xy /r6pt/:,
      (0,0)="b1",
      (.8660254,.5)="b2",
      "b1";"b2",{\ellipse(,.5){}},
      (8,0)="c1",
      "c1"-(.8660254,-.5)="c2",
%      "c1";"c2",{\ellipse(,.5){}},
      (4,8)="a1",
      (4,7)="a2",
 %     "a1";"a2",{\ellipse(,.5){}},
      "b1"+(-.125,.21650635)="b1u",
      "b1"+(.125,-.21650635)="b1d",
      "b1u"+(-.8660254,-.5)+(-.125,.21650635)="b1ul",
      "b1d"+(-.8660254,-.5)+(.125,-.21650635)="b1dl",
      "b1u"+(.08660254,.05)="b1ur",
      "b1d"+(.08660254,.05)="b1dr",      
      "b1dl";"b1dr"**\dir{-},
      "b1ul";"b1ur"**\dir{-},
      "b2"+"b2"-"b1"+(-.125,.21650635)="b2u",
      "b2"+"b2"-"b1"+(.125,-.21650635)="b2d",
      "b2u"+(.8660254,.5)="b2ur",
      "b2d"+(.8660254,.5)="b2dr",
      "b2u"+(-.08660254,-.05)="b2ul",
      "b2d"+(-.08660254,-.05)="b2dl",
      "b2dl";"b2dr"**\dir{}?(.4)="b2dm",
      "b2ul";"b2ur"**\dir{}?(.4)="b2um",
      "b2dl";"b2dm"**\dir{-},
      "b2ul";"b2um"**\dir{-},      
%% done b      
      "c1"+(.125,.21650635)="c1u",
      "c1"+(-.125,-.21650635)="c1d",
      "c1u"+(.8660254,-.5)+(.125,.21650635)="c1ura",
      "c1d"+(.8660254,-.5)+(-.125,-.21650635)="c1dra",
%%%%  move C left and up one and a half units
      "c1"+(-.8660254,.5)+(-.4330127,.25)="c1",
      "c1"-(.8660254,-.5)="c2",
      "c1";"c2",{\ellipse(,.5){}},
      "c1"+(.125,.21650635)="c1u",
      "c1"+(-.125,-.21650635)="c1d",
      "c1u"+(.8660254,-.5)+(.125,.21650635)="c1ur",
      "c1d"+(.8660254,-.5)+(-.125,-.21650635)="c1dr",
      "c1u"+(-.08660254,.05)="c1ul",
      "c1d"+(-.08660254,.05)="c1dl",      
      "c1dra"+(.8660254,-.5);"c1dl"**\crv{"c1dr"+<-2pt,3pt>},
      "c1ura"+(.8660254,-.5)+(1,5);"c1ul"**\crv{"c1ur"+<-2pt,-2pt>},
      "c2"+"c2"-"c1"+(.125,.21650635)="c2u",
      "c2"+"c2"-"c1"+(-.125,-.21650635)="c2d",
      "c2u"+(-.8660254,.5)="c2ul",
      "c2d"+(-.8660254,.5)="c2dl",
      "c2u"+(.08660254,-.05)="c2ur",
      "c2d"+(.08660254,-.05)="c2dr",
      "c2dr";"c2dl"**\dir{}?(.4)="c2dm",
      "c2ur";"c2ul"**\dir{}?(.4)="c2um",
      "c2dr";"c2dm"**\dir{-},
      "c2ur";"c2um"**\dir{-},      
% done c
      "a1"+(.25,0)="a1u",
      "a1"+(-.25,0)="a1d",
      "a1u"+(0,1)+(.25,0)="a1ul",
      "a1d"+(0,1)+(-.25,0)="a1dl",
      "a1u"+(0,-.1)="a1ur",
      "a1d"+(0,-.1)="a1dr",      
      "a1dl";"a1dr"**\dir{},
      "a1ul";"a1ur"**\dir{},
      "a2"+"a2"-"a1"+(.25,0)="a2u",
      "a2"+"a2"-"a1"+(-.25,0)="a2d",
      "a2u"+(0,.1)="a2ul",
      "a2d"+(0,.1)="a2dl",
      "a2u"+(0,-1)="a2ur",
      "a2d"+(0,-1)="a2dr",
      "a2dl";"a2dr"**\dir{}?(.4)="a2dm",
      "a2ul";"a2ur"**\dir{}?(.4)="a2um",
      "a2dl";"a2dm"**\dir{},
      "a2ul";"a2um"**\dir{},
% done a      
      "b1";"a1"**\dir{}?(.5)="ba",
      "c1";"a1"**\dir{}?(.5)="ca",
      "b1";"c1"**\dir{}?(.5)="bc",      
      {"b1";"ca":"c1";"ba",x}="O",
      "ba";"O"**\dir{}?(.7)="bamo",
      "b2dm";"a1ul"**\crv{"b2dm"+"b2dm"-"b2dl"+
         "b2dm"-"b2dl"+"b2dm"-"b2dl" & "b2"+(4,-1)
        & "c1"+(-1,-2)+(.4330127,-.25) & "c1" + (1.25,-1.25) & "c1"+(4,3)
        +(.4330127,-.25) &"a2ur"},
      "bc";"O"**\dir{}?(.6)="bcmo",
      "b2um";"c2dm"**\crv{"bcmo"+(-.75,.5) & "c2dl"},
      "ca";"O"**\dir{}?(.7)="camo",
      "c2um";"a1dl"**\crv{"c2ul" & "camo" & "a2dr"},             
      "b2"*{\hbox{\scriptsize$b$}},
      "c2"*{\hbox{\scriptsize$c$}},
%      "a2"*{\hbox{\scriptsize$a$}},            
      "c1ura"+(.8660254,-.5)+(-1.75,.75)*!<-12pt,-13pt>{\hbox{\scriptsize$y$}},
      "c1dra"+(.8660254,-.5)+(-1.75,2.75)*!<5pt,10pt>{\hbox{\scriptsize$x$}},      
      "a1ul"*!<-4pt,0pt>{\hbox{\scriptsize$1$}},
      "a1dl"*!<4pt,0pt>{\hbox{\scriptsize$6$}},
      "c1ura"+(.8660254,-.5)+(1,5)*!<-4pt,-3pt>{\hbox{\scriptsize$2$}},
      "c1dra"+(.8660254,-.5)*!<1pt,5pt>{\hbox{\scriptsize$3$}},
      "b1dl"*!<-3pt,4pt>{\hbox{\scriptsize$4$}},
      "b1ul"*!<1pt,-5pt>{\hbox{\scriptsize$5$}},            
     \endxy}$
    &
%%%%%%%%%%%%%%%%%%%%%%%%%%%%%%%%%%%%%%%%%%%%%%%%%%%%%%%%%%%%%%%%%%
    $\vcenter{\xy /r6pt/:,
      (0,0)="b1",
      (.8660254,.5)="b2",
      "b1";"b2",{\ellipse(,.5){}},
      (13,5)="v2",
      (10,-1)="v3",
      (6,2.5)="c1",
      "c1"+(.8660254,.5)="c2",
      "c1";"c2",{\ellipse(,.5){}},
      (4,8)="a1",
      (4,7)="a2",
      "b1"+(-.125,.21650635)="b1u",
      "b1"+(.125,-.21650635)="b1d",
      "b1u"+(-.8660254,-.5)+(-.125,.21650635)="b1ul",
      "b1d"+(-.8660254,-.5)+(.125,-.21650635)="b1dl",
      "b1u"+(.08660254,.05)="b1ur",
      "b1d"+(.08660254,.05)="b1dr",      
      "b1dl";"b1dr"**\dir{-},
      "b1ul";"b1ur"**\dir{-},
      "b2"+"b2"-"b1"+(-.125,.21650635)="b2u",
      "b2"+"b2"-"b1"+(.125,-.21650635)="b2d",
      "b2u"+(.8660254,.5)="b2ur",
      "b2d"+(.8660254,.5)="b2dr",
      "b2u"+(-.08660254,-.05)="b2ul",
      "b2d"+(-.08660254,-.05)="b2dl",
      "b2dl";"b2dr"**\dir{}?(.4)="b2dm",
      "b2ul";"b2ur"**\dir{}?(.4)="b2um",
      "b2ul";"b2um"**\dir{-},      
%% done b      
      "c1"+(-.125,.21650635)="c1u",
      "c1"+(.125,-.21650635)="c1d",
      "c1u"+(-.8660254,-.5)+(-.125,.21650635)="c1ul",
      "c1d"+(-.8660254,-.5)+(.125,-.21650635)="c1dl",
      "c1u"+(.08660254,.05)="c1ur",
      "c1d"+(.08660254,.05)="c1dr",      
      "c2"+"c2"-"c1"+(-.125,.21650635)="c2u",
      "c2"+"c2"-"c1"+(.125,-.21650635)="c2d",
      "c2u"+(.8660254,.5)="c2ur",
      "c2d"+(.8660254,.5)="c2dr",
      "c2u"+(-.08660254,-.05)="c2ul",
      "c2d"+(-.08660254,-.05)="c2dl",
% done c
      "a1"+(.25,0)="a1u",
      "a1"+(-.25,0)="a1d",
      "a1u"+(0,1)+(.25,0)="a1ul",
      "a1d"+(0,1)+(-.25,0)="a1dl",
      "a1u"+(0,-.1)="a1ur",
      "a1d"+(0,-.1)="a1dr",      
      "a1dl";"a1dr"**\dir{},
      "a1ul";"a1ur"**\dir{},
      "a2"+"a2"-"a1"+(.25,0)="a2u",
      "a2"+"a2"-"a1"+(-.25,0)="a2d",
      "a2u"+(0,.1)="a2ul",
      "a2d"+(0,.1)="a2dl",
      "a2u"+(0,-1)="a2ur",
      "a2d"+(0,-1)="a2dr",
      "a2dl";"a2dr"**\dir{}?(.4)="a2dm",
      "a2ul";"a2ur"**\dir{}?(.4)="a2um",
      "a2dl";"a2dm"**\dir{},
      "a2ul";"a2um"**\dir{},
% done a      
      "b1";"a1"**\dir{}?(.5)="ba",
      "c1";"a1"**\dir{}?(.5)="ca",
      "b1";"c1"**\dir{}?(.5)="bc",      
      "b2ur";"b2ur"+"b2ur"-"b2ul"+"b2ur"-"b2ul"+"b2ur"-"b2ul"**\dir{}?(.6)="xx1",
      "b2u";"a1ul"**\crv{"xx1" & 
         "c1"+(1,-5) & "c1" + (7,.5) & "c2"+(3.5,3.5)  & "c2"+(-2.5,3)},
      "b2dl";"c1ur"**\crv{"b2dl"+(1,0) &  "c1ul"+(-1,0)},
      "c1dr";"v3"+(-.2,0)**\crv{"c1dr"+(-1,-1) },             
      "c2ul";"a1dl"+(-.2,0)**\crv{"c2ul"+(1,1) & "a1dl"+(0,-5) },             
      "c2dl";"v2"**\crv{"c2dl"+(1.5,.2) & "v2"+(-1,.2)},
      "b2"*{\hbox{\scriptsize$b'$}},
      "c2"*{\hbox{\scriptsize$c$}},
%      "a2"*{\hbox{\scriptsize$a$}},            
%      "c1ura"+(.8660254,-.5)+(-1.75,.75)*!<-12pt,-13pt>{\hbox{\scriptsize$y$}},
%      "c1dra"+(.8660254,-.5)+(-1.75,2.75)*!<5pt,10pt>{\hbox{\scriptsize$x$}},      
      "a1ul"*!<-4pt,0pt>{\hbox{\scriptsize$1$}},
      "a1dl"*!<4pt,0pt>{\hbox{\scriptsize$6$}},
      "v2"*!<-4pt,-3pt>{\hbox{\scriptsize$2$}},
      "v3"*!<1pt,5pt>{\hbox{\scriptsize$3$}},
      "b1dl"*!<-3pt,4pt>{\hbox{\scriptsize$4$}},
      "b1ul"*!<1pt,-5pt>{\hbox{\scriptsize$5$}},            
      "b2ur"+(1.1,.8)*{\hbox{\scriptsize$t$}},
      "v3"+(-1.5,1.4)*{\hbox{\scriptsize$x$}},
      "c2"+(5.5,1)*{\hbox{\scriptsize$y$}},
     \endxy}$
     \\
     \noalign{\vskip 6pt}
     (a)   & (b)  & (c)
   \end{tabular}
  \caption{}
  \label{induction for loner ots}
\end{figure}

Turn group $C$, and separate off a crossing $t$ from the group of size $b$, denoting
the result as a group of size $b'=b-1$ (note that $b'=0$ is possible in this
setting). Consider the triangle with nodes the crossing $t$, the group of size $c$, 
and the crossing $x$, as shown  in Figure \ref{induction for loner ots} (c).
OTS the arc with endpoints $t$ and $x$ across the group $C$, obtaining
the result shown in Figure \ref{more about induction for loner ots} (a).

\begin{figure}[ht]
 \centering
   \begin{tabular}{c@{\hskip20pt}c@{\hskip20pt}c}
%%%%%%%%%%%%%%%%%%%%%%%%%%%%%%%%%%%%%%%%%%%%%%%%%%%%%%%%%%%%%%%%%%
    $\vcenter{\xy /r6pt/:,
      (0,0)="b1",
      (.8660254,.5)="b2",
      "b1";"b2",{\ellipse(,.5){}},
      (13,5)="v2",
      (10,-1)="v3",
      (6,2.5)="c1",
      "c1"+(.8660254,.5)="c2",
      "c1";"c2",{\ellipse(,.5){}},
      (4,8)="a1",
      (4,7)="a2",
      "b1"+(-.125,.21650635)="b1u",
      "b1"+(.125,-.21650635)="b1d",
      "b1u"+(-.8660254,-.5)+(-.125,.21650635)="b1ul",
      "b1d"+(-.8660254,-.5)+(.125,-.21650635)="b1dl",
      "b1u"+(.08660254,.05)="b1ur",
      "b1d"+(.08660254,.05)="b1dr",      
      "b1dl";"b1dr"**\dir{-},
      "b1ul";"b1ur"**\dir{-},
      "b2"+"b2"-"b1"+(-.125,.21650635)="b2u",
      "b2"+"b2"-"b1"+(.125,-.21650635)="b2d",
      "b2u"+(.8660254,.5)="b2ur",
      "b2d"+(.8660254,.5)="b2dr",
      "b2u"+(-.08660254,-.05)="b2ul",
      "b2d"+(-.08660254,-.05)="b2dl",
      "b2dl";"b2dr"**\dir{}?(.4)="b2dm",
      "b2ul";"b2ur"**\dir{}?(.4)="b2um",
      "b2ul";"b2um"**\dir{-},      
%% done b      
      "c1"+(-.125,.21650635)="c1u",
      "c1"+(.125,-.21650635)="c1d",
      "c1u"+(-.8660254,-.5)+(-.125,.21650635)="c1ul",
      "c1d"+(-.8660254,-.5)+(.125,-.21650635)="c1dl",
      "c1u"+(.08660254,.05)="c1ur",
      "c1d"+(.08660254,.05)="c1dr",      
      "c2"+"c2"-"c1"+(-.125,.21650635)="c2u",
      "c2"+"c2"-"c1"+(.125,-.21650635)="c2d",
      "c2u"+(.8660254,.5)="c2ur",
      "c2d"+(.8660254,.5)="c2dr",
      "c2u"+(-.08660254,-.05)="c2ul",
      "c2d"+(-.08660254,-.05)="c2dl",
% done c
      "a1"+(.25,0)="a1u",
      "a1"+(-.25,0)="a1d",
      "a1u"+(0,1)+(.25,0)="a1ul",
      "a1d"+(0,1)+(-.25,0)="a1dl",
      "a1u"+(0,-.1)="a1ur",
      "a1d"+(0,-.1)="a1dr",      
      "a1dl";"a1dr"**\dir{},
      "a1ul";"a1ur"**\dir{},
      "a2"+"a2"-"a1"+(.25,0)="a2u",
      "a2"+"a2"-"a1"+(-.25,0)="a2d",
      "a2u"+(0,.1)="a2ul",
      "a2d"+(0,.1)="a2dl",
      "a2u"+(0,-1)="a2ur",
      "a2d"+(0,-1)="a2dr",
      "a2dl";"a2dr"**\dir{}?(.4)="a2dm",
      "a2ul";"a2ur"**\dir{}?(.4)="a2um",
      "a2dl";"a2dm"**\dir{},
      "a2ul";"a2um"**\dir{},
% done a      
      "b1";"a1"**\dir{}?(.5)="ba",
      "c1";"a1"**\dir{}?(.5)="ca",
      "b1";"c1"**\dir{}?(.5)="bc",      
      "b2ur";"b2ur"+"b2ur"-"b2ul"+"b2ur"-"b2ul"+"b2ur"-"b2ul"**\dir{}?(.6)="xx1",
      "b2u";"a1ul"**\crv{"xx1" & 
         "c1"+(-2,5) & "c1" + (7.5,-.5) & "c2"+(3.5,3.5)  & "c2"+(-2.5,3)},
      "b2dl";"c1ur"**\crv{"b2dl"+(1,0) &  "c1ul"+(-1,0)},
      "c1dr";"v3"+(-.2,0)**\crv{"c1dr"+(-1,-1) },             
      "c2ul";"a1dl"+(-.2,0)**\crv{"c2ul"+(1,1) & "a1dl"+(0,-5) },             
      "c2dl";"v2"**\crv{"c2dl"+(1.5,.2) & "v2"+(-1,.2)},
      "b2"*{\hbox{\scriptsize$b'$}},
      "c2"*{\hbox{\scriptsize$c$}},
      "a1ul"*!<-4pt,0pt>{\hbox{\scriptsize$1$}},
      "a1dl"*!<4pt,0pt>{\hbox{\scriptsize$6$}},
      "v2"*!<-4pt,-3pt>{\hbox{\scriptsize$2$}},
      "v3"*!<1pt,5pt>{\hbox{\scriptsize$3$}},
      "b1dl"*!<-3pt,4pt>{\hbox{\scriptsize$4$}},
      "b1ul"*!<1pt,-5pt>{\hbox{\scriptsize$5$}},            
      "c1"+(-2,2.8)*{\hbox{\scriptsize$x$}},
      "c2"+(3.1,.4)*{\hbox{\scriptsize$t$}},
      "c2"+(5.5,1)*{\hbox{\scriptsize$y$}},
     \endxy}$
&
  %%%%%%%%%%%%%%%%%%%%%%%%%%%%%%%%%%%%%%%%%%%%%%%%%%%%%%%%%%
    $\vcenter{\xy /r6pt/:,
      (0,0)="b1",
      (.8660254,.5)="b2",
      (13,5)="v2",
      (6,-1)="v3",
      (6,2.5)="c1",
      "c1"+(.8660254,.5)="c2",
      "c1";"c2",{\ellipse(,.5){}},
      (4,8)="a1",
      (4,7)="a2",
      "b1"+(-.125,.21650635)="b1u",
      "b1"+(.125,-.21650635)="b1d",
      "b1u"+(-.8660254,-.5)+(-.125,.21650635)="b1ul",
      "b1d"+(-.8660254,-.5)+(.125,-.21650635)="b1dl",
      "b1u"+(.08660254,.05)="b1ur",
      "b1d"+(.08660254,.05)="b1dr",      
%      "b1dl";"b1dr"**\dir{-},
%      "b1ul";"b1ur"**\dir{-},
      "b2"+"b2"-"b1"+(-.125,.21650635)="b2u",
      "b2"+"b2"-"b1"+(.125,-.21650635)="b2d",
      "b2u"+(.8660254,.5)="b2ur",
      "b2d"+(.8660254,.5)="b2dr",
      "b2u"+(-.08660254,-.05)="b2ul",
      "b2d"+(-.08660254,-.05)="b2dl",
      "b2dl";"b2dr"**\dir{}?(.4)="b2dm",
      "b2ul";"b2ur"**\dir{}?(.4)="b2um",
%     "b2ul";"b2um"**\dir{-},      
%% done b      
      "c1"+(-.125,.21650635)="c1u",
      "c1"+(.125,-.21650635)="c1d",
      "c1u"+(-.8660254,-.5)+(-.125,.21650635)="c1ul",
      "c1d"+(-.8660254,-.5)+(.125,-.21650635)="c1dl",
      "c1u"+(.08660254,.05)="c1ur",
      "c1d"+(.08660254,.05)="c1dr",      
      "c2"+"c2"-"c1"+(-.125,.21650635)="c2u",
      "c2"+"c2"-"c1"+(.125,-.21650635)="c2d",
      "c2u"+(.8660254,.5)="c2ur",
      "c2d"+(.8660254,.5)="c2dr",
      "c2u"+(-.08660254,-.05)="c2ul",
      "c2d"+(-.08660254,-.05)="c2dl",
% done c
      "a1"+(.25,0)="a1u",
      "a1"+(-.25,0)="a1d",
      "a1u"+(0,1)+(.25,0)="a1ul",
      "a1d"+(0,1)+(-.25,0)="a1dl",
      "a1u"+(0,-.1)="a1ur",
      "a1d"+(0,-.1)="a1dr",      
      "a1dl";"a1dr"**\dir{},
      "a1ul";"a1ur"**\dir{},
      "a2"+"a2"-"a1"+(.25,0)="a2u",
      "a2"+"a2"-"a1"+(-.25,0)="a2d",
      "a2u"+(0,.1)="a2ul",
      "a2d"+(0,.1)="a2dl",
      "a2u"+(0,-1)="a2ur",
      "a2d"+(0,-1)="a2dr",
      "a2dl";"a2dr"**\dir{}?(.4)="a2dm",
      "a2ul";"a2ur"**\dir{}?(.4)="a2um",
      "a2dl";"a2dm"**\dir{},
      "a2ul";"a2um"**\dir{},
% done a      
      "b1";"a1"**\dir{}?(.5)="ba",
      "c1";"a1"**\dir{}?(.5)="ca",
      "b1";"c1"**\dir{}?(.5)="bc",      
      "c2u";"a1ul"**\crv{"c2u"+(1,1) &
"c1" + (7,0) & "c2"+(3.5,3.5)  & "c2"+(-2.5,3)},
      "c1dr";"b1dl"+(2,0)**\crv{"c1dr"+(-1.5,-.5) },             
      "c2dl";"v2"**\crv{"c2dl"+(1.5,.2) & "v2"+(-1,.2)},
      "c2"*{\hbox{\scriptsize$b'$}},
      "a1ul"*!<-4pt,0pt>{\hbox{\scriptsize$1$}},
      "a1dl"*!<4pt,0pt>{\hbox{\scriptsize$6$}},
      "v2"*!<-4pt,-3pt>{\hbox{\scriptsize$2$}},
      "v3"*!<-5pt,5pt>{\hbox{\scriptsize$3$}},
      "b1dl"*!<-7pt,2pt>{\hbox{\scriptsize$4$}},
      "b1ul"*!<1pt,-5pt>{\hbox{\scriptsize$5$}},            
      "c1"+(-2,-2.5)*{\hbox{\scriptsize$x$}},
      "c2"+(3.1,.4)*{\hbox{\scriptsize$t$}},
      "c2"+(5.5,1)*{\hbox{\scriptsize$y$}},
      %%%%%%%%%%%%%%%%%%%%%%%
      (4,1.5)="r1",
      "r1"-(.8660254,-.5)="r2",
      "r1"+(-2.1,1.2)*{\hbox{\scriptsize$c$}},
      "r1"+(.125,.21650635)="r1u",
      "r1"+(-.125,-.21650635)="r1d",
      "r1u"+(.8660254,-.5)+(.125,.21650635)="r1ura",
      "r1d"+(.8660254,-.5)+(-.125,-.21650635)="r1dra",
%%%%  move C left and up one and a half units
      "r1"+(-.8660254,.5)+(-.4330127,.25)="r1",
      "r1"-(.8660254,-.5)="r2",
      "r1";"r2",{\ellipse(,.5){}},
      "r1"+(.125,.21650635)="r1u",
      "r1"+(-.125,-.21650635)="r1d",
      "r1u"+(.8660254,-.5)+(.125,.21650635)="r1ur",
      "r1d"+(.8660254,-.5)+(-.125,-.21650635)="r1dr",
      "r1u"+(-.08660254,.05)="r1ul",
      "r1d"+(-.08660254,.05)="r1dl",      
      "v3";"r1dl"**\crv{"r1dr"+(1,-.15)},
      "c1u";"r1ul"**\crv{"r1ur"+<-2pt,-2pt>},
      "r2"+"r2"-"r1"+(.125,.21650635)="r2u",
      "r2"+"r2"-"r1"+(-.125,-.21650635)="r2d",
      "r2u"+(-.8660254,.5)="r2ul",
      "r2d"+(-.8660254,.5)="r2dl",
      "r2u"+(.08660254,-.05)="r2ur",
      "r2d"+(.08660254,-.05)="r2dr",
      "r2dr";"r2dl"**\dir{}?(.4)="r2dm",
      "r2ur";"r2ul"**\dir{}?(.4)="r2um",
%      "r2dr";"r2dm"**\dir{-},
%      "r2ur";"r2um"**\dir{-},
      "r2u";"a1dl"+(-.2,0)**\crv{"r2ul"+(-1,1) & "a1dl"+(0,-5) },
      "r2d";"b1ul"+(.5,.25)**\crv{"r2d"+(-1,1) & "b1ul"+(1,.5)},
     \endxy}$
&
%%%%%%%%%%%%%%%%%%%%%%%%%%%%%%%%%%%%%%%%%%%%%%%%%%%%%%%%%%%%%%%%%%
    $\vcenter{\xy /r6pt/:,
      (4,-1)="b1",
      "b1"+(0,1)="b2",
      "b1";"b2",{\ellipse(,.5){}},
      "b1"+(-.25,0)="b1u",
      "b1"+(.25,0)="b1d",
      "b1u"+(0,1)+(-.25,0)="b1ul",
      "b1d"+(0,1)+(.25,0)="b1dl",
      "b1u"+(0,.1)="b1ur",
      "b1d"+(0,.1)="b1dr",      
      "b2"+(-.25,0)="b2u",
      "b2"+(.25,0)="b2d",
      "b2u"+(0,-.1)="b2ul",
      "b2d"+(0,-.1)="b2dl",
      "b2u"+(0,1)="b2ur",
      "b2d"+(0,1)="b2dr",
      (13,5)="v2",
      (7,-1)="v3",
      (1,-1)="v4",
      (-.5,1)="v5",
      (6,2.5)="c1",
      "c1"+(.8660254,.5)="c2",
      "c1";"c2",{\ellipse(,.5){}},
      (4,8)="a1",
      (4,7)="a2",
%% done b      
      "c1"+(-.125,.21650635)="c1u",
      "c1"+(.125,-.21650635)="c1d",
      "c1u"+(-.8660254,-.5)+(-.125,.21650635)="c1ul",
      "c1d"+(-.8660254,-.5)+(.125,-.21650635)="c1dl",
      "c1u"+(.08660254,.05)="c1ur",
      "c1d"+(.08660254,.05)="c1dr",      
      "c2"+"c2"-"c1"+(-.125,.21650635)="c2u",
      "c2"+"c2"-"c1"+(.125,-.21650635)="c2d",
      "c2u"+(.8660254,.5)="c2ur",
      "c2d"+(.8660254,.5)="c2dr",
      "c2u"+(-.08660254,-.05)="c2ul",
      "c2d"+(-.08660254,-.05)="c2dl",
% done c
      "a1"+(.25,0)="a1u",
      "a1"+(-.25,0)="a1d",
      "a1u"+(0,1)+(.25,0)="a1ul",
      "a1d"+(0,1)+(-.25,0)="a1dl",
      "a1u"+(0,-.1)="a1ur",
      "a1d"+(0,-.1)="a1dr",      
      "a1dl";"a1dr"**\dir{},
      "a1ul";"a1ur"**\dir{},
      "a2"+"a2"-"a1"+(.25,0)="a2u",
      "a2"+"a2"-"a1"+(-.25,0)="a2d",
      "a2u"+(0,.1)="a2ul",
      "a2d"+(0,.1)="a2dl",
      "a2u"+(0,-1)="a2ur",
      "a2d"+(0,-1)="a2dr",
      "a2dl";"a2dr"**\dir{}?(.4)="a2dm",
      "a2ul";"a2ur"**\dir{}?(.4)="a2um",
      "a2dl";"a2dm"**\dir{},
      "a2ul";"a2um"**\dir{},
% done a      
      "b1";"a1"**\dir{}?(.5)="ba",
      "c1";"a1"**\dir{}?(.5)="ca",
      "b1";"c1"**\dir{}?(.5)="bc",      
      "c2ul";"a1ul"**\crv{"c2u"+(.5,.2) & 
         "c2"+(1,1.7) & "c1" + (7,0) & "c2"+(3.5,3.5)  & "c2"+(-2.5,3)},
      "v5"+(.3,.7);"c1ur"**\crv{"c1ur"+(-1.6,.3)},
      "c1dr";"b2dr"**\crv{"c1dr"+(-1,-.5) },             
      "b1ur";"v4"+(.2,0)**\crv{"b1ur"+(0,-.5) & "v4"+(.25,0) },             
      "b2ur";"a1dl"+(-.2,0)**\crv{"b2ur"+(0,.5) & "b2ur"+(-1,2) & "a1dl"+(0,-5) },             
      "c2dl";"v2"**\crv{"c2dl"+(1.5,.2) & "v2"+(-1,.2)},
      "b1dr";"v3"**\crv{"b1dr"+(0,-.5) & "v3"+(-.25,0)},
      "b2"*{\hbox{\scriptsize$c$}},
      "c2"*{\hbox{\scriptsize$b'$}},
      "a1ul"*!<-4pt,0pt>{\hbox{\scriptsize$1$}},
      "a1dl"*!<4pt,0pt>{\hbox{\scriptsize$6$}},
      "v2"*!<-4pt,-3pt>{\hbox{\scriptsize$2$}},
      "v3"*!<1pt,5pt>{\hbox{\scriptsize$3$}},
      "v4"*!<3pt,4pt>{\hbox{\scriptsize$4$}},
      "v5"*!<1pt,-5pt>{\hbox{\scriptsize$5$}},            
      "b2"+(-1.5,3)*{\hbox{\scriptsize$x$}},
      "c2"+(2.9,.3)*{\hbox{\scriptsize$t$}},
      "c2"+(5.5,1)*{\hbox{\scriptsize$y$}},
     \endxy}$
     \\
     \noalign{\vskip 6pt}
     (a)   & (b) $b+1$ odd  & (c) $b+1$ even
   \end{tabular}
  \caption{}
  \label{more about induction for loner ots}
\end{figure}

\noindent After turning the group of size $c$ in Figure \ref{more about induction 
for loner ots} (a), we may apply the induction hypothesis to the ots triangle
whose nodes are the crossing $x$ together with the groups of size $b'$ and $c$.
We consider the two possibilities: $b+1$ odd, and $b+1$ even. Since $b'=b-1$ has 
the same parity as $b+1$, the result when $b+1$ is odd is as shown in 
Figure \ref{more about induction for loner ots} (b), and the result when $b+1$
is even (even in the case $b'=0$) is as shown in Figure \ref{more about induction 
for loner ots} (c). In both cases, the group of size $b'=b-1$ combines with the
group of 2 consisting of the crossings $t$ and $y$ to reconstitute the group of
size $b$, as required. The result now follows by induction.
\end{proof}

\begin{theorem}\label{b=1}
 Given any reduced alternating diagram $D$ of a link $L$, and any
 groups $A$, $B$, and $C$ in $D$ which form an $ots$-triangle
 $O$ in the condensation $G$ of $D$, there exists a finite
 sequence of $T$ and $OTS$ operations on $D$ that results in a reduced
 alternating diagram $D'$ of a prime alternating link $L'$ whose
 condensation is the result of applying $ots$ in $G$ to the $ots$-triangle 
 $O$. 

 Specifically, let $a$, $b$, and $c$ denote the sizes of $A$, $B$, and $C$,
 respectively, which form the configuration in $D$ shown in Figure
 \ref{group ots}, where the groups have been labelled so that if
 at least two of the groups have odd size, then $B$ and $C$ are taken to
 be groups of odd size.

 \kern-1\baselineskip
\begin{figure}[ht]
 \centerline{\begin{tabular}{c}
    $\vcenter{\xy /r6pt/:,
      (.5,8.5)="v1",
      (6.75,-.3)="v2",
      (6.3,-1.6)="v3",
      (-6.3,-1.6)="v4",
      (-6.75,-.3)="v5",
      (-.5,8.5)="v6",
      "v1"*!<-4pt,0pt>{\hbox{\scriptsize$1$}},
      "v6"*!<4pt,0pt>{\hbox{\scriptsize$6$}},
      "v2"*!<-4pt,-3pt>{\hbox{\scriptsize$2$}},
      "v3"*!<-3pt,4pt>{\hbox{\scriptsize$3$}},
      "v4"*!<3pt,4pt>{\hbox{\scriptsize$4$}},
      "v5"*!<4pt,-3pt>{\hbox{\scriptsize$5$}},            
%%%
      (0,5)="a1",
      "a1"+(0,1)="a2",
      "a1";"a2",{\ellipse(,.5){}},
      "a1"+(.128944965,.094858827)="a1u",
      "a1"+(-.128944965,.094858827)="a1d",
      "a2"+(0,1)+(.128944965,-.094858827)="a2u",
      "a2"+(0,1)+(-.128944965,-.094858827)="a2d",
%% do c
      (3,.75)="c1",
      "c1"+(.8660254,-.5)="c2",
      "c1";"c2",{\ellipse(,.5){}},
      "c1"+(-.02410254,-.158253179)="c1d",
      "c1"+(.02410254,.158253179)="c1u",
      "c2"+(.8660254,-.5)+(.02410254,.158253179)="c2u",
      "c2"+(.8660254,-.5)+(-.02410254,-.158253179)="c2d",
%% now do b
      (-3,.75)="b1",
      "b1"+(-.8660254,-.5)="b2",
      "b1";"b2",{\ellipse(,.5){}},
      "b1"+(.02410254,-.158253179)="b1d",
      "b1"+(.02410254,.158253179)="b1u",
      "b2"+(-.8660254,-.5)+(.02410254,-.158253179)="b2d",
      "b2"+(-.8660254,-.5)+(-.02410254,.158253179)="b2u",
      "a1u";"c1u"**\crv{(.75,1.75)},
      "a1d";"b1u"**\crv{(-.75,1.75)},
      "b1d";"c1d"**\crv{(0,1.55)},
      "a2d";"v6"**\crv{"v6"+(.4,-.4)},
      "a2u";"v1"**\crv{"v1"+(-.4,-.4)},
      "c2u";"v2"**\crv{"v2"+(-.6,-.3)},
      "c2d";"v3"**\crv{"v3"+(-.6,.65)},
      "b2d";"v4"**\crv{"v4"+(.6,.65)},
      "b2u";"v5"**\crv{"v5"+(.6,-.3)},
      "a2"*=0{\hbox{\tiny$A$}},
      "b2"*=0{\hbox{\tiny$B$}},
      "c2"*=0{\hbox{\tiny$C$}},
     \endxy}$
   \end{tabular}}
  \caption{}
  \label{group ots}
\end{figure}

\kern-\baselineskip
 \begin{alphlist}
  \item If $b$ and $c$ are both odd, then there is a finite sequence
  of $T$ and $OTS$ operations that transforms the configuration in
  Figure \ref{group ots} to that shown in Figure \ref{a even/a odd} (a)
  if $a$ is even, or Figure \ref{a even/a odd} (b) if $a$ is odd.

\begin{figure}[ht]
 \centerline{\begin{tabular}{c@{\hskip30pt}c}
    $\vcenter{\xy /r6pt/:,0;(-1,0):,  %/r6pt/::,
      (.5,8.5)="v1",
      (6.75,-.3)="v2",
      (6.3,-1.6)="v3",
      (-6.3,-1.6)="v4",
      (-6.75,-.3)="v5",
      (-.5,8.5)="v6",
      "v1"*!<4pt,0pt>{\hbox{\scriptsize$4$}},
      "v2"*!<4pt,3pt>{\hbox{\scriptsize$5$}},
      "v3"*!<3pt,-4pt>{\hbox{\scriptsize$6$}},
      "v4"*!<-3pt,-4pt>{\hbox{\scriptsize$1$}},
      "v5"*!<-4pt,3pt>{\hbox{\scriptsize$2$}},
      "v6"*!<-4pt,0pt>{\hbox{\scriptsize$3$}},
%%%
      (0,5)="a1",
      "a1"+(0,1)="a2",
      "a1";"a2",{\ellipse(,.5){}},
      "a1"+(.128944965,.094858827)="a1u",
      "a1"+(-.128944965,.094858827)="a1d",
      "a2"+(0,1)+(.128944965,-.094858827)="a2u",
      "a2"+(0,1)+(-.128944965,-.094858827)="a2d",
%% do c
      (3,.75)="c1",
      "c1"+(.8660254,-.5)="c2",
      "c1";"c2",{\ellipse(,.5){}},
      "c1"+(-.02410254,-.158253179)="c1d",
      "c1"+(.02410254,.158253179)="c1u",
      "c2"+(.8660254,-.5)+(.02410254,.158253179)="c2u",
      "c2"+(.8660254,-.5)+(-.02410254,-.158253179)="c2d",
%% now do b
      (-3,.75)="b1",
      "b1"+(-.8660254,-.5)="b2",
      "b1";"b2",{\ellipse(,.5){}},
      "b1"+(.02410254,-.158253179)="b1d",
      "b1"+(.02410254,.158253179)="b1u",
      "b2"+(-.8660254,-.5)+(.02410254,-.158253179)="b2d",
      "b2"+(-.8660254,-.5)+(-.02410254,.158253179)="b2u",
      "a1u";"c1u"**\crv{(.75,1.75)},
      "a1d";"b1u"**\crv{(-.75,1.75)},
      "b1d";"c1d"**\crv{(0,1.55)},
      "a2d";"v6"**\crv{"v6"+(.4,-.4)},
      "a2u";"v1"**\crv{"v1"+(-.4,-.4)},
      "c2u";"v2"**\crv{"v2"+(-.6,-.3)},
      "c2d";"v3"**\crv{"v3"+(-.6,.65)},
      "b2d";"v4"**\crv{"v4"+(.6,.65)},
      "b2u";"v5"**\crv{"v5"+(.6,-.3)},
      "a2"*=0{\hbox{\tiny$A$}},
      "b2"*=0{\hbox{\tiny$C$}},
      "c2"*=0{\hbox{\tiny$B$}},
     \endxy}$
    &
    $\vcenter{\xy /r6pt/:,0;(-1,0):,  %/r6pt/::,
      (.5,8.5)="v1",
      (6.75,-.3)="v2",
      (6.3,-1.6)="v3",
      (-6.3,-1.6)="v4",
      (-6.75,-.3)="v5",
      (-.5,8.5)="v6",
      "v1"*!<4pt,0pt>{\hbox{\scriptsize$4$}},
      "v2"*!<4pt,3pt>{\hbox{\scriptsize$5$}},
      "v3"*!<3pt,-4pt>{\hbox{\scriptsize$6$}},
      "v4"*!<-3pt,-4pt>{\hbox{\scriptsize$1$}},
      "v5"*!<-4pt,3pt>{\hbox{\scriptsize$2$}},
      "v6"*!<-4pt,0pt>{\hbox{\scriptsize$3$}},
%%%
      (0,5)="a1",
      "a1"+(0,1)="a2",
      "a1";"a2",{\ellipse(,.5){}},
      "a1"+(.128944965,.094858827)="a1u",
      "a1"+(-.128944965,.094858827)="a1d",
      "a2"+(0,1)+(.128944965,-.094858827)="a2u",
      "a2"+(0,1)+(-.128944965,-.094858827)="a2d",
%% do c
      (3,.75)="c1",
      "c1"+(.8660254,-.5)="c2",
      "c1";"c2",{\ellipse(,.5){}},
      "c1"+(-.02410254,-.158253179)="c1d",
      "c1"+(.02410254,.158253179)="c1u",
      "c2"+(.8660254,-.5)+(.02410254,.158253179)="c2u",
      "c2"+(.8660254,-.5)+(-.02410254,-.158253179)="c2d",
%% now do b
      (-3,.75)="b1",
      "b1"+(-.8660254,-.5)="b2",
      "b1";"b2",{\ellipse(,.5){}},
      "b1"+(.02410254,-.158253179)="b1d",
      "b1"+(.02410254,.158253179)="b1u",
      "b2"+(-.8660254,-.5)+(.02410254,-.158253179)="b2d",
      "b2"+(-.8660254,-.5)+(-.02410254,.158253179)="b2u",
      "a1u";"c1u"**\crv{(.75,1.75)},
      "a1d";"b1u"**\crv{(-.75,1.75)},
      "b1d";"c1d"**\crv{(0,1.55)},
      "a2d";"v6"**\crv{"v6"+(.4,-.4)},
      "a2u";"v1"**\crv{"v1"+(-.4,-.4)},
      "c2u";"v2"**\crv{"v2"+(-.6,-.3)},
      "c2d";"v3"**\crv{"v3"+(-.6,.65)},
      "b2d";"v4"**\crv{"v4"+(.6,.65)},
      "b2u";"v5"**\crv{"v5"+(.6,-.3)},
      "a2"*=0{\hbox{\tiny$A$}},
      "b2"*=0{\hbox{\tiny$B$}},
      "c2"*=0{\hbox{\tiny$C$}},
     \endxy}$   \\
     \noalign{\vskip 6pt}
     (a) $a$ even   & (b) $a$ odd
   \end{tabular}}
  \caption{}
  \label{a even/a odd}
\end{figure}

  \item If $b$ is odd, while $a$ and $c$ are even, then there is a
  finite sequence of $T$ and $OTS$ operations that transforms the
  configuration in Figure \ref{group ots} to that shown in Figure
  \ref{a even/b odd/b even} (a), while if $b$ is even and $a$ is even, then
  there is a finite sequence of $T$ and $OTS$ operations that transforms
  the configuration in Figure \ref{group ots} to that shown in Figure
  \ref{a even/b odd/b even} (b).

\begin{figure}[ht]
 \centerline{\begin{tabular}{c@{\hskip30pt}c}
    $\vcenter{\xy /r6pt/:,0;(-1,0):,  %/r6pt/::,
      (.5,8.5)="v1",
      (6.75,-.3)="v2",
      (6.3,-1.6)="v3",
      (-6.3,-1.6)="v4",
      (-6.75,-.3)="v5",
      (-.5,8.5)="v6",
      "v1"*!<4pt,0pt>{\hbox{\scriptsize$4$}},
      "v2"*!<4pt,3pt>{\hbox{\scriptsize$5$}},
      "v3"*!<3pt,-4pt>{\hbox{\scriptsize$6$}},
      "v4"*!<-3pt,-4pt>{\hbox{\scriptsize$1$}},
      "v5"*!<-4pt,3pt>{\hbox{\scriptsize$2$}},
      "v6"*!<-4pt,0pt>{\hbox{\scriptsize$3$}},
%%%
      (0,5)="a1",
      "a1"+(0,1)="a2",
      "a1";"a2",{\ellipse(,.5){}},
      "a1"+(.128944965,.094858827)="a1u",
      "a1"+(-.128944965,.094858827)="a1d",
      "a2"+(0,1)+(.128944965,-.094858827)="a2u",
      "a2"+(0,1)+(-.128944965,-.094858827)="a2d",
%% do c
      (3,.75)="c1",
      "c1"+(.8660254,-.5)="c2",
%      "c1";"c2",{\ellipse(,.5){}},
      "c1"+(-.02410254,-.158253179)="c1d",
      "c1"+(.02410254,.158253179)="c1u",
      "c2"+(.8660254,-.5)+(.02410254,.158253179)="c2u",
      "c2"+(.8660254,-.5)+(-.02410254,-.158253179)="c2d",
%% now do b
      (-3,.75)="b1",
      "b1"+(-.8660254,-.5)="b2",
      "b1";"b2",{\ellipse(,.5){}},
      "b1"+(.02410254,-.158253179)="b1d",
      "b1"+(.02410254,.158253179)="b1u",
      "b2"+(-.8660254,-.5)+(.02410254,-.158253179)="b2d",
      "b2"+(-.8660254,-.5)+(-.02410254,.158253179)="b2u",
      "a1u";"c1u"**\crv{(.75,1.75)},
      "a1d";"b1u"**\crv{(-.75,1.75)},
      "b1d";"c1d"**\crv{(0,1.55)},
      "a2d";"v6"**\crv{"v6"+(.4,-.4)},
      "a2u";"v1"**\crv{"v1"+(-.4,-.4)},
 %     "c2u";"v2"**\crv{"v2"+(-.6,-.3)},
 %     "c2d";"v3"**\crv{"v3"+(-.6,.65)},
      "c1d";"v2"**\crv{"v2"+(-.6,-.3)},
      "c1u";"v3"**\crv{"c1u"+(.5,-.3) & "v3"+(-.6,.65)},
      "b2d";"v4"**\crv{"v4"+(.6,.65)},
      "b2u";"v5"**\crv{"v5"+(.6,-.3)},
      "a2"*!<-20pt,0pt>{\hbox{\tiny$a+b-1$}},
      "b2"*=0{\hbox{\tiny$C$}},
%      "c2"*=0{\hbox{\tiny$B$}},
     \endxy}$
    &
    $\vcenter{\xy /r6pt/:,0;(-1,0):,  %/r6pt/::,
      (.5,8.5)="v1",
      (6.75,-.3)="v2",
      (6.3,-1.6)="v3",
      (-6.3,-1.6)="v4",
      (-6.75,-.3)="v5",
      (-.5,8.5)="v6",
      "v1"*!<4pt,0pt>{\hbox{\scriptsize$4$}},
      "v2"*!<4pt,3pt>{\hbox{\scriptsize$5$}},
      "v3"*!<3pt,-4pt>{\hbox{\scriptsize$6$}},
      "v4"*!<-3pt,-4pt>{\hbox{\scriptsize$1$}},
      "v5"*!<-4pt,3pt>{\hbox{\scriptsize$2$}},
      "v6"*!<-4pt,0pt>{\hbox{\scriptsize$3$}},
%%%
      (0,5)="a1",
      "a1"+(0,1)="a2",
      "a1";"a2",{\ellipse(,.5){}},
      "a1"+(.128944965,.094858827)="a1u",
      "a1"+(-.128944965,.094858827)="a1d",
      "a2"+(0,1)+(.128944965,-.094858827)="a2u",
      "a2"+(0,1)+(-.128944965,-.094858827)="a2d",
%% do c
      (3,.75)="c1",
      "c1"+(.8660254,-.5)="c2",
      "c1";"c2",{\ellipse(,.5){}},
      "c1"+(-.02410254,-.158253179)="c1d",
      "c1"+(.02410254,.158253179)="c1u",
      "c2"+(.8660254,-.5)+(.02410254,.158253179)="c2u",
      "c2"+(.8660254,-.5)+(-.02410254,-.158253179)="c2d",
%% now do b
      (-3,.75)="b1",
      "b1"+(-.8660254,-.5)="b2",
 %     "b1";"b2",{\ellipse(,.5){}},
      "b1"+(.02410254,-.158253179)="b1d",
      "b1"+(.02410254,.158253179)="b1u",
      "b2"+(-.8660254,-.5)+(.02410254,-.158253179)="b2d",
      "b2"+(-.8660254,-.5)+(-.02410254,.158253179)="b2u",
      "a1u";"c1u"**\crv{(.75,1.75)},
      "a1d";"b1u"**\crv{(-.75,1.75)},
      "b1d";"c1d"**\crv{(0,1.55)},
      "a2d";"v6"**\crv{"v6"+(.4,-.4)},
      "a2u";"v1"**\crv{"v1"+(-.4,-.4)},
      "c2u";"v2"**\crv{"v2"+(-.6,-.3)},
      "c2d";"v3"**\crv{"v3"+(-.6,.65)},
%      "b2d";"v4"**\crv{"v4"+(.6,.65)},
%      "b2u";"v5"**\crv{"v5"+(.6,-.3)},
      "b1u";"v4"**\crv{"b1u"+(-.5,-.3) & "v4"+(.6,.65)},
      "b1d";"v5"**\crv{"v5"+(.6,-.3)},
%      "c1d";"v2"**\crv{"v2"+(-.6,-.3)},
%      "c1u";"v3"**\crv{"c1u"+(.5,-.3) & "v3"+(-.6,.65)},
%      "a2"*=0{\hbox{\tiny$A$}},
       "a2"*!<-20pt,0pt>{\hbox{\tiny$a+c-1$}},
%      "b2"*=0{\hbox{\tiny$B$}},
      "c2"*=0{\hbox{\tiny$B$}},
     \endxy}$   \\
     \noalign{\vskip 6pt}
     (a) $a$ even, $b$ odd, $c$ even  & (b) $a$ even, $b$ even
   \end{tabular}}
  \caption{}
  \label{a even/b odd/b even}
\end{figure}

 \end{alphlist}
\end{theorem}

\begin{proof}
 The proof of (a) is a straightforward induction argument (utilizing
 Proposition \ref{ssc}) on the size of the group $A$, where $B$ and
 $C$ are both odd size groups.

 For (b), we begin with the case when $b$ is odd, while $a$ and $c$ are
 both even. The proof will be by induction on $n$, where $a=2n$. Our
 hypothesis is that for any odd group $B$ of size $b$ and even group $C$
 of size $c$ and $A$ a group of even size $a=2n$, there exists a finite
 sequence of $T$ and/or $OTS$ operations that will transform the diagram
 shown in Figure \ref{group ots} into the diagram shown in Figure
 \ref{a even/b odd/b even} (a). We shall provide an argument which will
 establish both the base case and the inductive step. Suppose now that we
 have such a situation with $n\ge1$. Apply Proposition \ref{group ots}
 to the diagram formed by $B$, $C$ and the bottom crossing of $A$. The
 result is as shown in Figure \ref{proof of a even/b odd/c even} (a).

\begin{figure}[ht]
 \centerline{\begin{tabular}{c@{\hskip30pt}c}
    $\vcenter{\xy /r6pt/:,
      (.5,8.5)="v1",
      (6.75,-.3)="v2",
      (2.3,-4.6)="v3",
      (-2.3,-4.6)="v4",
      (-6.75,-.3)="v5",
      (-.5,8.5)="v6",
      "v1"*!<-4pt,0pt>{\hbox{\scriptsize$1$}},
      "v6"*!<4pt,0pt>{\hbox{\scriptsize$6$}},
      "v2"*!<-4pt,-3pt>{\hbox{\scriptsize$2$}},
      "v3"*!<-3pt,4pt>{\hbox{\scriptsize$3$}},
      "v4"*!<3pt,4pt>{\hbox{\scriptsize$4$}},
      "v5"*!<4pt,-3pt>{\hbox{\scriptsize$5$}},            
%%%
      (0,5)="a1",
      "a1"+(0,1)="a2",
      "a1";"a2",{\ellipse(,.5){}},
      "a1"+(.128944965,.094858827)="a1u",
      "a1"+(-.128944965,.094858827)="a1d",
      "a2"+(0,1)+(.128944965,-.094858827)="a2u",
      "a2"+(0,1)+(-.128944965,-.094858827)="a2d",
%% do c
      (3,.75)="c1",
      "c1"+(.8660254,-.5)="c2",
      "c1";"c2",{\ellipse(,.5){}},
      "c1"+(-.02410254,-.158253179)="c1d",
      "c1"+(.02410254,.158253179)="c1u",
      "c2"+(.8660254,-.5)+(.02410254,.158253179)="c2u",
      "c2"+(.8660254,-.5)+(-.02410254,-.158253179)="c2d",
%% now do b
      (-3,.75)="b1",
      "b1"+(-.8660254,-.5)="b2",
      "b1";"b2",{\ellipse(,.5){}},
      "b1"+(.02410254,-.158253179)="b1d",
      "b1"+(.02410254,.158253179)="b1u",
      "b2"+(-.8660254,-.5)+(.02410254,-.158253179)="b2d",
      "b2"+(-.8660254,-.5)+(-.02410254,.158253179)="b2u",
      "a1u";"c1u"**\crv{(.75,1.75)},
      "a1d";"b1u"**\crv{(-.75,1.75)},
      "b1d";"c1d"**\crv{(0,1.55)},
      "a2d";"v6"**\crv{"v6"+(.4,-.4)},
      "a2u";"v1"**\crv{"v1"+(-.4,-.4)},
      "c2u";"v2"**\crv{"v2"+(-.6,-.3)},
      "c2d";"v4"**\crv{"v2"+(2.75,-3) & "v2"+(-5,-1)},
      "b2d";"v3"**\crv{"v5"+(-2.75,-3) & "v5"+(5,-1)},
      "b2u";"v5"**\crv{"v5"+(.6,-.3)},
      "a2"*!<-16pt,0pt>{\hbox{\tiny $a-1$}},
      "b2"*=0{\hbox{\tiny$C$}},
      "c2"*=0{\hbox{\tiny$B$}},
     \endxy}$
    &
    $\vcenter{\xy /r6pt/:,0;(-1,0):,  %/r6pt/::,
      (.5,9.5)="v1",
      (6.75,-.3)="v2",
      (6.3,-1.6)="v3",
      (-6.3,-1.6)="v4",
      (-6.75,-.3)="v5",
      (-.5,9.5)="v6",
      "v1"*!<4pt,0pt>{\hbox{\scriptsize$4$}},
      "v2"*!<4pt,3pt>{\hbox{\scriptsize$5$}},
%      "v3"*!<3pt,-4pt>{\hbox{\scriptsize$6$}},
%      "v4"*!<-3pt,-4pt>{\hbox{\scriptsize$1$}},
      "v5"*!<-4pt,3pt>{\hbox{\scriptsize$2$}},
      "v6"*!<-4pt,0pt>{\hbox{\scriptsize$3$}},
%%%
      (0,5)="a1",
      "a1"+(0,1)="a2",
      "a1";"a2",{\ellipse(,.5){}},
      "a1"+(.128944965,.094858827)="a1u",
      "a1"+(-.128944965,.094858827)="a1d",
      "a2"+(0,1)+(.128944965,-.094858827)="a2u",
      "a2"+(0,1)+(-.128944965,-.094858827)="a2d",
%% a-2
      (0,-5)="a11",
      "a11"+(0,-1)="a12",
      "a11";"a12",{\ellipse(,.5){}},
      "a11"+(.128944965,.094858827)="a11u",
      "a11"+(-.128944965,.094858827)="a11d",
      "a12"+(0,-1)+(.128944965,.094858827)="a12u",
      "a12"+(0,-1)+(-.128944965,.094858827)="a12d",
      "a12"+(.75,-3)="v7"*{.},
      "a12"+(-.75,-3)="v8"*{.},
      "a12u";"v7"**\crv{"v7"+(-.5,1)},
      "a12d";"v8"**\crv{"v8"+(.5,1)},
      "v8"*!<-4pt,-2pt>{\hbox{\scriptsize$1$}},
      "v7"*!<4pt,-2pt>{\hbox{\scriptsize$6$}},
      "a12"*!<-17pt,0pt>{\hbox{\tiny$a-2$}},
%% do c
      (3,.75)="c1",
      "c1"+(.8660254,-.5)="c2",
%      "c1";"c2",{\ellipse(,.5){}},
      "c1"+(-.02410254,-.158253179)="c1d",
      "c1"+(.02410254,.158253179)="c1u",
      "c2"+(.8660254,-.5)+(.02410254,.158253179)="c2u",
      "c2"+(.8660254,-.5)+(-.02410254,-.158253179)="c2d",
%% now do b
      (-3,.75)="b1",
      "b1"+(-.8660254,-.5)="b2",
      "b1";"b2",{\ellipse(,.5){}},
      "b1"+(.02410254,-.158253179)="b1d",
      "b1"+(.02410254,.158253179)="b1u",
      "b2"+(-.8660254,-.5)+(.02410254,-.158253179)="b2d",
      "b2"+(-.8660254,-.5)+(-.02410254,.158253179)="b2u",
      "a1u";"c1u"**\crv{(.75,1.75)},
      "a1d";"b1u"**\crv{(-.75,1.75)},
      "b1d";"c1d"**\crv{(0,1.55)},
      "a2d";"v1"**\crv{"v6"+(-.2,-1)},
      "a2u";"v6"**\crv{"v1"+(.2,-1)},
 %     "c2u";"v2"**\crv{"v2"+(-.6,-.3)},
 %     "c2d";"v3"**\crv{"v3"+(-.6,.65)},
      "c1d";"v2"**\crv{"v2"+(-.6,-.3)},
%      "c1u";"v3"**\crv{"c1u"+(.5,-.3) & "v3"+(-.6,.65)},
      "c1u";"a11u"**\crv{"c1u"+(7,-3) & "a11u"+(0,3)},
%      "b2d";"v4"**\crv{"v4"+(.6,.65)},
      "b2d";"a11d"**\crv{"b2d"+(-3,-1.75) & "a11d"+(0,3)},      
      "b2u";"v5"**\crv{"v5"+(.6,-.3)},
      "a2"*=0{\hbox{\tiny$B$}},
      "b2"*=0{\hbox{\tiny$C$}},
%      "c2"*=0{\hbox{\tiny$B$}},
     \endxy}$   \\
     \noalign{\vskip 6pt}
     (a) \vtop{\baselineskip=10pt\leftskip=0pt\hsize=30ex\noindent the first crossing of $A$
     has been $OTS$'ed across $B$ and $C$}
     & (b) \vtop{\leftskip=0pt\hsize=30ex\noindent the second crossing of
     $A$ has now been $OTS$'ed}
   \end{tabular}}
  \caption{}
  \label{proof of a even/b odd/c even}
\end{figure}

 We now apply Proposition \ref{group ots} to the diagram formed by group
 $C$, group $B$ and the bottom crossing of the group of size $a-1$. Since
 $C$ is a group of even size, the result is as shown in Figure
 \ref{proof of a even/b odd/c even} (b). If $a-2=0$, we have the proof
 of the base case, while if $a-2>0$, then the induction hypothesis applies
 and the remaining $a-2$ crossings can $OTS$ to join the group of size
 $b+1$ formed by merging $B$ and the first $OTS$'ed crossing of $A$. The
 result follows now by induction.

\begin{figure}[ht]
 \centerline{\begin{tabular}{c@{\hskip30pt}c}
    $\vcenter{\xy /r6pt/:,0;(-1,0):, 
      (.5,9.5)="v1",
      (6.75,-.3)="v2",
      (6.3,-1.6)="v3",
      (-6.3,-1.6)="v4",
      (-6.75,-.3)="v5",
      (-.5,9.5)="v6",
      "v1"*!<4pt,0pt>{\hbox{\scriptsize$4$}},
      "v2"*!<4pt,3pt>{\hbox{\scriptsize$5$}},
%      "v3"*!<3pt,-4pt>{\hbox{\scriptsize$6$}},
%      "v4"*!<-3pt,-4pt>{\hbox{\scriptsize$1$}},
      "v5"*!<-4pt,3pt>{\hbox{\scriptsize$2$}},
      "v6"*!<-4pt,0pt>{\hbox{\scriptsize$3$}},
%%%
      (0,5)="a1",
      "a1"+(0,1)="a2",
      "a1";"a2",{\ellipse(,.5){}},
      "a1"+(.128944965,.094858827)="a1u",
      "a1"+(-.128944965,.094858827)="a1d",
      "a2"+(0,1)+(.128944965,-.094858827)="a2u",
      "a2"+(0,1)+(-.128944965,-.094858827)="a2d",
%% a-2
      (0,-5)="a11",
      "a11"+(0,-1)="a12",
      "a11";"a12",{\ellipse(,.5){}},
      "a11"+(.128944965,.094858827)="a11u",
      "a11"+(-.128944965,.094858827)="a11d",
      "a12"+(0,-1)+(.128944965,.094858827)="a12u",
      "a12"+(0,-1)+(-.128944965,.094858827)="a12d",
      "a12"+(.75,-3)="v7"*{.},
      "a12"+(-.75,-3)="v8"*{.},
      "a12u";"v7"**\crv{"v7"+(-.5,1)},
      "a12d";"v8"**\crv{"v8"+(.5,1)},
      "v8"*!<-4pt,-2pt>{\hbox{\scriptsize$1$}},
      "v7"*!<4pt,-2pt>{\hbox{\scriptsize$6$}},
      "a12"*!<-17pt,0pt>{\hbox{\tiny$a-1$}},
%% do c
      (3,.75)="c1",
      "c1"+(.8660254,-.5)="c2",
%      "c1";"c2",{\ellipse(,.5){}},
      "c1"+(-.02410254,-.158253179)="c1d",
      "c1"+(.02410254,.158253179)="c1u",
      "c2"+(.8660254,-.5)+(.02410254,.158253179)="c2u",
      "c2"+(.8660254,-.5)+(-.02410254,-.158253179)="c2d",
%% now do b
      (-3,.75)="b1",
      "b1"+(-.8660254,-.5)="b2",
      "b1";"b2",{\ellipse(,.5){}},
      "b1"+(.02410254,-.158253179)="b1d",
      "b1"+(.02410254,.158253179)="b1u",
      "b2"+(-.8660254,-.5)+(.02410254,-.158253179)="b2d",
      "b2"+(-.8660254,-.5)+(-.02410254,.158253179)="b2u",
      "a1u";"c1u"**\crv{(.75,1.75)},
      "a1d";"b1u"**\crv{(-.75,1.75)},
      "b1d";"c1d"**\crv{(0,1.55)},
      "a2u";"v1"**\crv{"v1"+(-.4,-1)},
      "a2d";"v6"**\crv{"v6"+(.4,-1)},
      "c1d";"v2"**\crv{"v2"+(-.6,-.3)},
%      "c1u";"v3"**\crv{"c1u"+(.5,-.3) & "v3"+(-.6,.65)},
      "c1u";"a11u"**\crv{"c1u"+(7,-3) & "a11u"+(0,3)},
%      "b2d";"v4"**\crv{"v4"+(.6,.65)},
      "b2d";"a11d"**\crv{"b2d"+(-3,-1.75) & "a11d"+(0,3)},      
      "b2u";"v5"**\crv{"v5"+(.6,-.3)},
      "a2"*=0{\hbox{\tiny$C$}},
      "b2"*=0{\hbox{\tiny$B$}},
     \endxy}$
   &
    $\vcenter{\xy  /r6pt/:,0;(0,-1)::,  %/r6pt/::,
      (.5,9.5)="v1",
      (6.75,-.3)="v2",
      (6.3,-1.6)="v3",
      (-6.3,-1.6)="v4",
      (-6.75,-.3)="v5",
      (-.5,9.5)="v6",
      "v1"*!<-4pt,0pt>{\hbox{\scriptsize$3$}},
      "v2"*!<-4pt,3pt>{\hbox{\scriptsize$2$}},
      "v5"*!<4pt,3pt>{\hbox{\scriptsize$5$}},
      "v6"*!<4pt,0pt>{\hbox{\scriptsize$4$}},
%%%
      (0,5)="a1",
      "a1"+(0,1)="a2",
      "a1";"a2",{\ellipse(,.5){}},
      "a1"+(.128944965,.094858827)="a1u",
      "a1"+(-.128944965,.094858827)="a1d",
      "a2"+(0,1)+(.128944965,-.094858827)="a2u",
      "a2"+(0,1)+(-.128944965,-.094858827)="a2d",
%% a-2
      (0,-5)="a11",
      "a11"+(0,-1)="a12",
      "a11";"a12",{\ellipse(,.5){}},
      "a11"+(.128944965,.094858827)="a11u",
      "a11"+(-.128944965,.094858827)="a11d",
      "a12"+(0,-1)+(.128944965,.094858827)="a12u",
      "a12"+(0,-1)+(-.128944965,.094858827)="a12d",
      "a12"+(.75,-3)="v7"*{.},
      "a12"+(-.75,-3)="v8"*{.},
      "a12u";"v7"**\crv{"v7"+(-.5,1)},
      "a12d";"v8"**\crv{"v8"+(.5,1)},
      "v8"*!<4pt,-2pt>{\hbox{\scriptsize$6$}},
      "v7"*!<-4pt,-2pt>{\hbox{\scriptsize$1$}},
      "a12"*!<-17pt,0pt>{\hbox{\tiny$a-2$}},
%% do c
      (3,.75)="c1",
      "c1"+(.8660254,-.5)="c2",
%      "c1";"c2",{\ellipse(,.5){}},
      "c1"+(-.02410254,-.158253179)="c1d",
      "c1"+(.02410254,.158253179)="c1u",
      "c2"+(.8660254,-.5)+(.02410254,.158253179)="c2u",
      "c2"+(.8660254,-.5)+(-.02410254,-.158253179)="c2d",
%% now do b
      (-3,.75)="b1",
      "b1"+(-.8660254,-.5)="b2",
      "b1";"b2",{\ellipse(,.5){}},
      "b1"+(.02410254,-.158253179)="b1d",
      "b1"+(.02410254,.158253179)="b1u",
      "b2"+(-.8660254,-.5)+(.02410254,-.158253179)="b2d",
      "b2"+(-.8660254,-.5)+(-.02410254,.158253179)="b2u",
      "a1d";"c1u"**\crv{"a1d"+(-.6,-1) &(0,1.75)},
      "a1u";"b1u"**\crv{"a1u"+(.6,-1) & (0,1.75)},
      "b1d";"c1d"**\crv{(0,1.55)},
      "a2u";"v1"**\crv{"v1"+(-.4,-1)},
      "a2d";"v6"**\crv{"v6"+(.4,-1)},
      "c1d";"v2"**\crv{"v2"+(-.6,-.3)},
%      "c1u";"v3"**\crv{"c1u"+(.5,-.3) & "v3"+(-.6,.65)},
      "c1u";"a11u"**\crv{"c1u"+(7,-3) & "a11u"+(0,3)},
%      "b2d";"v4"**\crv{"v4"+(.6,.65)},
      "b2d";"a11d"**\crv{"b2d"+(-3,-1.75) & "a11d"+(0,3)},      
      "b2u";"v5"**\crv{"v5"+(.6,-.3)},
      "a2"*=0{\hbox{\tiny$C$}},
      "b2"*=0{\hbox{\tiny$B$}},
     \endxy}$ \\
     \noalign{\vskip 6pt}
     (a) \vtop{\baselineskip=10pt\leftskip=0pt\hsize=30ex\noindent the
     first crossing of $A$ has been $OTS$'ed across $B$ and $C$}
     & (b) \vtop{\leftskip=0pt\hsize=30ex\noindent the second crossing of
     $A$ has now been $OTS$'ed}
   \end{tabular}}
  \caption{}
  \label{proof of a even/b even}
\end{figure}

 Finally, consider the case when $b$ is even and $a$ is even. The proof
 will be by induction on $n$, where $a=2n$. Our hypothesis is that for any
 even group $B$ of size $b$, group $C$ of size $c$, and group $A$ of even
 size $a=2n$, there exists a finite sequence of $T$ and/or $OTS$ operations
 that will transform the diagram shown in Figure \ref{group ots} into the
 diagram shown in Figure \ref{a even/b odd/b even} (b). We shall provide
 an argument which will establish both the base case and the inductive
 step. Suppose now that we have such a situation with $n\ge1$. Apply
 Proposition \ref{group ots} to the diagram formed by $B$, $C$ and the
 bottom crossing of $A$. The result is as shown in Figure
 \ref{proof of a even/b even} (a). We may now apply the induction hypothesis
 (or the base case has been established in the case $a=2$) to $OTS$ the
 group of size $a-2$ across the arc joining $B$ and the single crossing so
 as to join up with the group which is the combination of $C$ and the first
 $OTS$'ed crossing of $A$, thereby forming a group of size
 $(a-2)+1-1+(c+1)=a+c-1$, as required. The result now follows by induction.
\end{proof}

\par
 \section{A 9-crossing example}
 
  In this final section, we give an example of a 9-crossing alternating
  link (actually, knot $9_{32}$) and its reduction to the 9-crossing
  torus link. We begin with a reduced alternating diagram of knot $9_{32}$
  and its graph, and then show one sequence of $t$ and $ots$ that will
  reduce the graph to the one vertex, 2 loop graph alongside the
  changes in the knot diagram that result from the corresponding
  sequence of $T$ and $OTS$ applications. The end result is a reduced
  alternating diagram of the 9-crossing torus knot.
  
\begin{figure}[ht]
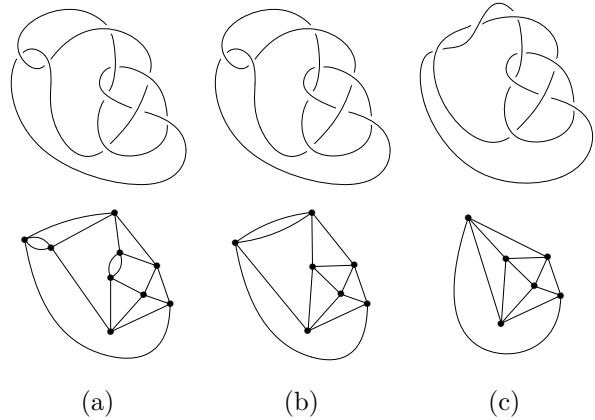

\centering
\begin{tabular}{c@{\hskip10pt}c@{\hskip10pt}c}
%% knot diagram no 1
$\vcenter{\xy /r10pt/:,
 (2,5.7)="1", 
 (1,6)="2", 
 (4.4,7)="3", 
 (4.6,5.5)="4", 
 (4.25,4.55)="5",
 (4.25,2.5)="6",
 (6.5,3.55)="7", 
 (6,5)="8",
 (5.5,3.9)="9", 
 "1"+(-.1,-.1);"3"+(-.15,.15)**\crv{"1"+(-.4,-.4) & "2"+(-.4,0)
  & "2"+(.18,1.55)  & "3"+(-1,1)},
% "2"*{o},
 "3"+(0.03,-.2);"5"+(.1,.1)**\crv{"3"+(.1,-.25) & "4"+(.1,0) & "5"+(.2,.2)},
% "4"*{o},     
 "5"+(-.1,-.275);"7"+(.05,-.2)**\crv{"5"+(-.4,-.7) & "6"+(-.6,-.9) & "7"+(.125,-1.25)},
% "6"*{o},     
 "7"+(.05,.15);"4"+(.15,.05)**\crv{"7"+(0,.75) & "8"+(.2,.2) & "4"+(.4,.1)},
% "8"*{o},   
 "4"+(-.2,.025);"9"+(-.3,.1)**\crv{"4"+(-.4,0) & "4"+(-.8,-.5) & "5"+(-.1,-.15) },
% "5"*{o},   
 "9"+(.2,-.05);"2"+(-.15,-.15)**\crv{"7"+(.35,0) & "7"+(1.15,-1.1) & "7"+(0,-3)
 & "2"+(0,-4) & "2"+(-.6,-1)},
% "7"*{o},   
 "2"+(.1,.1);"6"+(-.15,-.15)**\crv{"2"+(.35,.35) & "2"+(1.1,.1) & "2"+(1.2,-1)  
   &"6"+(-2.25,.9) & "6"+(-1,-.7)},
 "6"+(.15,.15);"8"+(-.2,-.15)**\crv{"6"+(.7,.7) & "9"+(-.05,-.15)},
 "8"+(-.1,.35);"1"+(.15,.15)**\crv{"8"+(.1,1) & "3"+(.6,.3) & "1"+(1,1.3)}, 
 \endxy}$
&
%% knot diagram no 2
%% initial contraction, no change in knot
$\vcenter{\xy /r10pt/:,
 (2,5.7)="1", 
 (1,6)="2", 
 (4.4,7)="3", 
 (4.6,5.5)="4", 
 (4.25,4.55)="5",
 (4.25,2.5)="6",
 (6.5,3.55)="7", 
 (6,5)="8",
 (5.5,3.9)="9", 
 "1"+(-.1,-.1);"3"+(-.15,.15)**\crv{"1"+(-.4,-.4) & "2"+(-.4,0)
  & "2"+(.18,1.55)  & "3"+(-1,1)},
% "2"*{o},
 "3"+(0.03,-.2);"5"+(.1,.1)**\crv{"3"+(.1,-.25) & "4"+(.1,0) & "5"+(.2,.2)},
% "4"*{o},     
 "5"+(-.1,-.275);"7"+(.05,-.2)**\crv{"5"+(-.4,-.7) & "6"+(-.6,-.9) & "7"+(.125,-1.25)},
% "6"*{o},     
 "7"+(.05,.15);"4"+(.15,.05)**\crv{"7"+(0,.75) & "8"+(.2,.2) & "4"+(.4,.1)},
% "8"*{o},   
 "4"+(-.2,.025);"9"+(-.3,.1)**\crv{"4"+(-.4,0) & "4"+(-.8,-.5) & "5"+(-.1,-.15) },
% "5"*{o},   
 "9"+(.2,-.05);"2"+(-.15,-.15)**\crv{"7"+(.35,0) & "7"+(1.15,-1.1) & "7"+(0,-3)
 & "2"+(0,-4) & "2"+(-.6,-1)},
% "7"*{o},   
 "2"+(.1,.1);"6"+(-.15,-.15)**\crv{"2"+(.35,.35) & "2"+(1.1,.1) & "2"+(1.2,-1)  
   &"6"+(-2.25,.9) & "6"+(-1,-.7)},
 "6"+(.15,.15);"8"+(-.2,-.15)**\crv{"6"+(.7,.7) & "9"+(-.05,-.15)},
 "8"+(-.1,.35);"1"+(.15,.15)**\crv{"8"+(.1,1) & "3"+(.6,.3) & "1"+(1,1.3)}, 
 \endxy}$
&
%% knot diagram no 3
%% second level of initial contraction
$\vcenter{\xy /r10pt/:,
% (2,5.7)="1",
(1.5,5.6)="1",
% (1,6)="2",
% "2"*{o},
 (2.75,6.5)="2", 
 (4.4,7)="3", 
 (4.6,5.5)="4", 
 (4.25,4.55)="5",
 (4.25,2.5)="6",
 (6.5,3.55)="7", 
 (6,5)="8",
 (5.5,3.9)="9", 
 "1"+(.1,.1);"3"+(-.15,.15)**\crv{"1"+(.4,.2) & "2"+(0,-1)
    & "3"+(-1,1)},
% "2"*{o},
 "3"+(0.03,-.2);"5"+(.1,.1)**\crv{"3"+(.1,-.25) & "4"+(.1,0) & "5"+(.2,.2)},
% "4"*{o},     
 "5"+(-.1,-.275);"7"+(.05,-.2)**\crv{"5"+(-.4,-.7) & "6"+(-.6,-.9) & "7"+(.125,-1.25)},
% "6"*{o},     
 "7"+(.05,.15);"4"+(.15,.05)**\crv{"7"+(0,.75) & "8"+(.2,.2) & "4"+(.4,.1)},
% "8"*{o},   
 "4"+(-.2,.025);"9"+(-.3,.1)**\crv{"4"+(-.4,0) & "4"+(-.8,-.5) & "5"+(-.1,-.15) },
% "5"*{o},   
 "9"+(.2,-.05);"1"+(-.3,0)**\crv{"7"+(.35,0) & "7"+(1.15,-1.1) & "7"+(0,-3)
 & "2"+(0,-6) & "1"+(-1,-2) & "1"+(-.8,-.1)},
% "7"*{o},   
 "2"+(-.1,.1);"6"+(-.15,-.15)**\crv{"2"+(-.15,0) & "2"+(-1.71,-.2)  
   &"6"+(-3,.9) & "6"+(-1,-.7)},
 "6"+(.15,.15);"8"+(-.2,-.15)**\crv{"6"+(.7,.7) & "9"+(-.05,-.15)},
 "8"+(-.1,.35);"2"+(.35,.15)**\crv{"8"+(.1,1) & "3"+(.6,.3)
    & "2"+(.5,.3)}, 
 \endxy}$\\
 \noalign{\vskip10pt}
%% graph no 1
$\vcenter{\xy /r10pt/:,
 (2,5.7)="1"*=0{\hbox{\tiny$\bullet$}}, 
 (1,6)="2"*=0{\hbox{\tiny$\bullet$}}, 
 (4.4,7)="3"*=0{\hbox{\tiny$\bullet$}}, 
 (4.6,5.5)="4"*=0{\hbox{\tiny$\bullet$}}, 
 (4.25,4.55)="5"*=0{\hbox{\tiny$\bullet$}},
 (4.25,2.5)="6"*=0{\hbox{\tiny$\bullet$}},
 (6.5,3.55)="7"*=0{\hbox{\tiny$\bullet$}}, 
 (6,5)="8"*=0{\hbox{\tiny$\bullet$}},
 (5.5,3.9)="9"*=0{\hbox{\tiny$\bullet$}},
 {\ar @{{}{-}{}} @/^.5ex/ "1";"2"},
 {\ar @{{}{-}{}} @/_.5ex/ "1";"2"},
 {\ar @{{}{-}{}} @/^.5ex/ "2";"3"},
 "1";"3"**\dir{-},
 "3";"4"**\dir{-},
 "3";"8"**\dir{-},
 {\ar @{{}{-}{}} @/^.5ex/ "4";"5"},
 {\ar @{{}{-}{}} @/_.5ex/ "4";"5"},
 "4";"8"**\dir{-},
 "8";"9"**\dir{-},
 "5";"9"**\dir{-},
 "8";"7"**\dir{-},
 "9";"7"**\dir{-},
 "6";"7"**\dir{-},
 "9";"6"**\dir{-},
 "5";"6"**\dir{-},
 "6";"1"**\dir{-},
 "2";"7"**\crv{"2"+(.5,-3) & "2"+(3,-5) & "7"+(0,-2)},
\endxy}$
&
%% graph no 2
$\vcenter{\xy /r10pt/:,
% (2,5.7)="1"*=0{\hbox{\tiny$\bullet$}},
 (1.5,5.85)="12"*=0{\hbox{\tiny$\bullet$}},
% (1,6)="2"*=0{\hbox{\tiny$\bullet$}}, 
 (4.4,7)="3"*=0{\hbox{\tiny$\bullet$}}, 
% (4.6,5.5)="4"*=0{\hbox{\tiny$\bullet$}},
 (4.425,4.95)="45"*=0{\hbox{\tiny$\bullet$}},  
% (4.25,4.55)="5"*=0{\hbox{\tiny$\bullet$}},
 (4.25,2.5)="6"*=0{\hbox{\tiny$\bullet$}},
 (6.5,3.55)="7"*=0{\hbox{\tiny$\bullet$}}, 
 (6,5)="8"*=0{\hbox{\tiny$\bullet$}},
 (5.5,3.9)="9"*=0{\hbox{\tiny$\bullet$}},
 {\ar @{{}{-}{}} @/^.5ex/ "12";"3"},
 {\ar @{{}{-}{}} @/_.5ex/ "12";"3"},
 "3";"45"**\dir{-},
 "3";"8"**\dir{-},
% {\ar @{{}{-}{}} @/^.5ex/ "4";"5"},
% {\ar @{{}{-}{}} @/_.5ex/ "4";"5"},
 "45";"8"**\dir{-},
 "8";"9"**\dir{-},
 "45";"9"**\dir{-},
 "8";"7"**\dir{-},
 "9";"7"**\dir{-},
 "6";"7"**\dir{-},
 "9";"6"**\dir{-},
 "45";"6"**\dir{-},
 "6";"12"**\dir{-},
 "12";"7"**\crv{"12"+(.5,-3) & "12"+(3,-5) & "7"+(0,-2)},
\endxy}$
&
%% graph no 3
$\vcenter{\xy /r10pt/:,
 (3,6.5)="123"*=0{\hbox{\tiny$\bullet$}},
 (4.425,4.95)="45"*=0{\hbox{\tiny$\bullet$}},  
 (4.25,2.5)="6"*=0{\hbox{\tiny$\bullet$}},
 (6.5,3.55)="7"*=0{\hbox{\tiny$\bullet$}}, 
 (6,5)="8"*=0{\hbox{\tiny$\bullet$}},
 (5.5,3.9)="9"*=0{\hbox{\tiny$\bullet$}},
 "123";"45"**\dir{-},
 "123";"8"**\dir{-},
 "45";"8"**\dir{-},
 "8";"9"**\dir{-},
 "45";"9"**\dir{-},
 "8";"7"**\dir{-},
 "9";"7"**\dir{-},
 "6";"7"**\dir{-},
 "9";"6"**\dir{-},
 "45";"6"**\dir{-},
 "6";"123"**\dir{-},
 "123";"7"**\crv{"123"+(-1,-3) & "123"+(0,-6) & "7"+(0,-2)},
\endxy}$\\
\noalign{\vskip10pt}
 (a) & (b) & (c)\\
\end{tabular}
\caption{The initial stage in the reduction of knot $9_{32}$ to the
  9-crossing torus knot.}
\label{stage 1 of 9_32 reduction}
\end{figure}

 In Figure \ref{stage 1 of 9_32 reduction} (a), we present a reduced alternating
 diagram for knot $9_{32}$, with the graph of the diagram shown below it. Then
 in Figure \ref{stage 1 of 9_32 reduction} (b), the first round of contraction
 in the graph is shown. Note that there is no change in the diagram.

\begin{figure}[ht]
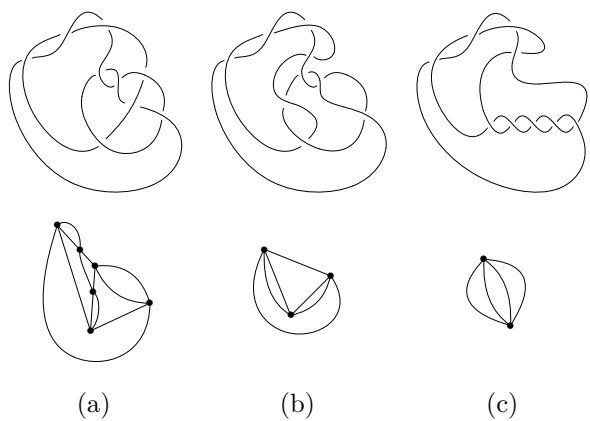

\centering
\begin{tabular}{c@{\hskip10pt}c@{\hskip10pt}c}
%% knot diagram no 4
%% diagram with first t to match ots graph
$\vcenter{\xy /r10pt/:,
% (2,5.7)="1",
(1.5,5.6)="1",
% (1,6)="2",
% "2"*{o},
 (2.75,6.5)="2", 
 (4.4,7)="3", 
 (4.6,5.5)="4", 
 (4.25,4.55)="5",
 (4.25,2.5)="6",
 (6.5,3.55)="7", 
 (6,5)="8",
 (5.5,3.9)="9",
 (4.8,4.75)="10",
 (4.2,5.25)="11",
 "1"+(.1,.1);"3"+(-.15,.15)**\crv{"1"+(.4,.2) & "2"+(0,-1)
    & "3"+(-1,1)},
% "2"*{o},
 "11"+(-.2,-.2);"7"+(0,-.2)**\crv{"11"+(-.4,-.2) & "11"+(-.9,-1) &"6"+(-.6,.5) &
    "7"+(-1.7,-1.8) & "7"+(-.15,-1.25)},
 "9"+(-.4,.1);"11"+(.2,0)**\crv{"9"+(-1,.25) & "9"+(-.2,1.6) },    
% "6"*{o},     
 "7"+(.05,.15);"10"+(.2,.2)**\crv{"7"+(.275,.95) & "8"+(-.3,.3)
 & "10"+(.4,.4) },
 "10"+(-.1,-.1);"3"+(.1,-.25)**\crv{"10"+(-.3,-.3)
      & "3"+(-.8,-1.7) & "3"+(.5,-.9)}, 
 "9"+(.2,-.05);"1"+(-.3,0)**\crv{"7"+(.35,0) & "7"+(1.15,-1.1) & "7"+(0,-3)
 & "2"+(0,-6) & "1"+(-1,-2) & "1"+(-.8,-.1)},
 "2"+(-.1,.1);"6"+(-.15,-.15)**\crv{"2"+(-.15,0) & "2"+(-1.71,-.2)  
   &"6"+(-3,.9) & "6"+(-1,-.7)},
 "6"+(.15,.15);"8"+(-.2,-.15)**\crv{"6"+(.7,.7) & "9"+(-.05,-.15)},
 "8"+(-.1,.35);"2"+(.35,.15)**\crv{"8"+(.1,1) & "3"+(.6,.3)
    & "2"+(.5,.3)}, 
 \endxy}$
 &
 %% diagram with both ots operations completed
$\vcenter{\xy /r10pt/:,
% (2,5.7)="1",
(1.5,5.6)="1",
% (1,6)="2",
% "2"*{o},
 (2.75,6.5)="2", 
 (4.4,7)="3", 
 (4.6,5.5)="4", 
 (4.25,4.55)="5",
 (4.25,2.5)="6",
 (6.5,3.55)="7", 
 (6,5)="8",
 (5.5,3.9)="9",
 (4.8,4.75)="10",
 (4.2,5.25)="11",
 (4.5,5.95)="12",
 (3.45,4.05)="13",
 "1"+(.1,.1);"3"+(-.15,.15)**\crv{"1"+(.4,.2) & "2"+(0,-1)
    & "3"+(-1,1)},
% "2"*{o},
 "11"+(-.2,-.2);"13"+(.1,.25)**\crv{"11"+(-.4,-.3) & "13"+(.05,.2) },
 "13"+(0,-.25);"7"+(0,-.2)**\crv{"13"+(-.2,-.8) &%"11"+(-.4,-.2) & "11"+(-.9,-1) &"6"+(-.6,.5) &
    "7"+(-1.7,-1.8) & "7"+(-.15,-1.25)},
% "9"+(-.4,.1);"11"+(.2,0)**\crv{"9"+(-1,.25) & "9"+(-.2,1.6) },    
% "6"*{o},     
 "7"+(.05,.15);"10"+(.2,.2)**\crv{"7"+(.275,.95) & "8"+(-.3,.3)
 & "10"+(.4,.4) },
 "10"+(-.1,-.1);"3"+(.1,-.25)**\crv{"10"+(-.3,-.3)
      & "3"+(-.8,-1.7) & "3"+(.5,-.9)},
 "11"+(.2,-.1);"1"+(-.3,0)**\crv{"11"+(.4,0) & "9"+(-.55,1) & "9"+(-.8,.25) &
    "9"+(.2,-.05) & "7"+(.35,0) & "7"+(1.15,-1.1) & "7"+(0,-3)
   & "2"+(0,-6) & "1"+(-1,-2) & "1"+(-.8,-.1)}, 
% "9"+(.2,-.05);"1"+(-.3,0)**\crv{"7"+(.35,0) & "7"+(1.15,-1.1) & "7"+(0,-3)
% & "2"+(0,-6) & "1"+(-1,-2) & "1"+(-.8,-.1)},
 "2"+(-.1,.1);"6"+(-.15,-.15)**\crv{"2"+(-.15,0) & "2"+(-1.71,-.2)  
   &"6"+(-3,.9) & "6"+(-1,-.7)},
% "6"+(.15,.15);"8"+(-.2,-.15)**\crv{"6"+(.7,.7) & "9"+(-.05,-.15)},
 "6"+(.15,.15);"12"+(-.2,-.15)**\crv{"6"+(.7,.7) & "13"+(.9,-.3)
   & "13"+(.5,-.1) & "13"+(-.5,.15) & "12"+(-1.4,-.5) &"12"+(-.5,0) },
 "12"+(.1,-.2);"2"+(.35,.15)**\crv{"12"+(.3,-.4) & "12"+(1.2,-.3) & "3"+(.6,.3)
    & "2"+(.8,.3)}, 
 \endxy}$
 &
 %% diagram showing next round of collapsing
$\vcenter{\xy /r10pt/:,
(1.5,5.6)="1",
 (2.75,6.5)="2", 
 (4.4,7)="3", 
 (3.25,2.75)="6",
% "6"*{o},
 (7,2.75)="7",
% "7"*{o},
 (4.5,5.95)="12",
 (3.25,3.75)="13",
% "13"*{o},
 (7,3.75)="14",
 % 3.25 to 6.8
 (3.55,3.25)="5-1",
 (4.35,3.25)="5-2",
 (5.15,3.25)="5-3",
 (5.95,3.25)="5-4",
 (6.75,3.25)="5-5",
 "13";"5-2"+(-.1,-.1)**\crv{ "5-2"+(-.6,-.8) },
 "5-1"+(.1,.1);"5-3"+(-.1,-.1)**\crv{"5-1"+(.6,.9) & "5-3"+(-.6,-.9)},
 "5-2"+(.1,.1);"5-4"+(-.1,-.1)**\crv{"5-2"+(.6,.9) & "5-4"+(-.6,-.9)},
 "5-3"+(.1,.1);"5-5"+(-.1,-.1)**\crv{"5-3"+(.6,.9) & "5-5"+(-.6,-.9)},
 "7";"5-4"+(.1,.1)**\crv{ "5-4"+(.6,.8) }, 
% "14"*{o},
 "1"+(.1,.1);"3"+(-.15,.15)**\crv{"1"+(.4,.2) & "2"+(0,-1)
    & "3"+(-1,1)},
% "2"*{o},
 "5-5"+(.1,.1);"3"+(.1,-.25)**\crv{ "5-5"+(1.1,2) &
    "5-5"+(-2.1,1.1)   & "3"+(-.3,-1.5) & "3"+(.3,-.6)},
 "7";"1"+(-.3,0)**\crv{"7"+(.3,-.8) & "7"+(-.5,-2.5)
       & "2"+(0,-5.5) & "1"+(-1,-2) & "1"+(-.8,-.1)}, 
 "2"+(-.1,.1);"5-1"+(-.15,-.15)**\crv{"2"+(-1.75,-.5)  
   &"6"+(-2,.9) & "5-1"+(-.7,-.9)},
 "13";"12"+(-.2,0)**\crv{"12"+(-1.75,-.5) &"12"+(-.7,.1) },
 "12"+(.2,0);"2"+(.35,.15)**\crv{"12"+(1.6,-.2)
    & "3"+(.6,.4) & "2"+(.5,.3)}, 
 \endxy}$
 \\
 \noalign{\vskip10pt}
%% graph no 4, showing ots
$\vcenter{\xy /r10pt/:,
 (3,6.5)="123"*=0{\hbox{\tiny$\bullet$}},
 (4.425,4.95)="45"*=0{\hbox{\tiny$\bullet$}},  
 (4.25,2.5)="6"*=0{\hbox{\tiny$\bullet$}},
 (6.5,3.55)="7"*=0{\hbox{\tiny$\bullet$}}, 
 (6,5)="8",%*=0{\hbox{\tiny$\bullet$}},
 (5.5,3.9)="9",%*=0{\hbox{\tiny$\bullet$}},
 "6";"45"**\dir{}?(.6)="645"*=0{\hbox{\tiny$\bullet$}},
 "123";"45"**\dir{}?(.6)="12345"*=0{\hbox{\tiny$\bullet$}},
 "123";"45"**\dir{-},
% "123";"8"**\dir{-},
% "45";"8"**\dir{-},
% "8";"9"**\dir{-},
% "45";"9"**\dir{-},
% "8";"7"**\dir{-},
% "9";"7"**\dir{-},
  "45";"7"**\crv{"45"+(.6,-1.4)},
  "45";"7"**\crv{"45"+(1.6,.2)},
 "6";"7"**\dir{-},
% "9";"6"**\dir{-},
 "6";"645"**\crv{"6"+(.4,.4) & "645"+(.3,-.2)},
 "645";"12345"**\crv{"645"+(-.1,.3) & "12345"+(-.1,-.3)},
 "123";"12345"**\crv{"123"+(.3,.3) & "12345"+(.1,.9)}, 
 "45";"6"**\dir{-},
 "6";"123"**\dir{-},
 "123";"7"**\crv{"123"+(-1,-3) & "123"+(0,-6) & "7"+(0,-2)},
\endxy}$
&
%% graph no 4, showing ots, first contractions
$\vcenter{\xy /r10pt/:,
 (0,0)="1"*=0{\hbox{\tiny$\bullet$}},
 (1.5,1.5)="2"*=0{\hbox{\tiny$\bullet$}},  
 (-1,2.5)="3"*=0{\hbox{\tiny$\bullet$}},
 "2";"3"**\dir{-},
 "1";"2"**\dir{-}?(.5)="x",
 "1";"2"**\crv{"x"+(.5,-.5) },
 "1";"3"**\dir{-}?(.5)="y",
 "1";"3"**\crv{"y"+(-.6,-.6) },
 "2";"3"**\crv{"2"+(1,-1.4) & "1"+(0,-1.75) & "3"+(-1,-1.9)},
\endxy}$
&
%% graph no 4, showing collapsing
 $\vcenter{\xy /r10pt/:,
   (0,0)="1"*=0{\hbox{\tiny$\bullet$}},
   (-1,2.5)="3"*=0{\hbox{\tiny$\bullet$}},
   "1";"3"**\dir{}?(.5)="y",
   "1";"3"**\crv{"y"+(-.6,-.5) },
   "1";"3"**\crv{"y"+(.6,.5) },
   "1";"3"**\crv{"y"+(-2.2,-.5)},
   "1";"3"**\crv{"y"+(2,1.1)},
\endxy}$\\
\noalign{\vskip10pt}
 (a) & (b) & (c) \\
%\noalign{\vskip8pt}
\end{tabular}
\caption{The second stage in the reduction of knot $9_{32}$ to the
  9-crossing torus knot.}
\label{stage 2 of 9_32 reduction}
\end{figure}

 Then in Figure \ref{stage 1 of 9_32 reduction} (c), we see that the second
 round of contraction results in a 4-regular simple graph. This required
 that a single $T$ operation be performed on a negative even (2) group in
 the diagram of Figure \ref{stage 1 of 9_32 reduction} (b), and, as
 described in Proposition \ref{all turn scenarios}, the number of components
 has gone up by one, so we now have a two component link. Since we
 have now arrived at a 4-regular simple graph, it is necessary to use
 $ots$ operations to empty out a min-2-region in the graph. There are eight
 $ots$-triangles in the graph shown in Figure \ref{stage 1 of 9_32 reduction}
 (c), and as it turns out, performing any one of the corresponding
 $ots$ operations will result in a 2-group. We choose an $ots$-triangle
 and perform an $ots$-operation on it. The result of our choice is shown
 in the graph of Figure \ref{stage 2 of 9_32 reduction} (a), while
 in the diagram above the graph, we show the result of the $T$ operation
 that must precede the two $OTS$ operations which, taken together,
 correspond to the $ots$ operation on the graph as shown in
 Figure \ref{stage 2 of 9_32 reduction} (a).

%% the second stage in the reduction of 9_32

%% the last stage in the reduction of 9_32
\begin{floatingfigure}
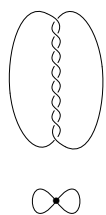
{.4\hsize}
%\begin{figure}[ht]
\centering
\begin{tabular}{c}
 %% diagram showing final round of collapsing
$\vcenter{\xy /r7pt/:,
(1.5,5.6)="1",
 (2.75,6.5)="2", 
 (4.4,7)="3", 
 (2.75,3.25)="bl",
 (3.75,3.1)="br",
 (2.75,10.3)="tl",
 (3.75,10.25)="tr",
 (3.25,3.55)="5-1",
 (3.25,4.35)="5-2",
 (3.25,5.15)="5-3",
 (3.25,5.95)="5-4",
 (3.25,6.75)="5-5",
 (3.25,7.55)="5-6",
 (3.25,8.35)="5-7",
 (3.25,9.15)="5-8",
 (3.25,9.95)="5-9",
 "br";"5-2"+(-.1,-.1)**\crv{ "5-2"+(-.45,-.75) },
 "5-1"+(.1,.1);"5-3"+(-.1,-.1)**\crv{"5-1"+(.6,.9) & "5-3"+(-.6,-.9)},
 "5-2"+(.1,.1);"5-4"+(-.1,-.1)**\crv{"5-2"+(.6,.9) & "5-4"+(-.6,-.9)},
 "5-3"+(.1,.1);"5-5"+(-.1,-.1)**\crv{"5-3"+(.6,.9) & "5-5"+(-.6,-.9)},
 "5-4"+(.1,.1);"5-6"+(-.1,-.1)**\crv{"5-4"+(.6,.9) & "5-6"+(-.6,-.9)},
 "5-5"+(.1,.1);"5-7"+(-.1,-.1)**\crv{"5-5"+(.6,.9) & "5-7"+(-.6,-.9)},
 "5-6"+(.1,.1);"5-8"+(-.1,-.1)**\crv{"5-6"+(.6,.9) & "5-8"+(-.6,-.9)},
 "5-7"+(.1,.1);"5-9"+(-.1,-.1)**\crv{"5-7"+(.6,.9) & "5-9"+(-.6,-.9)},
 "tl";"5-8"+(.1,.1)**\crv{ "5-8"+(.45,.75) },
 "tl";"5-1"+(-.1,-.1)**\crv{ "tl"+(-1.45,.75) & "tl"+(-3,-5) & "5-1"+(-1,-1.2) }, 
 "br";"5-9"+(.1,.1)**\crv{ "br"+(1.45,-.75) & "br"+(3,5) & "5-9"+(1,1.2) }, 
 \endxy}$
 \\
 \noalign{\vskip10pt}
%% graph no 4, showing collapsing
 $\vcenter{\xy /r10pt/:,
   (0,0)="1"*=0{\hbox{\tiny$\bullet$}},
   "1";"1"**\crv{"1"+(-1.2,-1.6) & "1"+(-1.2,1.6)},
   "1";"1"**\crv{"1"+(1.2,-1.6) & "1"+(1.2,1.6)},
\endxy}$\\
%\noalign{\vskip8pt}
\end{tabular}
\caption{The final stage in the reduction of knot $9_{32}$ to the
  9-crossing torus knot.}
\label{last stage of 9_32 reduction}
\end{floatingfigure}

 In the diagram shown in Figure \ref{stage 2 of 9_32 reduction} (b), the
 two $OTS$ operations have been performed, and the graph shown below the
 diagram displays the result of the first round of contractions in the
 graph of Figure \ref{stage 2 of 9_32 reduction} (a). No $T$ operations
 in the diagram are required for the first round of contractions.

 Next, the graph in Figure \ref{stage 2 of 9_32 reduction} (c) shows the
 result of the second round of contractions. This round of contractions
 require that two applications of $T$ be performed on the diagram
 shown in Figure \ref{stage 2 of 9_32 reduction} (b). One of the
 applications of $T$ is to a link 2-group and so the number of components
 is reduced by one, resulting in a diagram of a knot.

 In Figure \ref{last stage of 9_32 reduction}, the final collapse in
 the graph is shown. This requires one $T$ operation on the diagram, and the
 result is the 9-crossing torus knot.
 
\par
\section{References}

\noindent [1] J. A. Calvo, Knot enumeration through flypes
and twisted splices. J. Knot and its Ram., 6(1997), no. 6, 785--798.

\noindent [2] J. H. Conway, An enumeration of knots and links and
some of their related properties, Computational Problems in Abstract
Algebra (John Leech, ed.), Pergamon Press, Oxford and New York,
1969, 329--358.

\noindent [3] H. de Fraysseix and P. Ossona de Mendez, On a characterization
of Gauss codes, Discrete Comput. Geom. 22 (1999), no. 2, 267--295.

\noindent [4] W. B. R. Lickorish, Prime Knots and Tangles, Trans. Amer.
Math. Soc. 267(1981), no. 1, 321--332.

\noindent [5] Stuart Rankin, John Schermann, Ortho Smith,
 Enumerating the prime alternating knots, Part I,
 (to appear in J. Knot and its Ram.).

\noindent [6] Stuart Rankin, John Schermann, Ortho Smith,
 Enumerating the prime alternating knots, Part II,
 (to appear in J. Knot and its Ram.).

\end{document}